\documentclass[10pt,a4paper]{article}
\usepackage[english]{babel}
\usepackage{latexsym}
\usepackage{amssymb,amsbsy,amsmath,amsfonts,amssymb,amscd}
\usepackage{amsmath, amsthm}
\usepackage[utf8]{inputenc}
\usepackage{epsfig, graphicx}%,multirow}
\usepackage{color}
\usepackage{fullpage}
\usepackage{graphics}
\usepackage{wrapfig}
\usepackage{enumitem}
\usepackage{mathrsfs}
\usepackage{environ}
\usepackage{lipsum}
\usepackage{enumitem} 
\usepackage{sectsty}
\usepackage{tcolorbox}
\usepackage{dsfont}
\usepackage[colorlinks=true,linkcolor=red,citecolor=blue]{hyperref}
\usepackage{bm}
\usepackage{mathtools,accents}
\usepackage[marginparwidth=2cm]{geometry}

\def\X{\mathrm{X}}
\def\Y{\mathrm{Y}}
\def\V{\mathrm{V}}
\def\U{\mathrm{U}}
\def\W{\mathrm{W}}
\def\J{\mathrm{J}}

\def\R{\mathbb R}

\def\N{\mathbb N}
\def\P{\mathbb P}
\def\T{\mathbb T}

\def\det{\mathrm{det}}

\def\Ld{\mathrm{L}}
\def\H{\mathrm{H}}

\def\B{\mathrm{B}}
\def\E{\mathrm{E}}
\def\W{\mathrm{W}}
\def\D{\mathrm{D}}

\newcommand{\eps}{\varepsilon}

\newtcolorbox{dev}{arc=0pt}
\newcounter{compteur}
\counterwithin{compteur}{section}

\definecolor{thmcolor}{rgb}{0.8,0.14,0.2}
\definecolor{defcolor}{rgb}{0.0,0.50,0.0}
\definecolor{excolor}{rgb}{0.50,0.0,0.990}
\definecolor{applicolor}{rgb}{0.50,0.0,0.990}

\newtheoremstyle{thm}% name of the style to be used
  {\topsep}% measure of space to leave above the theorem. E.g.: 3pt
  {\topsep}% measure of space to leave below the theorem. E.g.: 3pt
  {\itshape}% name of font to use in the body of the theorem
  {0pt}% measure of space to indent
  {\bfseries}% name of head font
  {.}% punctuation between head and body
  { }% space after theorem head; " " = normal interword space
  {\textcolor{black!100}{\thmname{#1}\thmnumber{ #2}}\thmnote{ (#3)}}

\newtheoremstyle{def}% name of the style to be used
  {3pt}% measure of space to leave above the theorem. E.g.: 3pt
  {3pt}% measure of space to leave below the theorem. E.g.: 3pt
  {}% name of font to use in the body of the theorem
  {0pt}% measure of space to indent
  {\bfseries}% name of head font
  {.}% punctuation between head and body
  { }% space after theorem head; " " = normal interword space
  {\textcolor{black!100}{\thmname{#1}\thmnumber{ #2}}\thmnote{ (#3)}}
  
\newtheoremstyle{ex}% name of the style to be used
  {1pt}% measure of space to leave above the theorem. E.g.: 3pt
  {1pt}% measure of space to leave below the theorem. E.g.: 3pt
  {}% name of font to use in the body of the theorem
  {0pt}% measure of space to indent
  {\bfseries}% name of head font
  {.}% punctuation between head and body
  { }% space after theorem head; " " = normal interword space
  {\textcolor{black!100}{\thmname{#1}\thmnumber{ #2}\thmnote{ (#3)}}}
  
\theoremstyle{thm}
\newtheorem{thm}[compteur]{Theorem}
\newtheorem{coro}[compteur]{Corollary}
\newtheorem{propo}[compteur]{Proposition}
\theoremstyle{appli}

\newtheorem{lem}[compteur]{Lemma}
\theoremstyle{def}
\newtheorem{defi}[compteur]{Definition}
\newtheorem{nota}[compteur]{Notation}
\theoremstyle{ex}

\newtheorem{rem}[compteur]{Remark}

\date{}

\numberwithin{equation}{section}

\subsectionfont{\normalfont\Large\bfseries}
\sectionfont{\normalfont\LARGE\bfseries}
\newcommand{\nocontentsline}[3]{}
\newcommand{\tocless}[2]{\bgroup\let\addcontentsline=\nocontentsline#1{#2}\egroup}

\title{Decay and absorption for the Vlasov-Navier-Stokes system with gravity in a half-space}
\author{Lucas Ertzbischoff 
 \thanks{Department of Mathematics, Imperial College London, London, SW7 2AZ, United-Kingdom (\href{mailto:l.ertzbischoff@imperial.ac.uk}{l.ertzbischoff@imperial.ac.uk})}
}

\begin{document}

\maketitle
\begin{abstract}
This paper is devoted to the large time behavior of weak solutions to the three-dimensional Vlasov-Navier-Stokes system set on the half-space, with an external gravity force. This fluid-kinetic coupling arises in the modeling of sedimentation phenomena. 
Our main result establishes the convergence of the local density of particles and the fluid velocity to $0$ in large time, with a polynomial rate of decay. In order to overcome the effect of the gravity, we rely on a fine analysis of the absorption phenomenon at the boundary. We obtain a family of decay estimates for the moments of the kinetic distribution,  provided that the initial distribution function has a sufficient decay in the phase space.
%We prove that the local density of particles enjoys a polynomial decay to $0$ in large time in various $\Ld^p$ norms while the fluid velocity tends to $0$ in $\Ld^2$ norm, also at a polynomial rate. 
%In order to overcome the possible injection of energy coming from the gravity force, 
%we rely on a fine analysis of the absorption effect at the boundary. This one naturally provides a family of decay estimates for the moments of the kinetic distribution,  provided that the initial distribution function has a sufficient decay in the phase space.
%allowing to capture the behavior of both fluid and dispersed phases.
\end{abstract}

\renewcommand{\contentsname}{Contents}

{ \hypersetup{linkcolor=black}
\tableofcontents
}

\section{Introduction}

Fluid-kinetic systems aim at modelling the collective motion of a dispersed phase of small particles immersed within a fluid. In such systems, also called spray models, the dispersed phase is described at a mesoscopic level by a distribution function solving a kinetic equation while the evolution of macroscopic quantities for the fluid is governed by fluid mechanics equations.

Among the wide family of fluid-kinetic systems (see the pioneering works \cite{will,oro}), one can consider the so-called \textit{thin spray} models where the volume fraction of the particles is small compared to that of the surrounding fluid. In this context, an interesting prototype is the incompressible Vlasov-Navier-Stokes system, coupling in a nonlinear way the fluid and kinetic equations through a drag term. The later one depends on the fluid unknowns and on the density function and allows for an exchange of momentum between the fluid and the particles. Beyond its mathematical interest, this system also appears in the study of the transport and deposition of a therapeutic aerosol in the airflows contained in the human upper airways (see \cite{BGLM}).

In this paper, we are interested in the following Vlasov-Navier-Stokes system set in $\R^3_+\times \R^3$:
\begin{align}
\partial_t f +v \cdot \nabla_x f  & + {\rm div}_v [f(u-v) + fG ]=0, \ \ &&(t,x,v) \in \R^{*}_{+} \times \R^3_+ \times \R^3, \label{eq:Vlasov} \\
\partial_t u + (u \cdot\nabla)u- \Delta u + \nabla p &=\int_{\R^3} f(t,x,v)(v-u(t,x)) \, \mathrm{d}v,  \ \ &&(t,x) \in \R^{*}_{+} \times \R^3_+, \label{eq:NS} \\ 
\mathrm{div} \, u &=0, \ \ &&(t,x) \in \R^{*}_{+} \times \R^3_+. 
\label{eq:NS2}
\end{align}
Here $\R^3_+ := \R^2 \times (0,+\infty)$ is the tridimensional half-space while $G:=(0,0,-g) \in \R^3 $ is a given vector, with $g>0$ the constant gravitational acceleration. In these equations, $u=u(t,x) \in \R^3$ and $p=p(t,x) \in  \R$ stand for the velocity field and pressure of the fluid, while $f=f(t,x,v)\in  \R^+$ is the distribution function of the particles in the phase space $\R^3_+ \times \R^3$. Here, the particles undergo the friction force produced by the surrounding fluid, as well as the effect of gravity. Thus, using the Stokes law, the resultant force exerted on the particles is the sum of the drag and weight/buoyancy, that is 
\begin{align*}
u(t,x)-v + G,
\end{align*}
and the Vlasov equation \eqref{eq:Vlasov} is thus coupled to the Navier-Stokes equations \eqref{eq:NS}-\eqref{eq:NS2}. Note that we have implicitely considered spherical particles of same radius and that the mass density of the particles is greater than that of the fluid (because of the positive coefficient before the vector $G$ - see e.g. \cite{CB}).

A coupling term is also added in the Navier-Stokes equations \eqref{eq:NS}-\eqref{eq:NS2}, where a forcing term appears in the right-hand side and stems from the retroaction of the particles on the fluid. This source term is usually called the \textit{Brinkman force} and can be rewritten as
\begin{align}\label{eq:Brinkman}
\int_{\R^3} f(v-u) \, \mathrm{d}v=j_f-\rho_f u ,
\end{align}
where
\begin{align*}
\rho_f(t,x)&:= \int_{\R^3} f(t,x,v) \, \mathrm{d}v, \\
j_f(t,x)&:= \int_{\R^3} v f(t,x,v) \, \mathrm{d}v.
\end{align*}
Note that in the previous Navier-Stokes equations \eqref{eq:NS}-\eqref{eq:NS2}, the density and viscosity of the fluid are assumed to be constant and both chosen equal to $1$, while the external gravity force $G=-\nabla \Phi$, with $\Phi=gz$, has been absorbed in the pressure term. In short, the equations \eqref{eq:Vlasov}-\eqref{eq:NS}-\eqref{eq:NS2} account for the description of a cloud of fine particles sedimenting in an ambient incompressible viscous fluid.

\bigskip

%Unlike the \textit{thick spray} case where a very more involved system may arise from the mathematical point of view,

The main goal of this article is to study the asymptotics in large time of small-data solutions to the Vlasov-Navier-Stokes system, relying on some specific boundary conditions that we shall detail below. The analysis of this system has been explored in different directions over the past two decades.

\bigskip

\textbf{The Cauchy problem.} Concerning the existence theory of global weak solutions to the system \eqref{eq:Vlasov}-\eqref{eq:NS}-\eqref{eq:NS2} (without the gravity force), different settings and boundary conditions for the distribution function have been adressed, depending on the spatial domain: a fixed bounded domain with specular reflexion in \cite{ABdM}, the flat torus $\T^3$ in \cite{BDGM}, time-dependent domain with absorption boundary condition in \cite{BGM,BMM}, or a 2D rectangle with partly absorbing boundary condition in \cite{GHKM}. Local strong solutions can also be considered as in \cite{ChKw} (for inhomogeneous fluid equations), as well as blow-up in finite time of classic solutions in \cite{choi2017finite}. 

Note that the additional term involving the gravity has only been taken into account for the Vlasov-Stokes system on bounded domain with specular reflexion in \cite{Ham}, or for the Vlasov equation coupled to the stationary Stokes system with a regular and compactly supported initial distribution function on $\R^3$ in \cite{HoferInertia}.
\bigskip

\textbf{The 2D case.} \ In two dimensions, more results are available for this fluid-kinetic system, essentially because of the study of the Navier-Stokes system which is more favorable in this context: for instance, uniqueness of 2D-global weak solutions is proven for the whole space or the torus case in \cite{HKM3}. The controllabillity of the system is also explored in dimension $2$ in \cite{Moyano}. 

\bigskip

\textbf{Link to other models.} In the spirit of Hilbert's 6th problem of axiomatization of physics, fluid-kinetic models can be linked to other systems of ODEs and PDEs.
Deriving rigorously the Vlasov-Navier-Stokes system from “first laws" appears as an important issue, but remains essentially an outstanding open problem for the whole system. Two main strategies have been proposed so far. The \textit{mean-field} limit of a $\mathrm{N}$-solid particle system coupled with a fluid equation has been considered in \cite{DGR,H,HMS,CH} and allows one to recover the Brinkman force in a quasi-static framework thanks to homogenization techniques. Some partial results are also known for a dynamical but only macroscopic equation for the particles, in some dilute regime where they have no inertia (see \cite{Hof,Mech}).
Another direction has been taken in \cite{BDGR1,BDGR2} where, in a formal way, the coupling between two mixtures provides the derivation of the Vlasov-Navier-Stokes system in the same fashion as the hydrodynamic limits of the Boltzmann equation \cite{bardos1991fluid1,bardos1993fluid2}. However, a full and rigorous justification of this program based on kinetic theory is still open.

Through \textit{hydrodynamic limits} of the Vlasov-Navier-Stokes system, one can also seek to derive some systems involving only averaged quantities: more precisely, high friction regimes of the system have been shown to lead to Navier-Stokes type systems. These asymptotic regimes have been first considered in \cite{goudon2004hydrodynamic1,goudon2004hydrodynamic2} for the Vlasov-Fokker-Planck-Navier-Stokes equations, where the effect of Brownian motion is added in the equation of the distribution function. 

Without diffusion in velocity in the Vlasov equation, one of these limits has been handled in \cite{HoferInertia} for the Vlasov-(steady)Stokes system with gravity in the whole space. Very recently, the question raised by these different regimes has been adressed in \cite{HKM} for the full system \eqref{eq:Vlasov}-\eqref{eq:NS}-\eqref{eq:NS2} on the torus (without the gravity force).

\bigskip

\textbf{Large time behavior.} \ The \textit{large time dynamics} of the Vlasov-Navier-Stokes system, which is the main issue of this article, has very recently received a particular attention. This natural question is studied for the first time by Jabin in \cite{Jab} for a reduced kinetic model. 
In the absence of dissipative mechanism in the Vlasov equation (like a Fokker-Planck operator allowing to consider smooth equilibria, see \cite{GLMZ}), the sole effect of the drag force in the system should lead to nontrivial equilibria which are singular. More precisely, one expects a monokinetic behavior for the distribution function of the particles (that is to say, a convergence towards a Dirac mass in velocity). 
A conditional result accounting for this phenomenon in the Vlasov-Navier-Stokes system has been provided by Choi and Kwon in \cite{ChKw}. In short, it requires a global bound in time which is not \textit{a priori} satisfied by global weak solutions to the system. More recently, this extra assumption has been removed for initial data which are in some sense close to equilibrium: the first complete result stems from the article \cite{HKMM} of Han-Kwan, Moussa and Moyano where the authors work in a periodic setting and in a framework \textit{à la} Fujita-Kato. In the same spirit, such a monokinetic behavior of weak solutions has been obtained for small data in the whole space case by Han-Kwan in \cite{HK} and then extended to the case of bounded domains with absorption boundary conditions in \cite{EHKM}. In short, these results are all based on a remarkable energy-dissipation inequality satisfied by weak solutions to the system.

Unlike this series of works, the existence and stability of regular equilibria has been obtained by Glass, Han-Kwan and Moussa in \cite{GHKM} for a particular 2D bounded domain with partly absorbing and injection boundary conditions. 
%If one considers an additional dissipation operator in the kinetic equation ( that is a term of the form $- \Delta_v f$), 
%When a Fokker-Planck dissipation term (namely, a term of the form $- \Delta_v f$ ) is added in the kinetic equation (\ref{eq:Vlasov}), global classical solutions can be constructed for data close to Maxwellian equilibria and this non-singular steady states locally attract these solutions (see \cite{GLMZ}).
Finally, let us emphasize the fact that the high friction limit tackled by Han-Kwan and Michel in \cite{HKM} is closely related to the monokinetic behavior we mentioned earlier, and in particular to the techniques used in \cite{HKMM,HK}.

\medskip

\textbf{Main contribution of this paper.} \ In the continuation of these previous works, the main goal of this article is the study of the large time dynamics of global weak solutions to the system \eqref{eq:Vlasov}-\eqref{eq:NS}-\eqref{eq:NS2}. Its originality lies in dealing with a fluid-kinetic system on an unbounded domain with boundary, where one considers absorption boundary conditions for the distribution function, together with an additional gravity force term. 

Loosely speaking, the presence of a gravity force may ruin the decay of  the energy of the system. At first sight, this prevents the use of exactly the same techniques as in \cite{HKMM,HK,EHKM}. However, it is actually possible to take advantage of the absorption at the boundary to analyse the large time behavior of global weak solutions starting close to equilibrium.
%induces a continuous injection of momentum into the system resulting \textit{a priori} in a non-decaying total kinetic energy, while the gravity-less case studied in \cite{HKMM,HK,EHKM} is really based on the potential decay of the later. 
%Nevertheless, a careful analysis of the boundary effects will reveal that the absorption of the particles can be used to get some decay of the moments of the distribution function at any time, modulo an initial distribution function enjoying a sufficient decay in space and velocity itself.
%provided that the initial distribution function itself enjoys a sufficient decay in space and velocity.

\bigskip

Before going further, we give several definitions and set notations about the system that we will consider in this article. In what follows, we will sometimes refer to (\ref{eq:Vlasov})-(\ref{eq:NS})-(\ref{eq:NS2}) as the \textit{VNS system}. Along this paper, we will make a constant use of the notation
\begin{align*}
\R^3_+ := \R^2 \times (0,+\infty).
\end{align*}

First, the VNS system is supplemented with the following initial conditions for $u$ and $f$:
\begin{align}
u_{\mid t=0}&=u_0 \text{ in } \R^3_+,\\
f_{\mid t=0}&=f_0 \text{ in } \R^3_+ \times \R^3.
\end{align}
We prescribe the following Dirichlet boundary conditions for the fluid: 
\begin{align}\label{bcond-fluid}
u(t,\cdot)=0, \text{ on } \partial \R^3_+=\R^2 \times \lbrace 0 \rbrace.
\end{align}
We also need to introduce the following outgoing/incoming phase-space boundaries:
\begin{align}\label{def:Sigma}
\Sigma^{\pm}&:= \left\lbrace  (x,v) \in \partial \R^3_+ \times \R^3   \mid \pm v \cdot n(x)>0 \right\rbrace,\\
\Sigma_0&:= \left\lbrace  (x,v) \in \partial \R^3_+ \times \R^3   \mid v \cdot n(x)=0 \right\rbrace,\\
\Sigma &:= \Sigma^+ \sqcup \Sigma^- \sqcup \Sigma_0= \partial \R^3_+ \times \R^3,
\end{align}
where $n(x)$ stands for the normal vector to the boundary $\partial \R^3_+$ at point $x$. We observe that 
\begin{align*}
\Sigma^{\pm}&= \left\lbrace  (x,v) \in \R^3 \times \R^3   \mid x_3=0, \  \pm v_3<0 \right\rbrace,\\
\Sigma_0&= \left\lbrace  (x,v) \in  \R^3 \times \R^3   \mid x_3=0, \   v_3=0 \right\rbrace.
\end{align*}
Then, we prescribe the following absorption boundary conditions for the distribution function:
\begin{align}\label{bcond-f}
f(t,\cdot,\cdot)=0, \text{ on } \Sigma^{-},
\end{align}
meaning that particles reaching transversally the physical boundary $\lbrace x_3=0 \rbrace$ are absorbed.

\medskip

Several functionals play an important role in the study of the VNS system. We introduce the following ones. 
\begin{defi}
\begin{enumerate}
\item The \textbf{kinetic energy} of the Vlasov-Navier-Stokes system is defined for all $t \geq 0$ as:
\begin{align}\label{eq:Energy}
\mathrm{E}(t):= \dfrac{1}{2}\Vert u(t) \Vert_{\Ld^2(\R^3_+)}^2
+\dfrac{1}{2}\int_{\R^3_+ \times \R^3} f(t,x,v) \vert v \vert^2 \, \mathrm{d} x \, \mathrm{d} v.
\end{align}  

\item The \textbf{dissipation} of the Vlasov-Navier-Stokes system (without gravity) is defined for all $t \geq 0$ as:
\begin{align}\label{eq:Dissipation}
\mathrm{D}(t):= \int_{\R^3_+ \times \R^3} f(t,x,v) \vert u(t,x)-v\vert^2 \, \mathrm{d} x \, \mathrm{d} v 
 +\Vert \nabla u(t) \Vert_{\Ld^2(\R^3_+)}^2.
\end{align}

\item The \textbf{dissipation with gravity} of the Vlasov-Navier-Stokes system is defined for all $t \geq 0$ as:
\begin{align}\label{eq:Dissipation}
\mathrm{D}_G(t):= \int_{\R^3_+ \times \R^3} f(t,x,v) \vert u(t,x)-v\vert^2 \, \mathrm{d} x \, \mathrm{d} v 
 +\Vert \nabla u(t) \Vert_{\Ld^2(\R^3_+)}^2 -\int_{\R^3_+} G \cdot j_f(t,x) \, \mathrm{d}x.
\end{align}  
\end{enumerate}
\end{defi}

At a formal level, the VNS system enjoys an energy-dissipation structure involving the previous functionals. More precisely, smooth solutions to the system satisfy the following \textit{a priori} estimate
\begin{align*}
\dfrac{\mathrm{d}}{\mathrm{d}t}\mathrm{E}(t) + \mathrm{D}_G(t)\leq0.
\end{align*}
Thus, there is a variation of energy coming from the dissipation inside of the fluid and from the friction between the particles and the fluid, but the gravity force induces an additional term leading to a potential non-decay of the kinetic energy $\E$: indeed, we expect the term $G  \cdot j_f>0$ to be positive at least after some time because the particles should ultimately fall in the same direction as the gravity.
% Indeed, since the term $\int_{\R^3_+} G \cdot j_f(t,x) \, \mathrm{d}x$ does not have a fixed known sign, a possible decay of the total kinetic energy $t \mapsto \E(t)$ is not given \textit{a priori} by the previous inequality alone. Note that the case where $G  \cdot j_f>0$ (which is the non-favorable case with respect to the decay of the energy) shall be naturally appearing after some time.

\bigskip

We denote by $\mathscr{D}_{\mathrm{div}}(\R^3_+)$ the set of smooth $\R^3$ valued divergence free vector-fields having compact support in $\R^3_+$. The closures of $\mathscr{D}_{\mathrm{div}}(\R^3_+)$ in $\Ld^2(\R^3_+)$ and in $\H^1(\R^3_+)$ are respectively denoted by $\Ld^2_{\mathrm{div}}(\R^3_+)$ and by $\H^1_{0,\mathrm{div}}(\R^3_+)$. We write $\H^{-1}_{\mathrm{div}}(\R^3_+)$ for the dual of the later.

\bigskip

We now define the class of admissible initial data for the VNS system.
\begin{defi}[Initial condition]\label{CIadmissible}
We shall say that a couple $(u_0,f_0)$ is an admissible initial condition if: 
\begin{align}
u_0 &\in \Ld^2(\R^3_+), \ \ \mathrm{div}_x \, u_0=0,  \label{CI:fluid}\\
f_0 &\in \Ld^1 \cap \Ld^{\infty}(\R^3_+ \times \R^3), \label{CI-f1}\\
f_0 &\geq 0, \ \ \int_{\R^3_+ \times \R^3} f_0 \, \mathrm{d}x \, \mathrm{d}v =1, \label{CI-f2}\\
(x,v)& \mapsto f_0 (x,v)\vert v \vert^2 \in \Ld^1(\R^3_+ \times \R^3). \label{CI-f3}
\end{align}
\end{defi}
We then introduce some notations about the moments of any phase-space distribution function.
\begin{defi}\label{notation:moments}
For any $\alpha \geq 0$ and any measurable function $g:\R^+ \times \R^3_+ \times \R^3 \rightarrow \R^+$, we set
\begin{align*}
m_{\alpha}g(t,x)&:=\int_{\R^3} \vert v \vert^{\alpha} g(t,x,v)\,\mathrm{d}v, \\
M_{\alpha}g(t)&:=\int_{\R^3_+ \times \R^3} \vert v \vert^{\alpha} g(t,x,v) \, \mathrm{d}v \, \mathrm{d}x=\int_{\R^3_+} m_{\alpha} g(t,x) \, \mathrm{d}x.
\end{align*} 
\end{defi}

In our approach, we shall rely on some decay assumptions satisfied by the initial distribution function $f_0$. We thus introduce the following quantities.
\begin{defi}\label{decay_f0}
For any $q>0$, $m>0$ and $r \geq 1$, we set
\begin{align}
\label{def:Nq}N_q(f_0)&:= \Vert (1+\vert v \vert^q) f_0 \Vert_{\Ld^{\infty}(\R^3_+ \times \R^3)}, \\[1mm]
\label{def:Kp}K_{q,r}(f_0)&:= \left\Vert(1+\vert v \vert^q) \Vert f_0(\cdot,v) \Vert_{\Ld_x^r(\R^3_+)} \right\Vert_{{\Ld}^{\infty}_v(\R^3)}, \\[1mm]
%\label{def:Hq}H_q(f_0)&:= \underset{\substack{x \in \R^3_+ \\ v \in \R^3}}{\sup} (1+x_3^q) f_0(x,v), \\[1mm]
\label{def:Hq}H_{q,m}(f_0)&:=\left\Vert (1+\vert v \vert^{q}) \Vert (1+x_3^m) f_0(\cdot,v) \Vert_{\Ld_x^{\infty}(\R^3_+)} \right\Vert_{\Ld_v^1(\R^3)}, \\[1mm]
%\label{def:Fq}F_{q,r}(f_0)&:= \underset{v \in \R^3}{\sup} \, \Vert x_3^q \, f_0(\cdot, v) \Vert_{\Ld^r(\R^3_+)}.
\label{def:Fq}F_{q,m,r}(f_0)&:= \left\Vert (1+\vert v \vert^{q}) \Vert (1+x_3^m) f_0(\cdot,v) \Vert_{\Ld_x^{r}(\R^3_+)} \right\Vert_{\Ld_v^1(\R^3)}.
\end{align}
\end{defi}
%
%\begin{defi}\label{def:pointdecay}
%We say that an initial kinetic condition $f_0$ satisfies the pointwise decay assumption of order $q>0$ if 
%\begin{align}\label{def:Nq}
%N_q(f_0):= \underset{\substack{x \in \R^3_+ \\ v \in \R^3}}{\sup} (1+\vert v \vert^q) f_0(x,v) < \infty.
%\end{align}
%\end{defi}
%
%\begin{defi}\label{def:pointdecay2}
%We say that an initial kinetic condition $f_0$ satisfies the pointwise decay assumption of order $q>0$ if 
%\begin{align}\label{def:Hq}
%H_q(f_0):= \underset{\substack{x \in \R^3_+ \\ v \in \R^3}}{\sup} (1+x_3^q) f_0(x,v) < \infty.
%\end{align}
%\end{defi}
%
%\begin{defi}\label{def:mixdecay}
%We say that an initial kinetic condition $f_0$ satisfies a mixed decay assumption in velocity of order $(q,r)$ (with $q>0$ and $r \geq 1$) if 
%\begin{align}\label{def:Kp}
%K_{q,r}(f_0):= \underset{v \in \R^3}{\sup} \, (1+\vert v \vert^q) \Vert f_0(\cdot,v) \Vert_{\Ld^r(\R^3_+)} < \infty.
%\end{align}
%\end{defi}
%
%\begin{defi}\label{def:mixdecay2}
%We say that an initial kinetic condition $f_0$ satisfies a mixed decay in space assumption in velocity of order $(q,r)$ (with $q>0$ and $r \geq 1$) if 
%\begin{align}\label{def:Fq}
%F_{q,r}(f_0):= \underset{v \in \R^3}{\sup} \, \Vert x_3^q f_0(\cdot, v) \Vert_{\Ld^r(\R^3_+)} < \infty.
%\end{align}
%\end{defi}

We will consider weak solutions to the Vlasov equation with gravity force (\ref{eq:Vlasov}), with the boundary condition (\ref{bcond-f}) and the previous initial conditions, which are defined as follows.

\begin{defi}[Weak solutions to the Vlasov equation]\label{sol:Vlasov}
Let $\U \in \Ld^1_{\mathrm{loc}}(\R^+ \times \R^3_+)$ and an initial distribution $f_0$ satisfying (\ref{CI-f1})-(\ref{CI-f2})-(\ref{CI-f3}). We say that a nonnegative function $f \in \Ld^{\infty}_{\mathrm{loc}}(\R^+;\Ld^1 \cap \Ld^{\infty}(\R^3_+ \times \R^3))$ is a weak solution to the Vlasov equation (\ref{eq:Vlasov}) with force field $\U$, with boundary condition (\ref{bcond-f}) and with initial condition $f_0$ if, for all $\Psi \in \mathscr{D}([0,T] \times \overline{\R^3_+} \times \R^3)$ with $\Psi(T,\cdot)=0$ and vanishing on $\R^+ \times (\Sigma^+ \cup \Sigma_0)$, one has
\begin{multline*}
\int_0^T \int _{\R^3_+ \times \R^3} f(t,x,v) \left[  \partial_t \Psi + v \cdot \nabla_x \Psi +(\U-v) \cdot \nabla_v \Psi \right](t,x,v) \, \mathrm{d}x \, \mathrm{d} v \, \mathrm{d}t   =-\int_{\R^3_+ \times \R^3} f_0(x,v) \Psi(0,x,v) \, \mathrm{d}x \,\mathrm{d}v.
\end{multline*}
\end{defi}
Weak solutions to the Vlasov equations enter in the framework of the DiPerna-Lions theory for transport equations (in the phase space $\R^3_+ \times \R^3$). We refer to Section \ref{DiPernaLions} in the Appendix for more details, where we recall in particular the classic stability property of renormalized solutions that we will constantly use throughout this article. 

We will also consider weak solutions for the full Vlasov-Navier-Stokes with the boundary conditions (\ref{bcond-fluid})-(\ref{bcond-f}) and the initial conditions described in (\ref{CIadmissible}), in the following sense.
\begin{defi}[Weak solutions with strong energy inequality for the VNS system]\label{sol-faible}
Consider an admissible initial condition $(u_0,f_0)$ in the sense of Definition \ref{CIadmissible}. A global weak solution to the Vlasov-Navier-Stokes system with boundary condition (\ref{bcond-fluid})-(\ref{bcond-f}) and with initial condition $(u_0,f_0)$ is a pair $(u,f)$ such that:
\begin{align}
&u\in \Ld^{\infty}_{\mathrm{loc}}(\R^+;\Ld^2_{\mathrm{div}}(\R^3_+))\cap  \Ld^{2}_{\mathrm{loc}}(\R^+;\H^1_{0,\mathrm{div}}(\R^3_+)),\\
{\rm div}_x& \, u=0,  \\
f &\in \Ld^{\infty}_{\mathrm{loc}}(\R^+;\Ld^1 \cap \Ld^{\infty}(\R^3_+ \times \R^3)),\\
j_f-\rho_f u &\in \Ld^{2}_{\mathrm{loc}}(\R^+;\H^{-1}_{\mathrm{div}}(\R^3_+)),
\end{align}
and such that the following holds. The distribution function $f$ is a weak solution to the Vlasov equation with force field $G+u$ with initial condition $f_0$ in the sense of Definition \ref{sol:Vlasov} and the velocity field $u$ is a Leray solution to the Navier-Stokes equations with initial condition  $u_0$, that is for all $\Phi \in \mathscr{C}^1([0,T];\mathscr{D}_{\mathrm{div}}(\R^3_+))$ such that $\Phi(T,\cdot)=0$, we have
\begin{multline}\label{weak:formulNS}
\int_0^T\int_{\R^3_+}\left[ u \cdot \partial_t \Phi +(u \otimes u) : \nabla_x \Phi -\nabla_x u : \nabla_x \Phi \right](t,x) \, \mathrm{d} x \, \mathrm{d}t \\
= -\int_{\R^3_+} u_0(x) \cdot \Phi(0,x) \, \mathrm{d}x - \int_0^T \int_{\R^3_+ \times \R^3} f(t,x,v) \left[v-u(t,x)\right] \cdot \Phi (t,x)  \, \mathrm{d}x \, \mathrm{d} v \, \mathrm{d}t, 
\end{multline}
and the strong energy inequality holds for the Navier-Stokes equations: for any $t \geq 0$ and almost every $0 \leq s  \leq t$ (including $s=0$)
\begin{align}\label{ineq-energyNS}
\Vert u(t) \Vert_{\Ld^2(\R^3_+)}^2  + 2 \int_s ^t \Vert \nabla u(\tau) \Vert_{\Ld^2(\R^3_+)}^2 \mathrm{d}\tau \leq   \Vert u(s) \Vert^2_{\Ld^2(\R^3_+)} + 2 \int_s ^t \int_{\R^3_+} \left( j_f(\tau,x)-\rho_f u(\tau,x) \right) \cdot u(\tau,x) \, \mathrm{d}x \,  \mathrm{d}\tau.
\end{align}
Furthermore, the following energy estimate holds for the Vlasov-Navier-Stokes system: for any $t \geq 0$ and almost every $0 \leq s \leq t$ (including $s=0$) 
\begin{align}\label{ineq-energy}
\mathrm{E}(t) + \int_s ^t \mathrm{D}(\tau)\mathrm{d}\tau \leq \mathrm{E}(s)+\int_{s}^t \int_{\R^3_+ \times \R^3} G \cdot v f(\tau,x,v) \, \mathrm{d}x \, \mathrm{d}v \, \mathrm{d}\tau ,
\end{align}
where the energy $\E$ and dissipation $\mathrm{D}$ have been defined in \eqref{eq:Dissipation}.
\end{defi}

\medskip

Note that the last integral in the right-hand side of the inequality \eqref{ineq-energyNS} actually makes sense because of Sobolev embedding and the fact that $j_f-\rho_f u \in \Ld^2_{\mathrm{loc}}(\R^+;\Ld^{6/5}(\R^3_+))$ (see Section \ref{Section:PBCauchy} in the Appendix).

Such global weak solutions can be obtained through an approximation procedure which seems to be classic by now (see e.g. \cite{BGM, GHKM, BMM}). Since the half-space/gravity framework has not been explicitly treated in the former literature, we provide some rather sketchy elements of proof about the Cauchy problem in Section \ref{Section:PBCauchy}, with a particular insight on the obtention of the strong energy inequalities \eqref{ineq-energyNS} and \eqref{ineq-energy}.

\bigskip

In Section \ref{Section:Results+Strat}, we present the main result obtained in this paper. As we will explain later, the approach we will use to prove this result shares some similar features with the ones introduced in \cite{HKMM,HK,EHKM,GHKM}. We will detail the strategy set up in these works and in comparison, describe the method we need in our case.

\section{Main results}\label{Section:Results+Strat}
First, let us provide an informal statement of the main theorem of this article. Under some smallness assumption on the initial data $(u_0,f_0)$, we shall prove that the fluid velocity $u$ and the local density of particles $\rho_f$ decay to $0$ in large time, in the following sense: any global weak solution $(u,f)$ to the VNS system with small data satisfies for all $t>0$ and $r \in [1,+\infty]$
\begin{align*}
\left\Vert u(t) \right\Vert_{\Ld^{2}(\R^3_+)} \lesssim \frac{1}{(1+t)^{3/4}}, \ \ \Vert \rho_f(t) \Vert_{\Ld^r(\R^3_+)} \lesssim \dfrac{1}{(1+t)^{k}},
\end{align*}
for some $k>0$, where $\lesssim$ depends on the initial data. For instance, the exponent $k$ will be related to the decay of $f_0$ in space-velocity.
\medskip

Before stating our results, we define some quantities based on the regularity and decay of the initial data $(u_0,f_0)$. The notations we use here are introduced in Definitions \ref{notation:moments}--\ref{decay_f0} and in Section \ref{AnnexeMaxregStokes} of the Appendix. If $p,s>0$ are given, we set
\begin{align*}
\mathcal{E}_{p,s}&:= \W_0^{1,\frac{9}{7}}(\R^3_+) \cap \mathrm{D}_3^{\frac{1}{2},2}(\R^3_+) \cap \D_3^{1-\frac{1}{s},s}(\R^3_+) \cap  \mathrm{D}_p^{1-\frac{1}{p},p}(\R^3_+),\\[2mm]
\mathscr{E}(0)&:=\Vert u_0 \Vert_{\H^1 \cap \mathrm{D}_p^{1-\frac{1}{p},p}(\R^3_+)}\hspace{-4mm}+\mathrm{E}(0)+\Vert u_0 \Vert_{\Ld^1(\R^3_+)}+ M_{6} f_0,
\end{align*}
and for $q,m>0$
\begin{align*}
\mathscr{N}_{q,m}(f_0):=N_q (f_0)+H_{q,m}(f_0)+ \underset{r \in \left\lbrace 1,3 \right\rbrace}{\max} \left\lbrace K_{q,r}(f_0)  +F_{q,m,r}(f_0) \right\rbrace.
\end{align*}
The main result of this paper reads as follows.
\begin{thm}\label{thm1}
Let $s \in (2,3)$, $q \geq 7$ and $m \geq 4$. There exist $p_0 \in \left(3, \frac{3(2+s)}{4} \right)$ and $\mathrm{C}>0$ such that the following holds for all $p \in (3,p_0)$. Let $(u_0,f_0)$ be an admissible initial condition in the sense of Definition \ref{CIadmissible} satisfying 
\begin{align}\label{data:hyp}
\begin{split}
& u_0 \in \H^1_{0,\mathrm{div}}(\R^3_+) \cap \Ld^{1}(\R^3_+) \cap \mathcal{E}_{p,s} , \\
& M_6 f_0+ \mathscr{N}_{q,m}(f_0) <\infty. \\
%& N_q(f_0), \, H_{q,m}(f_0) <\infty,\\
%&K_{q,r}(f_0), \, F_{q,m,r}(f_0)<\infty, \ \ \text{for any} \ \ r \in \left\lbrace 1,3 \right\rbrace.
%& \mathrm{supp}_x \, f_0 \subset \R^2 \times (0,L), \ \ \text{for some finite} \ \ L>0.
\end{split}
\end{align}
There exists a constant $C_0=C_0(\mathscr{E}(0),\mathscr{N}_{q,m}(f_0))>0$
%\begin{multline}\label{smallness:condition:THM}
%\Phi_{\mathrm{C}_0} \Big( 1+\Vert u_0 \Vert_{\H^1 \cap \mathrm{D}_p^{1-\frac{1}{p},p}(\R^3_+)}\hspace{-4mm}+ N_q (f_0)+H_{\widetilde{q}}(f_0)+ K_{m,r}(f_0)  +F_{\widetilde{m},\widetilde{r}}(f_0)+ M_{6} f_0  + \mathrm{E}(0)+\Vert u_0 \Vert_{\Ld^1(\R^3_+)} \Big) \\ \times \left( \Vert u_0 \Vert_{\H^1(\R^3_+)}^2+\Vert u_0 \Vert_{\Ld^1(\R^3_+)}^2+N_q (f_0)+H_{\widetilde{q}}(f_0)\right)<\mathrm{C}_0,
%\end{multline}
such that if
\begin{align}\label{smallness:condition:THM}
C_0 \times  \left( \Vert u_0 \Vert_{\H^1\cap \Ld^1(\R^3_+)}+N_q (f_0)+H_{0,m}(f_0)\right)<\mathrm{C},
\end{align}
then the following holds: there exists a constant $\Lambda_0=\Lambda_0\left(\Vert u_0 \Vert_{\H^1\cap \Ld^1(\R^3_+)},N_q (f_0),H_{0,m}(f_0) \right)>0$
such that any global Leray solution $(u,f)$ to the Vlasov-Navier-Stokes system with initial data $(u_0,f_0)$ (in the sense of Definition \ref{sol-faible}) satisfies for all $t>0$, $\overline{q} \geq q$ and $k \in [0,\overline{q}-3)$
\begin{align}
\label{theo:estimu} \left\Vert u(t) \right\Vert_{\Ld^{2}(\R^3_+)} &\leq \frac{\Lambda_0}{(1+t)^{3/4}}, \\[2mm]
\label{theo:estim-rhofPONCT}\Vert \rho_f(t) \Vert_{\Ld^{\infty}(\R^3_+)} & \leq C_{k,\overline{q},g}  \dfrac{N_{\overline{q}}(f_0)+H_{0,k}(f_0)}{(1+t)^{k}},
\end{align}
for some constant $C_{k,\overline{q},g}>0$.
\end{thm} 

In particular, this shows that the more the initial data $f_0$ decays in the phase space, the more the local density $\rho_f$ enjoys some decay in time.

We can also prove the following result for the decay of $\rho_f(t)$ in $\Ld^r(\R^3_+)$.
\begin{propo}\label{PropoStatementMAIN:estimLP}
Consider the same assumptions \eqref{data:hyp} and \eqref{smallness:condition:THM} of Theorem \ref{thm1}, with the same set of exponents $(s,p,q,m)$. Let $r \in  [1, + \infty)$. Then for all $t>0$, $\overline{q} \geq q$ and $k \in [0,\overline{q}-3)$, we have
\begin{align}\label{theo:estim-rhoLP}
\Vert \rho_f(t) \Vert_{\Ld^r(\R^3_+)} & \leq C_{k,\overline{q},r,g} \dfrac{K_{\overline{q},r}(f_0)+F_{0,k,r}(f_0)}{(1+t)^{k}},
 \end{align}
for some constant $C_{k,q,r,g}>0$.
\end{propo}
Our result also holds for any moment $m_{\ell}f$ of $f$ (see Definition \ref{notation:moments}).
\begin{propo}\label{PropoStatementMAIN:Moments-l}
Consider the same assumptions \eqref{data:hyp} and \eqref{smallness:condition:THM} of Theorem \ref{thm1}, with the same set of exponents $(s,p,q,m)$. Let $r \in  [1, + \infty)$ and $\ell>0$. Then the following holds:
\begin{itemize}
\item for all $t>0$, $\overline{q} \geq q$ and $k \in [0,\overline{q}-3)$, we have 
\begin{align*}
\Vert m_{\ell}f(t) \Vert_{\Ld^{\infty}(\R^3_+)} & \leq C_{k,\overline{q},g,\ell}  \dfrac{N_{\overline{q}}(f_0)+H_{\ell,k}(f_0)}{(1+t)^{k}},
\end{align*}
for some constant $C_{k,\overline{q},g, \ell}>0$;
\item  for all $t>0$, $\overline{q} \geq q$ and $k \in [0,\overline{q}-\ell-3)$ such that $\ell +3 < \overline{q}$, we have 
\begin{align*}
\Vert m_{\ell}f(t) \Vert_{\Ld^r(\R^3_+)} & \leq C_{k,\overline{q},r,g, \ell} \dfrac{K_{\overline{q},r}(f_0)+F_{\ell,k,r}(f_0)}{(1+t)^{k}},
 \end{align*}
 for some constant $C_{k,\overline{q},r,g, \ell}>0$.
\end{itemize}
\end{propo}

%\begin{rem}
%\textcolor{red}{One can obtain some decay in time of the $\Ld^r(\R^3_+)$ norm of $\rho_f(t)$: under the same conditions and with the same exponents $q$ and $k \in [0,q-3)$ as in the previous statement, one has for all $r \in [1, + \infty)$ and $t>0$
%\begin{align}\label{theo:estim-rhoLP}
%\Vert \rho_f(t) \Vert_{\Ld^r(\R^3_+)} & \lesssim \dfrac{K_{q,r}(f_0)+F_{0,k,r}(f_0)}{(1+t)^{k}}.
% \end{align}}
%
%\end{rem}

%\begin{rem}
%%The two last estimates can be refined:
%%\begin{align}
%%\Vert \rho_f(t) \Vert_{\Ld^{\infty}(\R^3_+)} & \lesssim  \dfrac{N_q(f_0)}{(1+t)^{k_1}}+ \dfrac{H_{0,k_2}(f_0)}{(1+t)^{k_2}}, \\
%%\Vert \rho_f(t) \Vert_{\Ld^r(\R^3_+)} & \lesssim \dfrac{K_{q,r}(f_0)}{(1+t)^{k_1}}+\dfrac{F_{0,k_2,r}(f_0)}{(1+t)^{k_2}}, \ \ r \in [1, + \infty).
%%\end{align}
%%where $(k_1,k_2) \in \R^+ \times \R^+$ is choosen such that $q>k_1+3$ and  with $\lesssim$ only depending  on $k_1,k_2,q$ and $g$. 
%\textcolor{red}{One can also replace the quantity $\rho_f$ in \eqref{theo:estim-rhofPONCT}--\eqref{theo:estim-rhoLP} by any moment $m_{\ell}f$, with $k+ \ell + 3<q$ (see Definition \ref{notation:moments}) and an estimate involving $N_{q}(f_0)+H_{\ell,k}(f_0)$ for a pointwise decay in $\Ld^{\infty}(\R^3_+)$, or involving $K_{q,r}(f_0)+F_{\ell,k,r}(f_0)$ for a pointwise decay in $\Ld^r(\R^3_+)$}.
%\end{rem}
The two previous propositions will be direct consequences of our proof of Theorem \ref{thm1}.

\medskip

\begin{rem}
As a byproduct of our analysis, we will prove along the way that if $f_0$ is compactly supported in velocity and in the third direction in space, that is if
$$\mathrm{supp} \, f_0 \subset (\R^2 \times (0,L)) \times \B(0,R),$$
for some finite $L,R,>0$, then there exists a finite time $T=T(L,R,g)$ such that $f(t) \equiv 0$ for every $t>T$.
\end{rem}

%\begin{thm}\label{thm2}
%Under the same assumptions of Theorem (\ref{thm1}) and if $\eqref{data:hyp}$ is complemented with $q>k+5$ for some $k>2$, then for any weak solution $(u,f)$ to the Vlasov-Navier-Stokes system with admissible initial data $(u_0,f_0)$, there exists $\rho^{\infty} \in \Ld^{\infty}(\R^3_+)$ such that
%\begin{align}
%&\rho_f(t) \overset{t \rightarrow +\infty}{\longrightarrow} \rho^{\infty} \ \ in \ \H^{-1}(\R^3_+),
%\end{align}
%Moreover, the last convergence also occurs with a polynomial rate $k/2-1$.
%\end{thm}

%\begin{rem}
%\textcolor{red}{We will also prove that the fluid part $u$ of any global Leray solution $(u,f)$ to the VNS system starting at an initial data satisfying \eqref{data:hyp}--\eqref{smallness:condition:THM} is actually the unique 'strong' solution (see \cite{Rob}) to the Navier-Stokes equations \eqref{eq:NS}--\eqref{eq:NS2}.  Of course, this global in time information is non trivial because of the coupling between $u$ and $f$.}
%\end{rem}

\begin{rem}
In view of the results of \cite{HKMM,HK}, one could hope for a monokinetic behavior of the distribution function in large time: in other words, $f(t)$ should concentrate in velocity to a Dirac mass supported at $G$. However, the estimate \eqref{theo:estim-rhoLP} shows in particular that $\Vert f(t) \Vert_{\Ld^1(\R^3_+ \times \R^3)} \underset{t \rightarrow + \infty}{\longrightarrow} 0$ so that the previous singular behavior does not occur. This is due to the absorption of the particles at the boundary, which is combined to the presence of the gravity force.
\end{rem}

\bigskip

Let us explain the main strategy that has been already devised and used to study the large time behavior of the VNS system in \cite{HKMM,HK,EHKM} (which are gravity-less cases). Roughly speaking, under some smallness assumption on the initial data that we will detail below, \cite{HKMM,HK,EHKM} have proven that the fluid velocity $u$ tends to a constant when $t \rightarrow + \infty$, while the distribution function converges towards a Dirac mass in velocity. The later weak convergence is in particular measured thanks to the $1$-Wasserstein distance on the phase space. All of these works heavily rely on the decay of a well-chosen energy functional which essentially controls the convergence of $u$ and $f$. The choice of such a functional may depend on the domain. As explained before, two principal spatial frameworks have been explored.

\medskip

$\bullet$ \underline{\textit{The case of bounded domains}}: as already mentioned in the introduction, Han-Kwan, Moussa and Moyano have tackled the large time behavior of global weak solutions to the system set on the torus $\T^3$ in \cite{HKMM}, while the case of bounded domains with absorption boundary condition has been studied in \cite{EHKM}. The main strategy is the following. In the torus, the quantity which plays a crucial role in the large time dynamics is the so-called \textit{modulated energy} introduced by Choi and Kwon in \cite{ChKw}. It is defined as
\begin{align*}
\mathcal{E}(t):= 
\dfrac{1}{2}\int_{\T^3 \times \R^3} f(t,x,v) \vert v-\langle j_f(t) \rangle\vert^2 \, \mathrm{d} x \, \mathrm{d} v
+  \dfrac{1}{2}\int_{\T^3} \vert u(t,x)- \langle  u(t) \rangle \vert^2 \, \mathrm{d} x
+ \dfrac{1}{4} \vert \langle j_f(t) \rangle - \langle u(t) \rangle \vert ^2, 
\end{align*}
where $\langle \cdot \rangle$ stands for the spatial average on $\T^3$. On a smooth bounded domain $\Omega \subset \R^3$, the key functional is the kinetic energy $\E_{\Omega}$ itself, defined as in \eqref{eq:Energy} where $\R^3_+$ is replaced by $\Omega$. In both cases, the decay of the energies $\mathcal{E}$ and $\E_{\Omega}$ is based on the following formal energy-dissipation identities
\begin{align}\label{IdentityIntro}
\dfrac{\mathrm{d}}{\mathrm{d}t}\mathcal{E}(t) + \mathrm{D}_{\T^3}(t)=0, \ \ \ \dfrac{\mathrm{d}}{\mathrm{d}t}\E_{\Omega}(t) + \mathrm{D}_{\Omega}(t) = 0,
\end{align} 
where $\D_{\T^3}$ or $\D_{\Omega}$ are defined as in \eqref{eq:Dissipation} where $\R^3_+$ is replaced by $\T^3$ or $\Omega$. Under the assumption that the global bound $\rho_f \in \Ld^{\infty}(\R^+;\Ld^{\infty}(Q))$ holds, where $Q=\T^3$ or $Q=\Omega$, one can show that an exponential decay of the energy is satisfied, namely
\begin{align}\label{decayT3intro}
\forall t \geq 0, \ \ \ \mathcal{E}(t)\leq C_{\lambda} e^{- \lambda t}\mathcal{E}(0), \ \ \  \E_{\Omega}(t)\leq C'_{\lambda} e^{- \lambda t} \E_{\Omega}(0),
\end{align}
for some constants $\lambda, C_{\lambda}, C'_{\lambda}>0$. This mainly comes from Poincaré(-Wirtinger) inequality. Then, a straightening change of variable in velocity shows that a sufficient condition for obtaining the previous bound on $\rho_f$ is an estimate on the Lipschitz seminorm of $u$, that is
%allows to reduce the obtention of the previous bound to an estimate on the Lipschitz semi norm of $u$, that is
\begin{align*}
\int_0^{\infty} \Vert \nabla u(s) \Vert_{\Ld^{\infty}(Q)} \, \mathrm{d}s \ll 1.
\end{align*}
A bootstrap procedure has to be performed in order to ensure such a global control: the main idea is to interpolate the pointwise conditional decay \eqref{decayT3intro}
with higher order parabolic regularity estimates for the fluid velocity. In short, the previous approach requires a smallness assumption on the initial data of the type
\begin{align*}
\mathcal{E}(0) + \Vert u_0 \Vert_{\H^{1}(Q)} \ll 1, \ \ \text{or} \ \ \E_{\Omega}(0) + \Vert u_0 \Vert_{\H^{1}(Q)} \ll 1,
\end{align*}
and one can even replace the previous norm for $u_0$ by $\Vert u_0 \Vert_{\H^{1/2}(\T^3)}$ in the torus case. Furthermore, one can describe the structure of the final spatial density, which depends on the whole evolution of the system.

\medskip

$\bullet$ \underline{\textit{The case of the whole space}}:  in this case studied in \cite{HK}, Han-Kwan has shown that the crucial functional to consider is the kinetic energy $\E_{\R^3}$ itself, defined as in \eqref{eq:Energy} where $\R^3_+$ is replaced by $\R^3$. This energy satisfies the same formal energy identity \eqref{IdentityIntro} as above. A bound of the type $\rho_f \in \Ld^{\infty}(\R^+;\Ld^{\infty}(\R^3))$ now provides a decay of the form
\begin{align*}
\forall t \geq 0, \ \ \ \mathrm{E}_{\R^3}(t) \leq \dfrac{\varphi_{\alpha}\big(\mathrm{E}_{\R^3}(0)+ \Vert u_0 \Vert_{\Ld^1(\R^3)}\big)}{(1+t)^{\alpha}}, \ \ \text{for all} \ \ \alpha \in ]0,3/2[,
\end{align*}
for some function $\varphi_{\alpha}$. Here, the polynomial decay is the best one can hope for because of the absence of a Poincaré inequality (with respect to the Lebesgue measure) on domains unbounded in any direction. In short, it corresponds to the decay of the solutions to the heat equation on the whole space, with almost the same rate. As before, the second step of the analysis is a bootstrap analysis which aim is to obtain the control on $\rho_f$, relying on the same sufficient control of  $\Vert \nabla u \Vert_{\Ld^1(\R^+; \Ld^{\infty}(\R^3))}$. Since the energy decay is only polynomial and since the Brinkman force is not decaying, the use of dissipation functionals of higher order is needed\footnote{This family of identities has also found a powerful application in the study of the hydrodynamic limits of the VNS system in the torus performed by Han-Kwan and Michel in \cite{HKM}.}. It essentially leads to the study of weighted in time estimates for the second order derivatives (in space) of the fluid velocity. Here, this procedure is applied for small initial data, in the sense that 
 \begin{align*}
\E_{\R^3}(0) + \Vert u_0 \Vert_{\H^{1}(\R^3)} + \Vert f_0 \Vert_{\Ld^1_v(\R^3; \Ld^{\infty}_x(\R^3))} \ll 1.
\end{align*}

\medskip

\medskip
%Let us emphasize the fact that a particular geometry combined with specific boundary conditions has been considered by Glass, Han-Kwan and Moussa in \cite{GHKM}, and leads to a somewhat other type of asymptotic behavior. More precisely, the VNS system is set on a bidimensional rectangle $\Omega:=(-L,L) \times (-1,1)$: the fluid velocity satisfies a Dirichlet boundary condition corresponding to a Poiseuille flow while the distribution function obeys to absorption boundary conditions on the horizontal parts and to an injection boundary condition on the vertical left part. It has been shown that this particular framework allows to get existence and asymptotic stability of non-trivial smooth equilibria for the system. This means that the monokinetic behavior highlighted above does not occur. The strategy of proof heavily relies on the use of a key geometric control condition (the so-called \textit{exit geometric condition}) satisfied by the Poiseuille flow, and which compells the particle trajectory to leave the domain transversally before a fixed finite time. Here, the absorption phenomenon is thus used to avoid any concentration phenomena for the distribution function. 
Let us emphasize the fact that the case of a bounded domain with specific boundary conditions has been considered by Glass, Han-Kwan and Moussa in \cite{GHKM}. More precisely, the VNS system is set on a bidimensional rectangle $\Omega:=(-L,L) \times (-1,1)$: the fluid velocity satisfies a Dirichlet boundary condition corresponding to a Poiseuille flow while the distribution function obeys to mixed absorption/injection boundary conditions. Compared to the previous state of the art, this particular framework leads to a somewhat other type of asymptotic behavior. Indeed, relying on a key geometric control condition (the so-called \textit{exit geometric condition}), it has been shown that one can get the existence and asymptotic stability of non-trivial smooth equilibria for the system. 
%This means that the monokinetic behavior highlighted above does not occur. The strategy of proof heavily relies on the use of a key geometric control condition (the so-called \textit{exit geometric condition}) satisfied by the Poiseuille flow, and which compells the particle trajectory to leave the domain transversally before a fixed finite time. Here, the absorption phenomenon is thus used to avoid any concentration phenomena for the distribution function. 

\bigskip

\textbf{Main strategy}. \ Let us explain the main approach used in this article. As mentioned before, the study of the large time dynamics of the system \eqref{eq:Vlasov}-\eqref{eq:NS}-\eqref{eq:NS2} is in the same spirit as that of the previous works \cite{HKMM,HK,EHKM}. However, a main obstacle comes from the presence of the additional gravity force term in the Vlasov equation \eqref{eq:Vlasov}, which creates an extra term in the right-hand side of the energy inequality \eqref{ineq-energy}. This breaks one of the main structural tools of the analysis of the gravity-less case because it rules out the decay of the total kinetic energy $\E$. We thus need to base our study upon an additional mechanism which is at stake in the system. What comes into play here is the absorption of the particles at the boundary, coming from the boundary condition \eqref{bcond-f} for the distribution function $f$. To understand the crucial role of this phenomenon, we will use the Lagrangian structure of the Vlasov equation. We shall define the characteristic curves $(\X,\V)$ for the Vlasov equation as the solutions of the following differential system:
\begin{equation*}
\left\{
      \begin{aligned}
        \frac{\mathrm{d}}{\mathrm{d}s}\mathrm{X}(s;t,x,v) &=\mathrm{V}(s;t,x,v),\\
\frac{\mathrm{d}}{\mathrm{d}s}\mathrm{V} (s;t,x,v)&= u(s,\mathrm{X}(s;t,x,v))+G-\mathrm{V}(s;t,x,v),
      \end{aligned}
    \right.
\end{equation*}
with $(\mathrm{X}(t;t,x,v),\mathrm{V}(t;t,x,v))=(x,v)$ and where $u$ has been extended by $0$ outside the half-space. Introducing
$$\mathcal{O}^t=\left\lbrace (x,v) \in \R^3_+ \times \R^3 \mid \forall \sigma \in [0,t], \ \  \mathrm{X}(\sigma;t,x,v) \in \R^3_+ \right\rbrace,
$$
the method of characteristics shall provide the following representation formula:
\begin{align}\label{rep:formulaINTRO}
 f(t,x,v)= e^{3t} \mathbf{1}_{\mathcal{O}^t}(x,v) \, f_0\big(\mathrm{X}(0;t,x,v),\mathrm{V}(0;t,x,v)\big). 
\end{align}
In view of this expression, a certain decay in time of the moments of $f$ should be satisfied along the evolution of the system, provided that $f_0$ enjoys some decay in the phase space.

In order to take advantage of the absorption, we shall rely on an \textit{exit geometric condition}, reminiscent of the work of Glass, Han-Kwan and Moussa in \cite{GHKM}. In short, we ask that all the characteristic curves starting from a compact set leave the half-space before a fixed time. The main idea to propagate this condition will be to compare the coupled Vlasov equation to the Vlasov equation without fluid velocity and only governed by the gravity. The characteristic curves $(\mathrm{X}^g,\mathrm{V}^g)$ for this simplified equation are defined as the solutions of
\begin{equation*}
\left\{
      \begin{aligned}
        \frac{\mathrm{d}}{\mathrm{d}s}\mathrm{X}^g(s;t,x,v) &=\mathrm{V}^g(s;t,x,v),\\
\frac{\mathrm{d}}{\mathrm{d}s}\mathrm{V}^g (s;t,x,v)&= G-\mathrm{V}^g(s;t,x,v),
      \end{aligned}
    \right.
\end{equation*}
with $(\mathrm{X}(t;t,x,v),\mathrm{V}(t;t,x,v))=(x,v)$. We will then show that the exit geometric condition holds for all times for the VNS system, by using the particular geometry of the domain and the simple form of $(\mathrm{X}^g,\mathrm{V}^g)$.  This will essentially require a control of the form
\begin{align}\label{strat:Bound2}
\int_0^{\infty} \Vert u(s) \Vert_{\Ld^{\infty}(\R^3_+)} \, \mathrm{d}s \ll 1.
\end{align}

By combining the previous absorption phenomenon with a decay of the initial distribution function $f_0$ itself, we shall be able to obtain decay in time estimates of the moments of $f$. The argument will be based on the representation formula \eqref{rep:formulaINTRO} and on a change of variable in velocity, namely $v \mapsto \V(0;t,x,v)$ (which was already used in \cite{HKMM,HK,EHKM}). This procedure will be allowed if we can ensure a control of the Lipschitz seminorm of the fluid velocity, that is 
\begin{align}\label{strat:Bound1}
\int_0^{\infty} \Vert \nabla u(s) \Vert_{\Ld^{\infty}(\R^3_+)} \, \mathrm{d}s \ll 1.
\end{align}
Thus, it turns out that the presence of the gravity force is eventually favorable for our purpose and leads to the decay estimates of the moments stated in Theorem \ref{thm1} and Propositions \ref{PropoStatementMAIN:estimLP} and \ref{PropoStatementMAIN:Moments-l}.

Thanks to the specific form of the Brinkman force $j_f-\rho_f u$ in the Navier-Stokes equations, we shall obtain pointwise in time estimates in various norms for this source term. Hence, as a solution of the Navier-Stokes equations with a sufficiently decaying forcing term, the fluid velocity $u$ will enjoy a polynomial convergence towards $0$. This is essentially the result bearing on $u$ in Theorem \ref{thm1}.

Therefore, our main guiding line will be the obtention of decay in times estimates for the moments of $f$. We will base our proof on a bootstrap argument, mainly directly taken from \cite{HK}.

\bigskip

\textbf{Outline of the paper}. \ According to the previous strategy, let us describe how this paper is organised. 

\begin{itemize}
\item[$-$] In view of the arguments above, it makes sense to first consider the Navier-Stokes system having a source term (i.e. the Brinkman force \eqref{eq:Brinkman}) with polynomial decay in time. It turns out that this assumption falls within the scope of the work of Wiegner in \cite{Wieg} for the Navier-Stokes equations with a decaying source term on $\R^3$. \textit{Modulo} an adaptation to the half-space case, this entails a polynomial decay of the $\Ld^2$ norm of the fluid velocity. This conditional Theorem \ref{cond:decay} is contained in \textbf{Section \ref{Section:CondDecay}}. 
\item[$-$] Before going further, we shall need to state some preliminaries gathered in \textbf{Section \ref{Section:Prelim}}. They are necessary for a local in time analysis as well as for the subsequent bootstrap argument. We obtain rough bounds on the moments of $f$, ensuring short time controls. We also derive some $\H^1$ estimates for the fluid velocity and define the notion of \textit{strong existence time} for the Navier-Stokes system.
\item[$-$] As explained before, the absorption of the particles at the boundary will be the key effect leading to global decay in time for the moments of $f$. In \textbf{Section \ref{Section:EGCabs}}, we introduce the aforementioned crucial \textit{exit geometric condition} coming from \cite{GHKM} and analyse its effect on the system. In particular, this enables us to track which proportion of the support of the initial distribution has disappeared from the system at any given time.
\item[$-$] \textbf{Section \ref{Section:Bootstrap}} is devoted to the bootstrap argument, which aims to achieve the global controls \eqref{strat:Bound2} and \eqref{strat:Bound1}. Along the bootstrap, the previous absorption phenomenon is shown to lead to the desired decay estimates of the moments. We first show that the Brinkman force satisfies a suitable pointwise decay. Polynomial weighted in time estimates for the fluid velocity are then obtained. This will allow to close the bootstrap argument thanks to an interpolation procedure. 
\end{itemize}

In the rest of the article, we will use the standard notation $A \lesssim B$ for $A \leq c B$ for some $c>0$ which is independent of $A,B$ and that may change from line to line.

\section{Conditional large time behavior of the fluid velocity}\label{Section:CondDecay}
The main goal of this short section is to show some conditional results about the polynomial decay of the $\Ld^2$ norm of any Leray solution $u$ to the Navier-Stokes equation with a source term $F=F(t,x)$, that is
\begin{equation}\label{systNS:decay}
\left\{
      \begin{aligned}
        \partial_t u +(u\cdot \nabla_x )u-\Delta_x u + \nabla_x p&=F,\\
        \mathrm{div}_x \, u&=0, \\
        u_{\mid x_3=0}&=0, \\
        u_{\mid t=0}&=u_0.
      \end{aligned}
    \right.
\end{equation}

We shall require that $u$ satisfies the strong energy inequality, that is for any $t \geq 0$ and almost every $0 \leq s  \leq t$ (including $s=0$), we have
\begin{align}\label{ineq-energyNS-decay}
\Vert u(t) \Vert_{\Ld^2(\R^3_+)}^2  + 2 \int_s ^t \Vert \nabla u(\tau) \Vert_{\Ld^2(\R^3_+)}^2 \mathrm{d}\tau \leq   \Vert u(s) \Vert^2_{\Ld^2(\R^3_+)} + 2 \int_s ^t \int_{\R^3_+} F(\tau,x) \cdot u(\tau,x) \, \mathrm{d}x \,  \mathrm{d}\tau.
\end{align}
The decay of such a solution will hold if $F$ satisfies a conditional pointwise decay in $\Ld^2(\R^3_+)$ and is somewhat imposed by the decay of the Stokes semigroup on $\R^3_+$. The main result reads as follows.
%\begin{itemize}
%\item the polynomial decay of the $\Ld^2$ norm of the fluid velocity, which is somewhat imposed by the decay of the Stokes semigroup on $\R^3_+$ and which holds if the Brinkman force $j_f-\rho_f u$ satisfies a conditional pointwise decay in $\Ld^2(\R^3_+)$;
%\item the weak convergence of the distribution function $f$ which reveals a monokinetic behavior and which is measured thanks to the Wasserstein distance on the phase-space, with a control by the second moment in velocity of $f$.
%\end{itemize}
%We first focus on the conditional large time behavior of the $\Ld^2$ norm of the fluid velocity.
\begin{thm}\label{cond:decay}
Let $u_0 \in \Ld^2_{\mathrm{div}}(\R^3_+) \cap \Ld^1(\R^3_+)$ and $F \in \Ld^{\infty}_{\mathrm{loc}}(\R^+; \Ld^2(\R^3_+))$. Let $u$ be a global Leray solution to the Navier-Stokes system \eqref{systNS:decay} with strong energy inequality \eqref{ineq-energyNS-decay}, associated to the initial data $u_0$ and the source term $F$.
Let $T>0$ and assume that
\begin{align}\label{HYPdecay:BRINK}
\forall t \in [0,T], \ \ \Vert F(t) \Vert_{\Ld^2(\R^3_+)} \leq \dfrac{\mathrm{C}}{(1+t)^{7/4}},
\end{align}
for some constant $\mathrm{C}>0$ independent of $T$. Then there exists a continuous nonnegative function $\Psi$ cancelling at $0$ and independent of $T$
%(and increasing in its parameter $\mathrm{C}$),
such that
\begin{align}
\label{eq:decaythmcond}
\forall t \in [0,T], \ \ \Vert u(t) \Vert_{\Ld^2(\R^3_+)}^2 \leq \frac{\Psi\left (\Vert  u_0 \Vert_{\Ld^{2}(\R^3_+)}^2+\Vert  u_0 \Vert_{\Ld^{1}(\R^3_+)}^2+ \mathrm{C} \right)}{(1+t)^{3/2}}.
\end{align}
\end{thm}
Here, we will follow Wiegner \cite{Wieg} and Borchers and Miyakawa \cite{BM}, relying on the Fourier splitting method of Schonbeck \cite{Schon}. More precisely, assuming a decay of the type \eqref{HYPdecay:BRINK} for the source term $F$ means that the desired result for the large time behavior of the fluid velocity enters in the framework of \cite{Wieg}. Since we work in an unbounded domain with boundaries, we shall adapt this method written in the whole space case, thanks to a spectral decomposition of the Stokes operator and the use of its fractional powers, as in \cite{BM}. This will entail a polynomial decay of the fluid velocity similar to the one without a source term in the equations, and whose rate is roughly speaking the same as that of the unsteady Stokes equations. 

Coming back to the Vlasov-Navier-Stokes system, we shall consider the Brinkman force $j_f-\rho_f u$ as a fixed source term in the Navier-Stokes equations. This means that we shall use Theorem \ref{cond:decay} with $$F=j_f-\rho_f u.$$ Note that compared to \cite{HK} where the conditional decay of the energy $\E$ was related to the whole VNS system, our result concerns the decay of $u$ as a solution of the Navier-Stokes equations only and is independent of the coupling with the Vlasov equation: we only use the strong energy inequality \eqref{ineq-energyNS} for the Navier-Stokes system with a given source term (even if it may depend on $u$ and $f$). Therefore, the assumption \eqref{HYPdecay:BRINK} on the decay of this  source term makes the situation simpler. We will be able to prove that this strong decay does occur thanks to the absorption phenomenon along a bootstrap procedure in Section \ref{Section:Bootstrap}. Here, we do not require a bound of the type $\rho_f \in \Ld^{\infty}(0,T;\Ld^{\infty}(\R^3_+))$ as in \cite{HK} and the rate of convergence is slightly better.
%we assume that the Brinkman force is a forcing term with a sufficient decay in the Navier-Stokes equations and that we will be able to prove that this strong decay indeed occurs thanks to the absorption phenomenon along a bootstrap procedure in Section \ref{Section:Bootstrap}. This assumption entails a polynomial decay of the fluid velocity similar to the one without a source term, and whose rate is roughly speaking the same as that of the unsteady Stokes equations.

The combination of \cite{Wieg} and \cite{BM} for the proof of Theorem \ref{cond:decay} may appear as a classic result for the Navier-Stokes system: for the reader's convenience, we only write the proof in Section \ref{Section:CondDecayAPPENDIX} of the Appendix.
In view of the conditional Theorem \ref{cond:decay}, obtaining decay in time for the moments of $f$ will be the main goal of the rest of this paper.  
\section{Preliminaries for the bootstrap procedure}\label{Section:Prelim}
Thanks to a bootstrap argument, we will prove that there exists $\mathrm{M}>0$ such that 
\begin{align*}
%\Vert \rho_f \Vert_{\Ld^{\infty}(\R^+;\Ld^{\infty}(\R^3_+))} &\leq \mathrm{C},\\
\forall t \in \R^+, \ \ \Vert j_f(t)-\rho_f u(t) \Vert_{\Ld^2(\R^3_+)} &\leq \dfrac{\mathrm{M}}{(1+t)^{7/4}}.
%\forall t \in \R^+, \ \ \M_2 f(t) &\lesssim \dfrac{\mathrm{C}}{(1+t)^k}, \ \text{ for some} \ k>0
\end{align*}
Then, by the conditional Theorem \ref{cond:decay}, this will imply the first statement of Theorem \ref{thm1}. Along the way, we shall rely on the following bounds and estimates for the local density $\rho_f$:
\begin{align*}
\Vert \rho_f \Vert_{\Ld^{\infty}(\R^+;\Ld^{\infty}(\R^3_+))} & \lesssim 1,\\
\forall t \in \R^+, \forall r \geq 1, \ \ \Vert \rho_f (t) \Vert_{\Ld^r(\R^3_+)}  &\lesssim \dfrac{1}{(1+t)^k}, \ \text{ for some} \ k>0,
\end{align*}
which will essentially lead to the second part of Theorem \ref{thm1}.

\medskip

In this section, we collect several useful information in order to be able to set up a boostrap procedure in Section \ref{Section:Bootstrap}. We start by recalling some basic facts about the Lagrangian structure of the Vlasov equation \eqref{eq:Vlasov}. A careful analysis of the characteristic curves will indeed be required to deal with the absorption at the boundary. It also enables us to consider a straightening change of variable in velocity.
Then, we derive local in time estimates for the moments $\rho_f$ and $j_f$, as well as for the fluid velocity in $\Ld^1 \Ld^{\infty}$. This will offer short time controls on these quantities. Thanks to a smoothing property of the Navier-Stokes system, we finally obtain $\Ld^2 \H^2 \cap \Ld^{\infty} \H^1$ estimates for the fluid velocity. This requires the introduction of the so-called \textit{strong existence times} and eventually entails a local in time $\Ld^1 \W^{1,\infty}$ regularity for the fluid velocity.
%which will rely on the boundary effects.
\subsection{Characteristic curves for the Vlasov equation: representation formula and change of variable}
Given a time-dependent vector field $u$ on $\R^+ \times \R^3_+$, a time $t \in \R^+$ and a point $(x,v) \in \R^3_+ \times \R^3$, we define the characteristic curves $s \in \R^+ \mapsto (\mathrm{X}(s;t,x,v), \mathrm{V}(s;t,x,v)) \in \R^3 \times \R^3$ for the Vlasov equation (associated to $u$) as the solution of the following system of ordinary differential equations
\begin{equation}\label{EDO-charac}
\left\{
      \begin{aligned}
        \frac{\mathrm{d}}{\mathrm{d}s}\mathrm{X}(s;t,x,v) &=\mathrm{V}(s;t,x,v),\\
\frac{\mathrm{d}}{\mathrm{d}s}\mathrm{V} (s;t,x,v)&= (Pu)(s,\mathrm{X}(s;t,x,v))+G-\mathrm{V}(s;t,x,v),\\
	\mathrm{X}(t;t,x,v)&=x,\\
	\mathrm{V}(t;t,x,v)&=v.
      \end{aligned}
    \right.
\end{equation}
Here, $P$ is the linear extension operator continuous from $\Ld^{\infty}(\R^3_+)$ to $\Ld^{\infty}(\R^3)$ and from $\W^{1,\infty}_0(\R^3_+)$ to $\W^{1,\infty}(\R^3)$ defined by 
\begin{equation}
\forall x \in \R^d, \ \ (Pw)(x):=\left\{
      \begin{aligned}
        & w(x) \ &&\text{if} \ x \in \R^3_+, \\
	& 0\ &&\text{if} \ x \in \R^3 \setminus \R^3_+,
      \end{aligned}
    \right.
\end{equation} and which satisfies
\begin{align}
& \forall w \in \Ld^{\infty}(\R^3_+), \ \Vert Pw\Vert_{\Ld^{\infty}(\R^3)} = \Vert w \Vert_{\Ld^{\infty}(\R^3_+)}, \label{opP:Linfini} \\[2mm]
& \forall w \in \W^{1,\infty}_0(\R^3_+), \ \ \Vert \nabla(Pw)\Vert_{\Ld^{\infty}(\R^3)} \leq \Vert \nabla w \Vert_{\Ld^{\infty}(\R^3_+)}. \label{opP:grad}
\end{align}
We refer to \cite[Appendix]{EHKM} for a simple proof. Also, we will use the convention
\begin{align*}
(Pu)(t,\cdot)=P(u(t,\cdot)),
\end{align*}
as well as the notation 
\begin{align*}
\X^s_t(x,v):=\X(s;t,x,v), \ \ \V^s_t(x,v):=\V(s;t,x,v).
\end{align*}

\medskip

Let $T>0$ be fixed and suppose $$u \in \Ld^2(0,T;\H^1_0(\R^3_+)) \cap \Ld^1(0,T;\W^{1,\infty}(\R^3_+)).$$  We can apply the Cauchy-Lipschitz theorem to show the following proposition.
\begin{propo}\label{Propo:diffeoZ}
Given $(x,v) \in \R^3 \times \R^3$ and a time $t \in [0,T]$, the system (\ref{EDO-charac}) admits a unique solution $s \mapsto \mathrm{Z}^s_{t}(x,v) \in \R^3 \times \R^3$ on $[0,T]$ and
\begin{equation*}
\mathrm{Z}^s_{t} : \left\{
   \begin{aligned}
      & \R^3 \times \R^3 &&\longrightarrow \R^3 \times \R^3\\
      & \ (x,v) &&\longmapsto \mathrm{Z}^s_{t}(x,v):= (\mathrm{X}(s;t,x,v),\mathrm{V}(s;t,x,v))
      \end{aligned}
    \right.
\end{equation*}
is a (bi-Lipschitz) diffeomorphism of $\R^3 \times \R^3$ whose inverse is given by $(\mathrm{Z}^s_{t})^{-1}=\mathrm{Z}^t_{s}$ and whose Jacobian determinant is $e^{3(t-s)}$.
\end{propo}
In this context, the characteristic curves for the Vlasov equation are classically defined (at least) until time $T$ and are given for all $s,t \in [0,T]$ by
\begin{equation}\label{expr:Zt}
\left\{
      \begin{aligned}
        \mathrm{X}(s;t,x,v)&=x+(1-e^{-s+t})v+(s-t+e^{-s+t}-1)G+\int_t^s (1-e^{\tau-s}) (Pu)(\tau,\mathrm{X}(\tau;t,x,v))  \, \mathrm{d}\tau,  \\
        \mathrm{V}(s;t,x,v)&= e^{-s+t}v+(1-e^{-s+t})G+\int_t ^s e^{\tau-s}(Pu)(\tau,\mathrm{X}(\tau;t,x,v))   \, \mathrm{d}\tau.   \\
      \end{aligned}
    \right.
\end{equation}
Starting from a point $(x,v) \in \R^3_+ \times \R^3$ at time $t$, the curve $\X(s;t,x,v)$ remains during a certain interval of time in the half-space $\R^3_+$. This naturally leads to the following definitions, already considered in \cite{EHKM}.
\begin{defi}\label{def:tau}
For $(x,v) \in \R^3_+ \times \R^3$ and for any $t \geq 0$, we set
\begin{align}
\tau^{+}(t,x,v)&:=\sup \left\lbrace s \geq t \ \mid \forall \sigma \in [t,s], \  \mathrm{X}(\sigma;t,x,v) \in \R^3_+ \right \rbrace. \label{def:tau+}
%\tau^{-}(t,x,v)&:=\inf \left\lbrace s \leq t \ \mid \forall \sigma \in [s,t], \  \mathrm{X}(\sigma;t,x,v) \in \R^3_+ \right \rbrace, \label{def:tau-}
\end{align}
We also define
\begin{align}\label{def:Ot}
\mathcal{O}^t=\left\lbrace (x,v) \in \R^3_+ \times \R^3 \mid \forall \sigma \in [0,t], \ \  \mathrm{X}(\sigma;t,x,v) \in \R^3_+ \right\rbrace.
\end{align}
\end{defi}

We state two basic results whose proof can be found in \cite[Appendix]{EHKM}. The second one is a representation formula for the weak solution to the Vlasov equation \eqref{eq:Vlasov}, where one has to take into account the absorption boundary condition \eqref{bcond-f} satisfied by the distribution function.
\begin{lem}\label{LM:tau+-}
For $z=(x,v) \in \R^3_+ \times \R^3$, if $\tau^{+}(0,z)< +\infty$ then we have
\begin{align}\label{traj:exit}
\mathrm{Z}^{\tau^{+}(0,z)}_{0}(z)=\left(\mathrm{X}^{\tau^{+}(0,z)}_{0}(z),\mathrm{V}^{\tau^{+}(0,z)}_{0}(z) \right) \in \Sigma^+ \cup \Sigma^0,
\end{align}
where $\Sigma^+$ and $ \Sigma^0$ are defined in \eqref{def:Sigma}. Furthermore, we have for all $t \geq 0$
\begin{align*}
\mathrm{Z}^0_t(\mathcal{O}^t)=\left\lbrace  (x,v) \in \R^3_+ \times \R^3  \mid  \tau^{+}(0,x,v) >t  \right\rbrace.
\end{align*}
\end{lem}

\begin{propo}\label{Prop:formulerep}
Let $f$ be the weak solution to the Vlasov equation
\begin{equation*}
\left\{
      \begin{aligned}
      \partial_t f +v\cdot \nabla_x f + \mathrm{div}_v \, (f(u-v)+fG)&=0,\\
f_{\mid t=0}&=f_0,\\
f &=0 \ \mathrm{on} \ \Sigma^{-},
      \end{aligned}
    \right.
\end{equation*}
associated to a velocity field $u \in  \Ld^2_{\mathrm{loc}}(\R^+;\H^1_0(\R^3_+)) \cap \Ld^1_{\mathrm{loc}}(\R^+;\W^{1,\infty}(\R^3_+))$ with initial condition $f_0$ and with absorption boundary condition. There holds
\begin{align}\label{formule-rep}
 f(t,x,v)= e^{3t} \mathbf{1}_{\mathcal{O}^t}(x,v) \, f_0(\mathrm{Z}^0_{t}(x,v)) \ \ \text{a.e.}
\end{align}
\end{propo}

Recall that we aim at obtaining a sufficient decay in time of the Brinkman force $j_f-\rho_f u$ and more generally of the moments of $f$. The representation formula \eqref{formule-rep} will be our main starting point: elaborating on the same strategy as that of \cite{HKMM,HK,EHKM,HKM}, we shall rely on a straightening change of variable in velocity (and then in space) in this formula. Nevertheless, such a procedure requires a smallness assumption on the quantity $\Vert \nabla u \Vert_{\Ld^{1}(\R^+;\Ld^{\infty}(\R^3_+))}$ and obtaining this control will be at the core of Section \ref{Section:Bootstrap}. We emphasize the fact that we also need the help of the absorption at the boundary in order to recover the desired decay in time of the moments.
%Therefore the previous change of variable will only be used in Section \ref{Section:Bootstrap} (see in particular Lemma \ref{LM:decaymom-abs}).
%Elaborating on the same strategy than that of \cite{HKMM,HK,EHKM,HKM}, we provide a sufficient condition ensuring global bounds on the moments, relying on a straightening change of variable in velocity. Nevertheless, such a procedure requires a smallness assumption on the quantity $\Vert \nabla u \Vert_{\Ld^{1}(\R^+;\Ld^{\infty}(\R^3_+))}$ and obtaining this control will be at the core of Section \ref{Section:Bootstrap}. We emphasize the fact that these changes of variable shall entail decay estimates of the moments of $f$ (see in particular Lemma \ref{LM:decaymom-abs}).

\medskip

In view of the following formulas
\begin{align*}
\V(0;t,x,v)&=e^{t}v+(1-e^t)G-\int_0^t e^{\tau}(Pu)(\tau,\mathrm{X}(\tau;t,x,v))   \, \mathrm{d}\tau, \\
\mathrm{D}_v \V(0;t,x,v)-e^t \mathrm{Id}&=-\int_0^t e^{\tau} \nabla(Pu)(\tau,\mathrm{X}(\tau;t,x,v)) \mathrm{D}_v \mathrm{X}(\tau;t,x,v))    \, \mathrm{d}\tau,
\end{align*}
 and following closely the arguments of \cite{HKMM}, we infer several statements which read as follows.
\begin{lem}\label{chgmt-var-prepa}
Suppose $u \in \Ld^2_{\mathrm{loc}}(\R^+,\H^1_0(\R^3_+)) \cap\Ld^1_{\mathrm{loc}}(\R^+;\Ld^{\infty}(\R^3_+))$. Fix $\delta >0$ satisfying $\delta e^{\delta} < 1/9$. Then, for all times $t \in \R^{+}$ satisfiying
\begin{align}\label{borne-gradient0}
\int_0 ^t \Vert \nabla u(s) \Vert_{\Ld^{\infty}(\R^3_+)}  \, \mathrm{d}s < \delta,
\end{align}
and for all $x \in \R^3_+$, the map  
$$ \Gamma_{t,x} : v \mapsto \mathrm{V}(0;t,x,v),$$ is a global $\mathscr{C}^{1}$-diffeomorphism from $\R^3$ to itself satisfying furthermore 
\begin{align}\label{control:jacob}
\forall v \in \R^3, \ \ \vert \det \, \mathrm{D}_v \Gamma_{t,x} (v)  \vert \geq \dfrac{e^{3t}}{2}.
\end{align} 
\end{lem}

Thanks to 
\begin{align}
\vert \Gamma_{t,x}^{-1}(w) \vert &\leq e^{-t} \left[ \vert w\vert +(e^t-1) \vert G \vert + \int_0^t e^{\tau} \Vert Pu(\tau) \Vert_{\Ld^{\infty}(\R^3)} \, \mathrm{d}\tau \right], \label{ineqGamma-1} \\
\vert v \vert  & \leq \vert \V(0;t,x,v) \vert +(1-e^{-t}) \vert G \vert + \int_0^t \Vert Pu(\tau) \Vert_{\Ld^{\infty}(\R^3)} \, \mathrm{d}\tau, \label{ineqV^0_t}
\end{align}
we also have the following result.
%\begin{propo}
%Suppose $u \in \Ld^2_{\mathrm{loc}}(\R^+,\H^1_0(\R^3_+)) \cap\Ld^1_{\mathrm{loc}}(\R^+;\Ld^{\infty}(\R^3_+))$. If the assumption (\ref{borne-gradient0}) is satisfied at a time $t\geq 0$, then we have
%\begin{align*}
%\Vert \rho_f(t) \Vert_{\Ld^{\infty}(\R^3_+)}  & \leq 2 I_q N_q(f_0), \\
%\Vert j_f(t) \Vert_{\Ld^{\infty}(\R^3_+)} & \leq 2 I_q N_q(f_0) \left[ \vert G \vert + e^{-t} \left(\int_0 ^t e^s \Vert u(s) \Vert_{\Ld^{\infty}(\R^3_+)}   \, \mathrm{d}s +1 \right) \right] ,
%\end{align*}
%where we recall that $N_q(f_0)$ has been defined in Definition \ref{def:Nq} and where
%$$I_q:= \int_{\R^3} \dfrac{1+\vert v\vert}{1+\vert v \vert ^q}  \, \mathrm{d}v.$$
%\end{propo}
\begin{lem}\label{inegdecal}
Let $t_0 >0$ and $q>0$. If $N_q(f_0) < \infty$ and if $u \in \mathrm{L}_{\mathrm{loc}}^1(\R_+ ; \mathrm{H}^1_0 \cap \mathrm{L}^{\infty}(\R^3_+))$ then $N_q(f(t_0)) < \infty $ with 
\begin{align*}
N_q(f(t_0))  \lesssim e^{3 t_0} (1+ \vert G \vert^q+ \Vert u \Vert_{\mathrm{L}^1 (0,t_0 ; \mathrm{L}^{\infty}(\R^3_+))} ^q) N_q(f_0).
\end{align*}
\end{lem}

%\begin{lem}\label{lastLMgradient}
%Let $t_0 \geq 0$. Let $u \in \Ld^2_{\mathrm{loc}}(\R^+,\H^1_0(\R^3_+)) \cap\Ld^1_{\mathrm{loc}}(\R^+;\Ld^{\infty}(\R^3_+))$. Let $\delta$ be fixed such that $\delta e^{\delta} < 1/9$. For all times $t \geq t_0 \geq0$ such that
%\begin{align}\label{borne-gradient}
%\int_{t_0} ^t \Vert \nabla u(s) \Vert_{\Ld^{\infty}(\R^3_+)}  \, \mathrm{d}s \leq \delta,
%\end{align}
%we have
%\begin{align*}
%\Vert \rho_f(t) \Vert_{\Ld^{\infty}(\R^3_+)}  & \lesssim e^{3t_0}N_q(f_0) (1+\vert G \vert^q+  \Vert u \Vert_{\mathrm{L}^1 (0,t_0 ; \mathrm{L}^{\infty}(\R^3_+))} ^q).
%\Vert j_f(t) \Vert_{\Ld^{\infty}(\R^3_+)} & \lesssim  e^{3t_0-t} \left(\int_{t_0} ^t e^s \Vert u(s) \Vert_{\Ld^{\infty}(\R^3_+)}   \, \mathrm{d}s +1 \right) N_q(f_0) (1+\vert G \vert^q+  \Vert u \Vert_{\mathrm{L}^1 (0,t_0 ; \mathrm{L}^{\infty}(\R^3_+))} ^q).
%\end{align*}
%\end{lem}
As we shall see later, we will also perform a change of variable in space after the previous change of variable in velocity.
\begin{lem}\label{chgmt:varSPACE}
Consider the same assumptions as in Lemma \ref{chgmt-var-prepa}. For any $t \in \R^+$ satisfying
\begin{align*}
\int_{0} ^t \Vert \nabla u(s) \Vert_{\Ld^{\infty}(\R^3_+)}  \, \mathrm{d}s \leq \delta,
\end{align*}
%for any $(x,v) \in \R^3 \times \R^3$, we define
%\begin{align}
%\widetilde{\X}^{\cdot}_t(x,v) : s \longrightarrow \widetilde{\X}^{s}_t(x,v) :=\X^{s}_t(x, \Gamma_{t,x}^{-1}(v)).
%\end{align}
and for any $v \in \R^3$, the map
\begin{align}
\Lambda_{t,v} : x \longrightarrow \X^{0}_t(x, \Gamma_{t,x}^{-1}(v))
\end{align}
 is a global $\mathscr{C}^{1}$-diffeomorphism from $\R^3$ to itself satisfying 
\begin{align}
\forall x \in \R^3, \ \ \vert \det \, \mathrm{D}_x \Lambda_{t,v} (x)  \vert \geq \dfrac{1}{2}.
\end{align} 
\end{lem}
\begin{proof}
From
\begin{align*}
 \Gamma_{t,x}^{-1}(w) = e^{-t}  w +(1-e^t) G+ \int_0^t e^{\tau-t} (Pu)(\tau,\X^{\tau}_t(x, \Gamma_{t,x}^{-1}(w))\, \mathrm{d}\tau,
\end{align*}
we infer that
\begin{align*}
\Lambda_{t,w}(x)=x+(e^{-t}-1)w+(1-t+e^{-t})G + \int_0^{t} (e^{\tau-t}-1) (Pu)(\tau,\X^{\tau}_t(x,\Gamma_{t,x}^{-1}(w)))\, \mathrm{d}\tau.
\end{align*}
We then refer to \cite[Lemma 3.26]{HKM}, the proof of which can be exactly adapted to the case of the characteristic curves with an additional gravity term.
\end{proof}

\subsection{Local in time estimates}
In this subsection, we derive some local in time estimates for the moments and velocity field. Recalling the quantities of Definition \ref{decay_f0}, we introduce the following useful notations, which allows us to track down the dependency on the initial data in the later estimates. In view of the smallness condition \eqref{smallness:condition:THM}, this will enable us to set up a bootstrap strategy in the proof of Theorem \ref{thm1}.

\begin{nota}\label{nota:lesssim_0}
The notation $m \lesssim_0 M$ means that there exist some exponent $p>3$ and a continuous increasing function $\varphi : \R^+ \rightarrow \R^+$ such that
%\begin{align}\label{Lesssim0}
%m \leq \varphi \Big( 1+\Vert u_0 \Vert_{\H^1\cap \mathrm{D}_p^{1-\frac{1}{p},p}(\R^3_+)}+ N_q(f_0) + H_{\widetilde{q}}(f_0)+ K_{m,r}(f_0)+F_{\widetilde{m},\widetilde{r}}(f_0) + M_{6} f_0  + \mathrm{E}(0) + \Vert u_0 \Vert_{\Ld^1(\R^3_+)} \Big)M,
%\end{align}
\begin{align}\label{Lesssim0}
m \leq \varphi \Big( 1+\mathscr{E}(0) + \mathscr{N}_{q,m}(f_0)  \Big)M,
\end{align}
where $(q,m)$ are the exponents introduced in Theorem \ref{thm1} and where
\begin{align*}
\mathscr{E}(0)&=\Vert u_0 \Vert_{\H^1 \cap \mathrm{D}_p^{1-\frac{1}{p},p}(\R^3_+)}\hspace{-4mm} \ \  +\mathrm{E}(0)+\Vert u_0 \Vert_{\Ld^1(\R^3_+)}+ M_{6} f_0, \\
\mathscr{N}_{q,m}(f_0)&=N_q (f_0)+H_{q,m}(f_0)+ \underset{r \in \left\lbrace 1,3 \right\rbrace}{\max} \left\lbrace K_{q,r}(f_0)  +F_{q,m,r}(f_0) \right\rbrace.
\end{align*}
Here, $\mathrm{D}_p^{1-\frac{1}{p},p}(\R^3_+)$ refers to the space defined in \eqref{domaineDqs}, while the semi-norms involved in $\mathscr{N}_{q,m}(f_0)$ have been set in Definition \ref{decay_f0}. Note that in the end of the bootstrap argument, we shall be able to consider only the largest exponents $q$ and $m$ which are involved in the estimates.
%Furthermore, we will denote by $\mathrm{A}_0$ a generic constant which refers to an expression of the type $\varphi \Big( 1+\Vert u_0 \Vert_{\H^1\cap \mathrm{D}_p^{1-\frac{1}{p},p}(\R^3_+)}+ N_q(f_0) + H_{\widetilde{q}}(f_0)+ K_{m,r}(f_0)+F_{\widetilde{m},\widetilde{r}}(f_0) + M_{6} f_0  + \mathrm{E}(0) \Big)$ and which may vary from line to line.
\end{nota}
\textbf{Until the end of this Section \ref{Section:Prelim}, we consider a fixed global weak solution $(u,f)$ to the Vlasov-Navier-Stokes system, in the sense of Definition \ref{sol-faible} and associated to an admissible initial data $(u_0,f_0)$ satisfying \eqref{data:hyp}.}

%We first state the following Lemma, claiming that there is \textit{a priori} no conservation of the total mass of the distribution $f$ because of the absorption boundary condition \eqref{bcond-f}. 

\medskip

%We first state the following lemma, which is a direct consequence of the absorption boundary condition \eqref{bcond-f}. The proof is similar to that of \cite[Lemma 2.1]{EHKM}.
We first state the following lemma, entailing some rough bounds on the kinetic distribution. Note that the first one is a direct consequence of the absorption boundary condition \eqref{bcond-f}.
\begin{lem}
For all $t \geq  0$, we have
\begin{align}
\label{loose:mass}\int_{\R^3_+ \times \R^3} f(t,x,v) \, \mathrm{d}x \, \mathrm{d}v &\leq \int_{\R^3_+ \times \R^3} f_0(x,v) \, \mathrm{d}x \, \mathrm{d}v, \\[2mm]
\label{max:principle} \Vert f(t) \Vert_{\Ld^{\infty}(\R^3_+ \times \R^3)} & \leq e^{3t} \Vert f_0\Vert_{\Ld^{\infty}(\R^3_+ \times \R^3)}.
\end{align}
\end{lem}
\begin{proof}
We rely on the strong stability results from DiPerna-Lions theory about transport equations on $\R^3_+ \times \R^3$ (see \ref{DiPernaLions} in the Appendix). In short, it allows to prove the desired estimate for a sequence of distributions $(f_n)$ associated to a sequence of regularized initial data $(f_{0,n})_n$ and an approximating sequence of fluid velocities $(u_n)_n$. In particular, the associated characteristic curves \eqref{EDO-charac} are defined in a classic way. The strong stability property of renormalized solutions to the Vlasov equation is then used to recover the estimate for the original solution $(u,f)$. We do not detail the argument and we write the proof as if $u$ and $f_0$ were smooth.
%\footnote{We mention the fact that the proof of \cite[Lemma 2.1]{EHKM} contains an inacurracy: since $f$ may not be smooth for regular data and vector field $u$, one must rely on the characteristic curves and the representation formula \eqref{formule-rep}. The same procedure has to be performed in the proof of \cite[Lemma 4.5]{EHKM}, as it is done in the proof of Lemma \ref{interpo-estimate}.}
By Proposition \ref{Prop:formulerep}, we have
\begin{align*}
 f(t,x,v)= e^{3t} \mathbf{1}_{\mathcal{O}^t}(x,v) \, f_0(\mathrm{Z}^0_{t}(x,v)),
\end{align*}
therefore 
\begin{align*}
\int_{\R^3_+ \times \R^3} f(t,x,v) \, \mathrm{d}x \, \mathrm{d}v &=\int_{\R^3_+ \times \R^3} e^{3t} \mathbf{1}_{\mathcal{O}^t}(x,v) \, f_0(\mathrm{Z}^0_{t}(x,v)) \, \mathrm{d}x \, \mathrm{d}v \\
&=\int_{\R^3_+ \times \R^3}  \mathbf{1}_{\tau^{+}(0,x,v)>t} \, f_0(x,v) \, \mathrm{d}x \, \mathrm{d}v,
\end{align*}
thanks to the change of variable $(x,v) \mapsto \mathrm{Z}^t_{0}(x,v)$ (see Proposition \ref{Propo:diffeoZ} and Lemma \ref{LM:tau+-}). The inequality \eqref{loose:mass} follows, as well as \eqref{max:principle}.
\end{proof}
%\begin{nota}
%We set
%\begin{align*}
%F&:=j_f-\rho_f u,\\
%S&:=F- (u \cdot\nabla)u.
%\end{align*}
%\end{nota}
An application of Hölder's inequality implies the following result.
\begin{lem}
\label{lem-Dp}
Let $p>1$. For all $t \geq0$, we have 
\begin{equation}
\|( j_f - \rho_f u )(t)\|_{\mathrm{L}^p(\R^3_+)} \leq   \| \rho_f(t) \|_{\Ld^\infty(\R^3_+))}^{\frac{p-1}{p}} \left(\int_{\R^3_+ \times \R^3} f(t,x,v) |v-u(t,x)|^p \, \mathrm{d} x \, \mathrm{d} v\right)^{1/p}.
%\|( j_f - \rho_f u )(t)\|_{\mathrm{L}^p(\R^3_+)}^p \leq \int_{\R^3_+} \rho_f(t,x) \left( \int_{\R^3} f(t,x,v) |v-u(t,x)|^p \, \mathrm{d} v \right) \, \mathrm{d} x.
\end{equation}
\end{lem}
Along the way, we will also need the simple following sublinear Grönwall's lemma.
\begin{lem}\label{gronwall-sublinear}
Let $y \in \mathscr{C}(\R^+;\R^+)$ and $\mathrm{h} \in \mathscr{C}(\R^+;\R^+)$ such that for all $t \geq 0$
$$y(t) \leq y_0 +\int_0^t \mathrm{h}(s) y(s)^{\beta} \, \mathrm{d}s,$$
where $\beta \in (0,1)$ and $y_0 \in \R^+$. Then for all $t \geq 0$
\begin{align*}
y(t) \leq  \left( y_0^{1-\beta}+(1-\beta)\int_0^t \mathrm{h}(s)\, \mathrm{d}s  \right)^{\frac{1}{1-\beta}}.
\end{align*}
\end{lem}
%\begin{proof}
%We define 
%$$\widetilde{y}(t):= \varepsilon+y_0 +\int_0^t \mathrm{h}(s) y(s)^{\beta} \, \mathrm{d}s,$$
%for some fixed $\varepsilon>0$. Since $\widetilde{y} \geq y$, we have
%$$\widetilde{y} \, '(t)=\mathrm{h}(t) y(t)^{\beta} \leq \mathrm{h}(t) \widetilde{y}(t)^{\beta}.$$
%Since $ \widetilde{y}$ does not vanish, we get 
%$$\dfrac{\widetilde{y} \, '(t)}{\widetilde{y}(t)^{\beta}} \leq \mathrm{h}(t),$$
%so that 
%$$y(t) \leq \widetilde{y}(t) \leq \left( (\varepsilon+y_0)^{1-\beta}+ (1-\beta)\int_0^t \mathrm{h}(s)\, \mathrm{d}s \right)^{\frac{1}{1-\beta}}.$$
%We conclude by letting $\varepsilon \rightarrow 0$. 
%\end{proof}

We now state several rough bounds on the moments of the distribution function $f$.
\begin{lem}\label{j_f:L^1L^1}
For all $t \geq0$, we have 
$$ \Vert j_f \Vert_{\Ld^1(0,t;\Ld^1(\R^3_+))} \lesssim \E(0)^{\frac{1}{2}}t + t^2,$$
%where $\lesssim$ only depends on $g$. 
In particular, the map $t \mapsto \Vert j_f \Vert_{\Ld^1(0,t;\Ld^1(\R^3_+))}$ belongs to $\Ld^{\infty}_{\mathrm{loc}}(\R^+)$.
\end{lem}
\begin{proof}
We use the energy inequality \eqref{ineq-energy} and the Cauchy-Schwarz inequality to write that for all $t \geq 0$
\begin{align*}
\E(t) \leq \E(0) + g\int_0^t \int_{\R^3_+ \times \R^3} \vert v \vert f(s,x,v) \, \mathrm{d}x \, \mathrm{d}v  \, \mathrm{d}s  \leq \E(0) + g  \int_0^t \E(s)^{\frac{1}{2}} \, \mathrm{d}s,
\end{align*}
thanks to \eqref{loose:mass} and the normalization of $f_0$ given by \eqref{CI-f2}. The sub-linear version of Grönwall's Lemma \ref{gronwall-sublinear} with $\beta=1/2$ implies
\begin{align*}
\E(t) \leq \left( \E(0)^{\frac{1}{2}}+\frac{g}{2} t \right)^2.
\end{align*}
The Cauchy-Schwarz inequality again reads
\begin{align*}
\int_0^t \int_{\R^3_+ \times \R^3} \vert v \vert f(s,x,v)  \, \mathrm{d}x \, \mathrm{d}v  \, \mathrm{d}s \leq  \int_0^t \E(s)^{\frac{1}{2}} \, \mathrm{d}s \leq t \E(0)^{\frac{1}{2}}+g\frac{t^2}{4},
\end{align*}
and this concludes the proof.
\end{proof}

\begin{lem}\label{ineg-Brinkman}
If $F:=j_f-\rho_f u$ then for all $t \geq 0$, we have
\begin{align*}
\int_0 ^t \Vert F(s) \Vert_{\Ld^2(\R^3_+)}^2   \, \mathrm{d}s \lesssim \underset{s \in [0,t]}{\sup} \Vert \rho_f (s) \Vert_{\Ld^{\infty}(\R^3_+)}  \left[ \E(0) + \E(0)^{\frac{1}{2}}t + t^2 \right].
\end{align*}
\end{lem}
\begin{proof}
Using Lemma \ref{lem-Dp} with $p=2$, we obtain for all $s \in [0,t]$
$$ \Vert F(s) \Vert_{\Ld^2(\R^3_+)}^2 \leq  \Vert \rho_f (s) \Vert_{\Ld^{\infty}(\R^3_+)} \D(s),$$
where the dissipation $\D$ has been defined in \eqref{eq:Dissipation}. Then, we integrate the last inequality between $0$ and $t$ and use the energy inequality (\ref{ineq-energy}) to obtain
\begin{align*}
\int_0 ^t \Vert F(s) \Vert_{\Ld^2(\R^3_+)}^2 \, \mathrm{d}s \leq \underset{s \in [0,t]}{\sup} \Vert \rho_f (s) \Vert_{\Ld^{\infty}(\R^3_+)}  \left[ \E(0) + g \int_0^t \int_{\R^3_+} \vert j_f \vert \, \mathrm{d}\tau  \, \mathrm{d}x \right],
\end{align*}
which concludes the proof, thanks to Lemma \ref{j_f:L^1L^1}.
\end{proof}

We recall standard interpolation estimates on the moments of any kinetic distribution, where we use the notations introduced in Definition \ref{notation:moments}.
%see +ou - [Lions-Perthame (1991)] pour le premier et  [LM 2.2 Leopold-Saffirio (2022)] pour le second
\begin{propo}\label{interpo-moment}
Let $\mathrm{h}$ be a nonnegative function in $\Ld^{\infty}(\R^+ \times \R^3_+ \times \R^3)$. Then we have for all $t \geq0$
\begin{align*}
%m_{\ell}g(t,x) &\lesssim \left( \Vert h\Vert_{\Ld^{\infty}(\R^+ \times \R^3_+ \times \R^3)} +1 \right) m_{k}h(t,x)^{\frac{\ell +3}{k+3}},\\
\forall  \ 0 \leq b \leq c, \ \ \forall \ell \in [b,c], \ \ M_{\ell} \mathrm{h}(t) &\leq M_{b} \mathrm{h}(t)^{\frac{c-\ell}{c-b}} M_{c} \mathrm{h}(t)^{\frac{\ell-b}{c-b}}, \\
\forall k>0, \ \ \forall \ell \in [0,k], \ \  \Vert m_{\ell}\mathrm{h}(t) \Vert_{\Ld^{\frac{k+3}{\ell+3}}(\R^3_+)} &\leq C_{k,\ell} \Vert \mathrm{h}(t) \Vert_{\Ld^{\infty}(\R^3_+ \times \R^3)}^{\frac{k-\ell}{k+3}}  M_{k}\mathrm{h}(t)^{\frac{\ell +3}{k+3}},
\end{align*}
for some universal constant $C_{k,\ell}>0$.
\end{propo}
We now provide a pointwise estimate for the moments of $f$, solution to the Vlasov equation.
\begin{lem}\label{interpo-estimate}
Suppose that $u \in \Ld^1_{\mathrm{loc}}(\R^+; \Ld^{a+3} \cap \W^{1,1}(\R^3_+))$ and $M_{a} f_0 < \infty$ for all $a \in [2, \alpha]$, for some $\alpha \geq 2$. Then for all $a \in [2,\alpha]$ and for all $t>0$, $M_{a} f(t) < \infty$ and $M_a f \in \Ld^{\infty}_{\mathrm{loc}}(\R^+)$.
Furthermore, if $T>0$ then for all $t \in [0,T]$ 
\begin{align}\label{interpo-moment-vlasov}
M_{\alpha} f(t) \leq \left( \left[ M_{\alpha} f_0 + \alpha g\int_0^T M_{\alpha-1} f(s) \, \mathrm{d}s \right]^{\frac{1}{\alpha+3}} + e^{\frac{3t}{\alpha+3}} \Vert f_0 \Vert_{\Ld^{\infty}(\R^3_+ \times \R^3)}^{\frac{1}{\alpha+3}} \int_0^t \Vert u(s) \Vert_{\Ld^{\alpha+3}(\R^3_+)}\mathrm{d}s \right)^{\alpha+3}.
\end{align}
\end{lem}
\begin{proof}
As in the proof of the bounds \eqref{loose:mass}--\eqref{max:principle}, we rely on the strong stability results from DiPerna-Lions theory about transport equations on $\R^3_+ \times \R^3$. We do not detail the argument and we write the proof as if $u$ and $f_0$ were smooth. By Proposition \ref{Prop:formulerep}, we have
\begin{align*}
 f(t,x,v)= e^{3t} \mathbf{1}_{\mathcal{O}^t}(x,v) \, f_0(\mathrm{Z}^0_{t}(x,v)),
\end{align*}
so that
\begin{align*}
M_{\alpha} f(t)= \int_{\R^3_+ \times \R^3}  \vert v \vert^{\alpha} e^{3t} \mathbf{1}_{\mathcal{O}^t}(x,v) \, f_0(\mathrm{Z}^0_{t}(x,v)) \, \mathrm{d}x \, \mathrm{d}v=\int_{\R^3_+ \times \R^3}  \vert \mathrm{V}^t_{0}(x,v) \vert^{\alpha} \mathbf{1}_{\tau^{+}(0,x,v)>t} \, f_0(x,v) \, \mathrm{d}x \, \mathrm{d}v,
\end{align*}
thanks to the change of variable $(x,v) \mapsto \mathrm{Z}^t_{0}(x,v)$ (see Proposition \ref{Propo:diffeoZ} and Lemma \ref{LM:tau+-}).

In view of the first estimate of Proposition \ref{interpo-moment}, it is sufficient to prove the formula \eqref{interpo-moment-vlasov} for $\alpha \geq 2$ being an integer. We then argue by induction on $ \alpha$. For $\alpha=2$, we observe 
\begin{align*}
\frac{\mathrm{d}}{\mathrm{d} s}& |\V^s_0(x,v))|^2 = 2 \frac{\mathrm{d}}{\mathrm{d} s} [\V^s_0(x,v))] \cdot \V^s_0(x,v) = 2 \left[ u(s,\X^s_0(x,v)) - \V^s_0(x,v) +G  \right]  \cdot \V^s_0(x,v),
\end{align*}
from which we infer that
\begin{align*}
 |\V^t_0(x,v))|^2 \leq \vert v \vert^2+2 \int_0^t \left[u(s,\X^s_0(x,v)) + G \right] \cdot \V^s_0(x,v)  \, \mathrm{d}s.
\end{align*}
By Fubini Theorem, we obtain
\begin{multline*}
M_{2} f(t) \leq \int_{\R^3_+ \times \R^3}  \vert v \vert^{2} \mathbf{1}_{\tau^{+}(0,x,v)>t} \, f_0(x,v) \, \mathrm{d}x \, \mathrm{d}v  
+ 2 g\int_0^t \int_{\R^3_+ \times \R^3}  \vert \V^s_0(x,v) \vert \mathbf{1}_{\tau^{+}(0,x,v)>t} \, f_0(x,v) \, \mathrm{d}x \, \mathrm{d}v \\ + 2 \int_0^t \int_{\R^3_+ \times \R^3}  \vert u(s,\X^s_0(x,v)) \vert \vert \V^s_0(x,v) \vert \mathbf{1}_{\tau^{+}(0,x,v)>t} \, f_0(x,v) \, \mathrm{d}x \, \mathrm{d}v.
\end{multline*}
Using the reverse change of variable $\mathrm{Z}^t_{0}(x,v) \mapsto (x,v)$ in the two last integrals, we get 
\begin{align*}
M_{2} f(t) & \leq M_{2} f_0+ 2 g \int_0^t M_1 f(s)  \, \mathrm{d}s+ 2 \int_0^t \int_{\R^3_+} \vert u(s,x) \vert  m_{1}f(s,x) \, \mathrm{d}x \, \mathrm{d}s \\
\quad & \leq M_{2} f_0+ 2 g \int_0^t M_1 f(s)  \, \mathrm{d}s+ 2 \int_0^t \Vert u(s) \Vert_{\Ld^{5}(\R^3_+)} \Vert m_{1}f(s) \Vert_{\Ld^{\frac{5}{4}}(\R^3_+)} \, \mathrm{d}s,
\end{align*}
thanks to Hölder's inequality. Furthermore, by Proposition \ref{interpo-moment} with $(\ell,k)=(1,2)$ and the rough control provided by \ref{max:principle}, we get
\begin{align*}
\Vert m_{1}f(s) \Vert_{\Ld^{\frac{5}{4}}(\R^3_+)} \lesssim  \Vert f(s) \Vert_{\Ld^{\infty}(\R^3_+ \times \R^3)}^{\frac{1}{5}} M_{2}f(s)^{\frac{4}{5}} \leq \Vert f_0 \Vert_{\Ld^{\infty}(\R^3_+ \times \R^3)}^{\frac{1}{5}} e^{\frac{3s}{5}}M_{2}f(s)^{\frac{4}{5}},
\end{align*}
where $\lesssim$ is independent of $s$. We thus infer that for all $t \in [0,T]$
\begin{align*}
M_2 f(t) \leq M_{2} f_0+ 2 g \int_0^T M_1 f(s)  \, \mathrm{d}s + 2 \Vert f_0 \Vert_{\Ld^{\infty}(\R^3_+ \times \R^3)}^{\frac{1}{5}} \int_0^t e^{\frac{3s}{5}} \Vert u(s) \Vert_{\Ld^{5}(\R^3_+)} M_{2}f(s)^{\frac{4}{5}} \, \mathrm{d}s.
\end{align*}
Using the Grönwall's lemma stated in Lemma \ref{gronwall-sublinear} with $\beta=4/5$ entails
$$M_{2} f(t) \lesssim \left( \left[ M_{2} f_0 + 2 g\int_0^T M_{1} f(s) \, \mathrm{d}s \right]^{\frac{1}{5}} + \dfrac{2}{5} e^{\frac{3t}{5}} \Vert f_0 \Vert_{\Ld^{\infty}(\R^3_+ \times \R^3)}^{\frac{1}{5}} \int_0^t \Vert u(s) \Vert_{\Ld^{5}(\R^3_+)}\mathrm{d}s \right)^{5}.$$
This yields the result for $\alpha=2$ (indeed, note that the previous right-hand side is finite because of Lemma \ref{j_f:L^1L^1}). Now, if $M_{\alpha-1} f(t) < \infty$ and $M_{\alpha-1} f \in \Ld^{\infty}_{\mathrm{loc}}(\R^+)$, we perform the same analysis as before. Since $\alpha \geq 2$, we have
\begin{align*}
\frac{\mathrm{d}}{\mathrm{d} s} |\V^s_0(x,v))|^{\alpha} &= \alpha \frac{\mathrm{d}}{\mathrm{d} s} [\V^s_0(x,v))] \cdot \V^s_0(x,v)  \times |
\V(s;0,x,v) |^{\alpha-2} \\[2mm]
&= \alpha \left[ u(s,\X^s_0(x,v)) - \V^s_0(x,v) +G  \right] \cdot \V^s_0(x,v)  |\V^s_0(x,v)|^{\alpha-2},
\end{align*}
which entails 
\begin{align*}
 |\V^t_0(x,v))|^{\alpha} \leq \vert v \vert^{\alpha}+\alpha\int_0^t \left[u(s,\X^s_0(x,v)) + G \right] \cdot \V^s_0(x,v)  \vert \V^s_0(x,v) \vert^{\alpha-2}\, \mathrm{d}s.
\end{align*}
As before, we obtain
\begin{align*}
M_{\alpha} f(t) \leq M_{\alpha}f_0 +\alpha g \int_0^t M_{\alpha-1}f(s) \, \mathrm{d}s + \alpha \int_0^t \Vert u(s) \Vert_{\Ld^{\alpha+3}(\R^3_+)} \Vert m_{\alpha-1}f(s) \Vert_{\Ld^{\frac{\alpha+3}{\alpha+2}}(\R^3_+)} \, \mathrm{d}s.
\end{align*}
Using Proposition \ref{interpo-moment} with $(\ell,k)=(\alpha-1,\alpha)$ and the inequality \ref{max:principle}, we get $$\Vert m_{\alpha-1}f(s) \Vert_{\Ld^{\frac{\alpha+3}{\alpha+2}}(\R^3_+)} \lesssim \Vert f_0 \Vert_{\Ld^{\infty}}^{\frac{1}{\alpha+3}} e^{\frac{3s}{\alpha+3}}  M_{\alpha}f(s)^{\frac{\alpha+2}{\alpha+3}}.$$
This yields for all $t \in [0,T]$
\begin{align*}
M_{\alpha} f(t) \leq M_{\alpha} f_0+ \alpha g \int_0^T M_{\alpha-1} f(s)  \, \mathrm{d}s + \alpha \Vert f_0 \Vert_{\Ld^{\infty}(\R^3_+ \times \R^3)}^{\frac{1}{\alpha+3}} \int_0^t e^{\frac{3s}{\alpha+3}} \Vert u(s) \Vert_{\Ld^{\alpha+3}(\R^3_+)} M_{\alpha}f(s)^{\frac{\alpha+2}{\alpha+3}} \, \mathrm{d}s.
\end{align*}
Thanks to Lemma \ref{gronwall-sublinear} with $\beta=\frac{\alpha+2}{\alpha+3}$, we obtain the conclusion.
\end{proof}

\begin{lem}\label{rough:bounds}
We have $M_3 f\in \Ld^{\infty}_{\mathrm{loc}}(\R^+)$. Moreover, for all finite $T>0$, there exists a continuous nonnegative and nondecreasing function $\varphi_{\E(0), M_3 f_0,T, \Vert f_0 \Vert_{\Ld^{\infty}_{x,v}}}$ (increasing in all its parameters) such that for all $t \in [0,T]$
\begin{align*}
\Vert \rho_f (t) \Vert_{\Ld^{2}(\R^3_+)} + \Vert j_f (t) \Vert_{\Ld^{3/2}(\R^3_+)} \leq \varphi_{\E(0), M_3 f_0,T, \Vert f_0 \Vert_{\Ld^{\infty}_{x,v}}}(t). 
\end{align*}
\end{lem}
\begin{proof}
Since $M_2 f_0<\infty$ and $M_6 f_0<\infty$ by \eqref{data:hyp}, we have $M_3 f_0 \lesssim M_2 f_0 +M_6 f_0 < \infty$. Furthermore, $u$ is a Leray solution so that by the Sobolev embedding, we can apply Lemma \ref{interpo-estimate} with $\alpha = 3$ and we deduce that $M_3 f\in \Ld^{\infty}_{\mathrm{loc}}(\R^+)$. Furthermore, the estimate \eqref{interpo-moment-vlasov} yields for all $t \in [0,T]$
\begin{align*}
M_{3} f(t) \lesssim \left( \left[ M_{3} f_0 + 3 g\int_0^T M_{2} f(s) \, \mathrm{d}s \right]^{\frac{1}{6}} + \dfrac{e^{\frac{t}{2}}}{2}  \Vert f_0 \Vert_{\Ld^{\infty}(\R^3_+ \times \R^3)}^{\frac{1}{6}} \int_0^t \Vert u(s) \Vert_{\Ld^6(\R^3_+)} \mathrm{d}s \right)^6.
\end{align*}
The Sobolev embedding on $\H^1_0(\R^3_+)$ and the Cauchy-Schwarz inequality lead to
\begin{align*}
\int_0^t \Vert u(s) \Vert_{\Ld^6(\R^3_+)} \mathrm{d}s \lesssim t^{1/2} \left( \int_0^t \Vert \nabla u(s) \Vert_{\Ld^2(\R^3_+)}^2 \mathrm{d}s\right)^{1/2} \lesssim T^{1/2} \left[ \mathrm{E}(0) +\E(0)^{\frac{1}{2}}T + T^2 \right]^{1/2} ,
\end{align*}
where we have used the energy inequality (\ref{ineq-energy}) and Lemma \ref{j_f:L^1L^1}. For the same reasons, we also have
\begin{align*}
\int_0^T M_{2} f(s) \, \mathrm{d}s \lesssim \mathrm{E}(0) +\E(0)^{\frac{1}{2}}T + T^2.
\end{align*}
This implies that there exists a continuous nonnegative and nondecreasing function $\varphi=\varphi_{\E(0), M_3 f_0,T, \Vert f_0 \Vert_{\Ld^{\infty}_{x,v}}}$ (increasing in all its parameters) such that for all $t \in [0,T]$, we have
\begin{align}\label{bound:M3}
M_3 f(t) \leq \varphi_{\E(0), M_3 f_0,T, \Vert f_0 \Vert_{\Ld^{\infty}_{x,v}}}(t).
\end{align}
Using Proposition \ref{interpo-moment} on interpolation of moments of the distribution function $f$ with $k=3$ and $\ell \in \{0,1\}$, together with \ref{max:principle}, we get
\begin{align*}
\Vert \rho_f (t) \Vert_{\Ld^{2}(\R^3_+)} &\lesssim \Vert f(t) \Vert_{\Ld^{\infty}(\R^3_+)}^{1/2} M_3 f(t)^{1/2} \leq \Vert f_0 \Vert_{\Ld^{\infty}(\R^3_+)}^{1/2} e^{3t/2}M_3 f(t)^{1/2}, \\
\Vert j_f (t) \Vert_{\Ld^{3/2}(\R^3_+)} &\lesssim \Vert f(t) \Vert_{\Ld^{\infty}(\R^3_+)}^{1/3} M_3 f(t)^{2/3} \leq \Vert f_0 \Vert_{\Ld^{\infty}(\R^3_+)}^{1/3} e^{t}  M_3 f(t)^{2/3},
\end{align*}
which yields the conclusion, thanks to the bound \eqref{bound:M3}.
\end{proof}

We now prove that the source term in the Navier-Stokes equations, namely the Brinkman force $j_f-\rho_f u$, belongs to $\Ld^2 \Ld^2$ (locally in time). The strategy of proof is very similar to that of Lemma \cite[Lemma 4.7]{EHKM} (with a minor adaptation to the half-space case) and details are thus omitted. Let us only emphasize that the result mainly follows from the combination of Lemma \ref{j_f:L^1L^1}, Lemma \ref{rough:bounds}, Lemma \ref{interpo-moment} and the maximal regularity property for the Stokes system (see Section \ref{AnnexeMaxregStokes} in the Appendix) which can be applied in that case because $u_0 \in \W_0^{1,9/7}(\R^3_+)$ and $M_6 f_0<\infty$ (see the assumption \eqref{data:hyp}).
%but we write down the whole argument for the sake of completeness.
\begin{lem}\label{integBrinkman-1}
We have
\begin{align*}
j_f - \rho_f u \in \Ld^2_{\mathrm{loc}}(\R^+;\Ld^2(\R^3_+)).
\end{align*}
\end{lem}
We are then in position to state the following local in time integrability results of the Leray solution $(u,f)$.
\begin{propo}\label{controlLinfiniLOC}
We have 
\begin{align}
\label{borneL1Linfty-u} u \in \Ld^{1}_{\mathrm{loc}}(\R^+;\Ld^{\infty}(\R^3_+)), 
\end{align}
and if $N_d(f_0)<\infty$ for some $d>4$, then
\begin{align*}
\rho_f \in  \Ld^{\infty}_{\mathrm{loc}} (\R^+; \Ld^{\infty}(\R^3_+)), \\
j_f \in  \Ld^{\infty}_{\mathrm{loc}} (\R^+; \Ld^{\infty}(\R^3_+)).
\end{align*}
More precisely, there exists a continuous nondecreasing function $\eta : \R^+\rightarrow \R^+$ such that
\begin{align}
\Vert u \Vert_{\Ld^{1}(0,t;\Ld^{\infty}(\R^3_+))} &\lesssim \eta(t), \\
\Vert \rho_f \Vert_{\Ld^{\infty}(0,t;\Ld^{\infty}(\R^3_+))} + \Vert j_f \Vert_{\Ld^{\infty}(0,t;\Ld^{\infty}(\R^3_+))}  & \lesssim N_d(f_0) \eta(t). \label{borneLinftyLinfty-rhoetj}
\end{align}
\end{propo}
\begin{proof}
Let $T>0$. The proof of the fact that $u \in \Ld^{1}(0,T;\Ld^{\infty}(\R^3_+))$ is mostly directly taken from the arguments used in \cite[Proposition 4.8]{EHKM} and mainly relies on the theory of epochs of regularity for the Leray solutions to the Navier-Stokes equations. 
Owing to \cite[Theorem 8]{Hey1} and \cite[Remark 4]{Hey2} (which are valid since the strong energy inequality \eqref{ineq-energyNS} is satisfied by the weak solutions that we consider), we know there exists a subset $\sigma_T \subset [0,T]$ of full measure in $[0,T]$ with $\sigma_T=\bigsqcup_i ]a_i, b_i[$ (the union being countable) and for which $u \in \Ld^{\infty}_{\mathrm{loc}}(a_i,b_i;\H^1(\R^3_+)) \cap \Ld^2_{\mathrm{loc}}(a_i,b_i; \H^2(\R^3_+))$ and $\Vert \nabla u(t) \Vert_{\Ld^2(\R^3_+)} \underset{t \rightarrow b_i^-}{\longrightarrow}+\infty$ for all $i$. Furthermore, the function $t \mapsto \Vert \nabla u(t) \Vert^2$ is absolutely continuous on each interval $]a_i, b_i[$ (see e.g. \cite{Rob}).

Now, we can take the $\Ld^2(\R^3_+)$ inner product of \eqref{eq:NS} with $Au$ on each interval $]a_i, b_i[$, where $A=-\P\delta_x u$ stands for the Stokes operator on $\Ld^2_{\mathrm{div}}(\R^3_+)$ and $\P$ is the Leray projection on divergence-free vector field (see in the Appendix \ref{AnnexeMaxregStokes}). We obtain
\begin{align*}
\dfrac{\mathrm{d}}{\mathrm{d}t}\Vert \nabla u \Vert_{\Ld^2(\R^3_+)}^2 + 2\Vert  Au \Vert_{\Ld^2(\R^3_+)}^2+ 2\langle  \mathbb{P} (u \cdot \nabla)u , Au\rangle=2\langle \mathbb{P}( j_f-\rho_f u), Au\rangle \ \ \text{on} \ \ ]a_i, b_i[,
\end{align*}
where we have dropped the time variable. In order to estimate the term $ \langle \mathbb{P} (u \cdot \nabla)u , Au\rangle$, we use the Gagliardo-Nirenberg-Sobolev inequality for the function $\nabla u$ with the exponents $(p,j,r,m,\alpha)=(3,0,2,1,1/2)$, which reads as
\begin{align*}
\Vert \nabla u \Vert_{\Ld^3(\R^3_+)} \lesssim \Vert \D^2 u \Vert_{\Ld^2(\R^3_+)}^{1/2}  \Vert \nabla u \Vert_{\Ld^2(\R^3_+)}^{1/2},
\end{align*}
and we combine this inequality with the Hölder's inequality to write 
\begin{align*}
\vert \langle \mathbb{P} (u \cdot \nabla)u , Au\rangle \vert  &\leq  \Vert  u \Vert_{\Ld^6(\R^3_+)} \Vert \nabla u \Vert_{\Ld^3(\R^3_+)} \Vert A u \Vert_{\Ld^2(\R^3_+)} \\
& \lesssim \Vert \nabla u \Vert_{\Ld^2(\R^3_+)}^{3/2} \Vert  A u \Vert_{\Ld^2(\R^3_+)}^{3/2}.
\end{align*}
Note that we have used \cite[Theorem IV.3.2]{Galdi} on each $(a_i,b_i)$. Thanks to the Young inequality, we infer that 
\begin{align}\label{ineg-enstrophy}
\dfrac{\mathrm{d}}{\mathrm{d}t}\Vert \nabla u \Vert_{\Ld^2(\R^3_+)}^2 + \Vert  Au \Vert_{\Ld^2(\R^3_+)}^2\lesssim \Vert j_f-\rho_f u \Vert^2_{\Ld^2(\R^3_+)} + \Vert \nabla u \Vert_{\Ld^{2}(\R^3_+)}^{6},
\end{align}
on each interval $]a_i, b_i[$, where $\lesssim$ is independent of the time variable and independent of $i$. Dividing this inequality by $(1+\Vert \nabla u \Vert_{\Ld^2(\R^3_+)}^2)^2$, we get
\begin{align*}
\dfrac{1}{(1+\Vert \nabla u \Vert_{\Ld^2(\R^3_+)}^2)^2}\dfrac{\mathrm{d}}{\mathrm{d}t}\Vert \nabla u \Vert_{\Ld^2(\R^3_+)}^2 + \dfrac{1}{(1+\Vert \nabla u \Vert_{\Ld^2(\R^3_+)}^2)^2}\Vert  Au \Vert_{\Ld^2(\R^3_+)}^2  \lesssim
\Vert j_f-\rho_f u \Vert^2_{\Ld^2(\R^3_+)} +\Vert \nabla u \Vert_{\Ld^{2}(\R^3_+)}^{2},
\end{align*}
on each interval $]a_i, b_i[$. Integrating and summing over the previous epoch of regularities, an using in particular the fact that $\Vert \nabla u(t) \Vert_{\Ld^2(\R^3_+)} \underset{t \rightarrow b_i^-}{\longrightarrow}+\infty$ for all $i$, we can perform the same exact computations as in the proof of \cite[Proposition 4.8]{EHKM} and end up with
\begin{align*}
\int_0^T\dfrac{\Vert  Au(s) \Vert_{\Ld^2(\R^3_+)}^2}{(1+\Vert \nabla u(s) \Vert_{\Ld^2(\R^3_+)}^2)^2}\, \mathrm{d}s 
  \leq \int_0^T \left(\Vert j_f(s)-\rho_f u(s) \Vert^2_{\Ld^2(\R^3_+)} + \Vert \nabla u(s) \Vert_{\Ld^{2}(\R^3_+)}^{2}\right) \, \mathrm{d}s +1,
\end{align*}
from which we infer that
\begin{align*}
\int_0^T \Vert  Au(s) \Vert_{\Ld^2(\R^3_+)}^{2/3}\, \mathrm{d}s \leq \left( \int_0^T\dfrac{\Vert  Au(s) \Vert_{\Ld^2(\R^3_+)}^2}{(1+\Vert \nabla u(s) \Vert_{\Ld^2(\R^3_+)}^2)^2}\, \mathrm{d}s \right)^{1/3}\left( \int_0^T (1+\Vert \nabla u(s) \Vert_{\Ld^2(\R^3_+)}^2)\, \mathrm{d}s \right)^{2/3}<\infty,
\end{align*}
due to Proposition \ref{integBrinkman-1} and to the fact that $u$ is a Leray solution to the Navier-Stokes equations. Using the Gagliardo-Nirenberg-Sobolev inequality with exponents $(p,j,r,m,\alpha)=(\infty, 0,2,2,1/2)$ and Sobolev embedding, we deduce that
\begin{align*}
\int_0^T \Vert  u(s) \Vert_{\Ld^{\infty}(\R^3_+)}\, \mathrm{d}s &\lesssim \int_0^T \Vert  \nabla u(s) \Vert_{\Ld^2(\R^3_+)}^{2}\, \mathrm{d}s +\int_0^T \Vert  Au(s) \Vert_{\Ld^2(\R^3_+)}^{2/3}\, \mathrm{d}s <\infty,
\end{align*}
from the same reasons as before, therefore this proves \eqref{borneL1Linfty-u}. The last estimate \eqref{borneLinftyLinfty-rhoetj} is eventually obtained by observing that for $d>4$
\begin{align*}
\Vert  \rho_f(t) \Vert_{\Ld^{\infty}(\R^3_+)} + \Vert j_f(t) \Vert_{\Ld^{\infty}(\R^3_+)} & \lesssim N_d (f(t)) \\
& \lesssim e^{3 t} (1+ \vert G \vert^q+ \Vert u \Vert_{\mathrm{L}^1 (0,t ; \mathrm{L}^{\infty}(\R^3_+))} ^q) N_d(f_0) \\
& \lesssim N_q(f_0) \eta(t),
\end{align*}
thanks to Lemma \ref{inegdecal} and \eqref{borneL1Linfty-u}.
\end{proof}

\subsection{Strong existence times and higher order energy estimates}
Along the bootstrap procedure, we shall need $\H^1$ energy estimates for the fluid velocity $u$, which is \textit{a priori} only a Leray solution to the Navier-Stokes equations. In order to consider higher regularity for this solution, we rely on a parabolic smoothing property of the (Vlasov-)Navier-Stokes system. We will be able to propagate this regularity if the contribution of the source term, that is the Brinkman force $F:=j_f-\rho_f u$, and the initial data, are small enough.
\begin{propo}\label{propdatasmall:VNSreg}
There exists a universal constant $\mathrm{C}_{\star}$ such that the following holds. Assume that for some $T>0$, one has
\begin{align}\label{datasmall:VNSreg}
\Vert  u_0 \Vert_{\H^1(\R^3_+)}^2 +  \int_0^T \Vert F(s) \Vert_{\Ld^2(\R^3_+)}^2 \, \mathrm{d}s +\int_0^T \Vert F(s) \Vert_{\Ld^2(\R^3_+)} \, \mathrm{d}s< \mathrm{C}_{\star}.
\end{align}
Then one has 
$$u \in \Ld^{\infty}(0,T; \H^1(\R^3_+)) \cap \Ld^{2}(0,T; \H^2(\R^3_+)),  $$
and for all $t \in [0,T]$
\begin{align}\label{ineq:VNSreg}
\Vert \nabla u(t) \Vert_{\Ld^2(\R^3_+)}^2 + \dfrac{1}{2}\int_0^t \Vert \mathrm{D}^2 u(s) \Vert_{\Ld^2(\R^3_+)}^2 \, \mathrm{d}s   \lesssim \Vert \nabla u_0 \Vert_{\Ld^2(\R^3_+)}^2 + \int_0^t \Vert F(s) \Vert_{\Ld^2(\R^3_+)}^2 \, \mathrm{d}s,
\end{align}
where $\lesssim$ only depends on $\mathrm{C}_{\star}$.
\end{propo}
\begin{proof}
The estimate is a direct consequence of the parabolic regularization for the Navier-Stokes system with source $F=j_f-\rho_f u$, that we state in Theorem \ref{RegParabNS} in Section \ref{AnnexeParabNS} of the Appendix.
%, together with the estimate on the Brinkman force of Lemma \ref{ineg-Brinkman}.
\end{proof}
%\begin{rem}
%By choosing an appropriate function $\varphi$ in (\ref{smallness:condition}), we can actually ensure that
%\begin{align}\label{Hyp:H1small}
%\Vert \nabla u_0 \Vert_{\Ld^2(\R^3_+)}^2 \leq \dfrac{1}{2\sqrt{8C_1 C_2}}.
%\end{align}
%Furthermore, again in view of (\ref{smallness:condition}), we can also reduce $\Vert \nabla u_0 \Vert_{\Ld^2(\R^3_+)}^2$ and $\mathrm{E}(0)$ in the sequel if necessary.
%\end{rem}
%so that the conclusion of Proposition \ref{datasmall:VNSreg} turns into
%\begin{align}\label{ineq:VNSreg}
%\Vert \nabla u(t) \Vert_{\Ld^2(\R^3_+)}^2 + \dfrac{1}{2}\int_0^t \Vert Au(s) \Vert_{\Ld^2(\R^3_+)}^2 \, \mathrm{d}s \lesssim \mathrm{E}(0) \left( 1 +\underset{s \in [0,t]}{\sup} \Vert \rho_f (s) \Vert_{\Ld^{\infty}(\R^3_+)} \right).
%\end{align}
\begin{rem}\label{control:parab}
By Lemma \ref{ineg-Brinkman}, the right-hand side of (\ref{ineq:VNSreg}) is finite. In particular, if the condition (\ref{datasmall:VNSreg}) holds for some $T$, then for all $t \in [0,T]$
%Using the elliptic estimate $\Vert u \Vert_{\H^2(\R^3_+)} \lesssim \Vert Au \Vert_{\Ld^2(\R^3_+)}$,  we infer that 
%$$ u \in \Ld^{\infty}(0,T; \H^1(\R^3_+)) \cap \Ld^{2}(0,T;\H^2(\R^3_+)),
%$$
%and in particular $$ \nabla u \in \Ld^{2}(0,T;\H^1(\R^3_+)).$$
%In particular, we have
%$$u \in \Ld^{\infty}(0,T; \H^1(\R^3_+)),$$
%In particular, we have the following estimate
\begin{align*}
\Vert u \Vert_{\Ld^{\infty}(0,t; \Ld^6(\R^3_+))}^2 \lesssim \Vert \nabla u \Vert_{\Ld^{\infty}(0,t; \Ld^2(\R^3_+))}^2 &\lesssim \Vert \nabla u_0 \Vert_{\Ld^2(\R^3_+)}^2 + \underset{s \in [0,t]}{\sup} \Vert \rho_f (s) \Vert_{\Ld^{\infty}(\R^3_+)}  \left[ \E(0) + \E(0)^{\frac{1}{2}}t + t^2 \right].
%\\
%\Vert \nabla u - \left\langle \nabla u \right\rangle\Vert_{\Ld^2(0,t;\Ld^{6}(\R^3_+))}^2 & \lesssim \Vert \nabla u_0 \Vert_{\Ld^2(\R^3_+)}^2 + \mathrm{E}(0)\underset{s \in [0,t]}{\sup} \Vert \rho_f (s) \Vert_{\Ld^{\infty}(\R^3_+)}.\label{controlparab2}
\end{align*}
%where $\langle g \rangle$ stands for the average of a function $g$ on $\R^3_+$.
\end{rem}

\medskip

In order to ensure that the smallness condition (\ref{datasmall:VNSreg}) is satisfied for all times, we now introduce the following terminology, which has been already used in \cite{HKMM, HK} to take advantage of the parabolic regularization for the fluid velocity.
\begin{defi}[Strong existence time]\label{strongtime}
A real number $T\geq 0$ is a \emph{strong existence} time whenever (\ref{datasmall:VNSreg}) holds.
\end{defi}

%We first have the following result.
%\begin{lem}\label{1:strongtime}
%The smallness condition (\ref{smallness:condition}) of Theorem \ref{theoreme1} ensures that $T=1$ is a strong existence time in the sense of Definition \ref{strongtime}.
%\end{lem}
%\begin{proof}
%Recall the meaning of the notation $\lesssim_0$ from Definition \ref{nota:lesssim_0}. We use Lemma \ref{ineg-Brinkman} and the local estimate (\ref{borneLinftyLinfty-rhoetj}) to write
%\begin{align*}
%\int_0^1 \Vert F(s) \Vert_{\Ld^2(\R^3_+)}^2 \, \mathrm{d}s \leq  \mathrm{E}(0) \underset{s \in [0,1]}{\sup} \Vert \rho_f (s) \Vert_{\Ld^{\infty}(\R^3_+)} \lesssim_0  \mathrm{E}(0) \leq \dfrac{1}{2C_1 \sqrt{8C_1 C_2}},
%\end{align*}
%where we have used the assumption (\ref{smallness:condition}). Combining this inequality with (\ref{Hyp:H1small}) leads to the result.
%\end{proof}

%\subsection{Local estimates for the second derivatives of the fluid velocity}\label{Section:estimateforces}

%We first introduce the following useful notations involving the moments of the kinetic distribution.
%\begin{defi}
%We set for all $t \geq 1$
%\begin{align*}
%\M_{\rho_f}(t) &:= \underset{s \in [1,t]}{\sup} \Vert \rho_f(s) \Vert_{\Ld^{\infty}(\R^3_+)}, \\ \M_{j_f}(t) &:= \underset{s \in [1,t]}{\sup} \Vert j_f(s) \Vert_{\Ld^{\infty}(\R^3_+)} ,\\
%\M_{\rho_f,j_f}(t) &:=\M_{\rho_f}(t)+\M_{j_f}(t).
%\end{align*}
%\end{defi}
In the remaining part of this section, we state a local in time $\Ld^1 \W^{1, \infty}$ regularity result for the fluid velocity. Note that for the moment, we are only interested in obtaining non-uniform in time estimates. Of course, quantitative and uniform in time estimates based on the polynomial decay of the kinetic energy will require an additional analysis.

\begin{coro}\label{BrinkLpLp}
For any finite strong existence time $t>0$ and for any $p \in [1,6]$, we have
\begin{align*}
j_f-\rho_f u \in \Ld^p(0,t;\Ld^p(\R^3_+)).
\end{align*}
\end{coro}
\begin{proof}
Thanks to Lemma \ref{rough:bounds} and \eqref{borneLinftyLinfty-rhoetj}, we have
\begin{align*}
j_f \in \Ld^{\infty}(0,t;\Ld^{3/2}(\R^3_+)) \cap \Ld^{\infty}(0,t;\Ld^{\infty}(\R^3_+)) \hookrightarrow \Ld^{r}(0,t;\Ld^{r}(\R^3_+)),
\end{align*}
for any $r \in [3/2,+\infty]$ by interpolation. For the same reasons, we also have
\begin{align*}
\rho_f \in \Ld^{\infty}(0,t;\Ld^{2}(\R^3_+)) \cap \Ld^{\infty}(0,t;\Ld^{\infty}(\R^3_+)) \hookrightarrow \Ld^{r}(0,t;\Ld^{r}(\R^3_+)),
\end{align*}
for any $r \in [2,+\infty]$. In addition, by \eqref{controlparab1}, we have $u \in \Ld^{\infty}(0,t;\Ld^q(\R^3_+))$ for any $q \in [2,6]$, because $t$ is a strong existence time. So, the Hölder's inequality yields
\begin{align*}
\rho_f u \in \Ld^p(0,t;\Ld^p(\R^3_+)), \ \ \dfrac{1}{p}=\dfrac{1}{r}+\dfrac{1}{q}, \ \ r \in [2,+\infty], \ \ q \in [2,6].
\end{align*}
This leads to the condition $p \in [1,6]$ and then to the conclusion.
\end{proof} 
 
The next result is similar to that of \cite[Lemma 3.28]{HKM}. We detail the proof for the sake of completeness, highlighting the role of Assumption \eqref{data:hyp}.
 \begin{propo}\label{Prop:convecLPLP}
Consider the exponent $s$ and $p_0$ given in Assumption \eqref{data:hyp}. For any $p \in (3, p_0)$ and any finite strong existence time $t > 0$, we have
 \begin{align}\label{convecLPLP}
(u \cdot  \nabla )u \in \Ld^p(0,t;\Ld^p(\R^3_+)).
 \end{align}
 \end{propo}
\begin{proof}
For any $(a,b,r_1,r_2) \in (1,+\infty)^4$, we can use interpolation inequalities to write
\begin{align}\label{interpo:convec}
\Vert (u \cdot  \nabla )u \Vert_{\Ld^a(0,t;\Ld^b(\R^3_+))} &\leq \Vert u \Vert_{\Ld^{\infty}(0,t;\Ld^6(\R^3_+))}  \Vert \nabla u \Vert_{\Ld^{\infty}(0,t;\Ld^2(\R^3_+))}^{1-\frac{r_1}{a}} \Vert \nabla u \Vert_{\Ld^{r_1}(0,t;\Ld^{r_2}(\R^3_+))}^{\frac{r_1}{a}},
\end{align}
provided that $r_1 \leq a$, $2 \leq b \leq r_2$ and
\begin{align}\label{interpo:convecREL}
\dfrac{1}{b}=\dfrac{2}{3}+\dfrac{r_1}{a}\left(\dfrac{1}{r_2}-\dfrac{1}{2} \right).
\end{align}
Taking $(a,b,r_1,r_2)=(2,3,2,6)$ in \eqref{interpo:convec}, the Sobolev embedding and Proposition \ref{propdatasmall:VNSreg} imply
\begin{align*}
(u \cdot  \nabla )u \in \Ld^2(0,t;\Ld^3(\R^3_+)).
\end{align*}
Owing to the maximal regularity for the Stokes system (see Section \ref{AnnexeMaxregStokes} in the Appendix) and to Corollary \ref{BrinkLpLp} which gives $j_f -\rho_f u \in \Ld^{2}(0,t;\Ld^3(\R^3_+))$, as well as on the assumption \eqref{data:hyp}, we obtain
\begin{align*}
u \in \Ld^2(0,t; \W^{2,3}(\R^3_+)).
\end{align*}
So, by the Sobolev embedding and since $\nabla u \in \Ld^2(0,t; \Ld^2(\R^3_+))$, we infer that for all $r \in [2,+\infty)$
\begin{align*}
\nabla u \in \Ld^2(0,t;\Ld^r(\R^3_+)).
\end{align*}
Coming back to the inequality \eqref{interpo:convec} with $(b,r_1)=(3,2)$, $r_2 \geq 2$ (which means $a=\widetilde{a}=3(r_2-2)r_2^{-1} \in [2,3)$), we now get 
\begin{align*}
(u \cdot  \nabla )u \in \Ld^{\widetilde{a}}(0,t;\Ld^3(\R^3_+)).
\end{align*}
By Corollary \ref{BrinkLpLp}, we also have
$j_f -\rho_f u \in \Ld^{\widetilde{a}}(0,t;\Ld^3(\R^3_+))$,
therefore an other application of the maximal regularity for the Stokes system implies that for all $\widetilde{a} \in [2,3)$
\begin{align*}
u \in \Ld^{\widetilde{a}}(0,t; \W^{2,3}(\R^3_+)),
\end{align*}
if $u_0 \in \D_3^{1-\frac{1}{\widetilde{a}},\widetilde{a}}(\R^3_+)$. So, under this assumption, the Sobolev embedding implies that for all $\widetilde{a} \in [2,3)$ and for all $r \in [2,+\infty]$, we have
\begin{align*}
\nabla u \in \Ld^{\widetilde{a}}(0,t;\Ld^r(\R^3_+)).
\end{align*}
This allows to apply the estimate \eqref{interpo:convec} with $a=b=p>3$ and $\widetilde{a}=r_1 \in [2,3)$ (and also $p \leq r_2$) so that
\begin{align*}
(u \cdot  \nabla )u \in \Ld^p(0,t;\Ld^p(\R^3_+)).
\end{align*}
The relation \eqref{interpo:convecREL} reads as
\begin{align*}
r_2=\dfrac{6 r_1}{3 r_1+6-4p}, \ \ p\leq r_2,
\end{align*}
therefore this turns into
\begin{align*}
p<\dfrac{3(r_1+2)}{4}.
\end{align*}
Since $r_1=2$ leads to the limiting case $p<3$, we can rely on the assumption \eqref{data:hyp} and use the exponents $s$ and $p_0$ (taking $r_1=s \in (2,3)$, and considering all $p \in (3,p_0)$) to conclude the proof.
\end{proof}

As in \cite{HKMM,HK,EHKM}, we finally deduce the following non-uniform in time result.

\begin{coro}\label{grad:nonunif}
For any finite strong existence time $t > 0$, we have
\begin{align*}
\nabla u \in \Ld^1(0,t;\Ld^{\infty}(\R^3_+)).
\end{align*}
\end{coro}

\begin{proof}
Let $t>0$ be a finite strong existence time. We consider the exponent $p>3$ given Proposition \ref{convecLPLP}, and which also appears in the Assumption \eqref{data:hyp}. We invoke the Gagliardo-Nirenberg-Sobolev inequality (see Theorem \ref{gagliardo-nirenberg} in the Appendix) which yields
\begin{align*}
\Vert \nabla u(s) \Vert_{\Ld^{\infty}(\R^3_+)} \lesssim \Vert \D^2 u(s) \Vert_{\Ld^{p}(\R^3_+)}^{\beta_p} \Vert u(s) \Vert_{\Ld^{2}(\R^3_+)}^{1-\beta_p}, \ \ s \in [0,t],
\end{align*}
where $\beta_p:=\dfrac{5p}{7p-6}$. Combining the energy inequality (\ref{ineq-energy}) and Lemma \ref{j_f:L^1L^1}, we get
\begin{align*}
\int_{0}^t \Vert \nabla u(s) \Vert_{\Ld^{\infty}(\R^3_+)} \, \mathrm{d}s& \lesssim \int_{0}^t \Vert \D^2 u(s) \Vert_{\Ld^{p}(\R^3_+)}^{\beta_p} \Vert u(s) \Vert_{\Ld^{2}(\R^3_+)}^{1-\beta_p} \, \mathrm{d}s \\
& \lesssim \Vert u \Vert_{\Ld^{\infty}(0,t;\Ld^{2}(\R^3_+))}^{1-\beta_p}\int_{0}^t \Vert \D^2 u(s) \Vert_{\Ld^{p}(\R^3_+)}^{\beta_p}  \, \mathrm{d}s \\
&\lesssim  \left( \E(0)+ \E(0)^{\frac{1}{2}}t + t^2  \right)^{\frac{(1-\beta_p)}{2}}\int_{0}^t \Vert \D^2 u(s) \Vert_{\Ld^{p}(\R^3_+)}^{\beta_p}  \, \mathrm{d}s \\
&\lesssim  \left( \E(0)+ \E(0)^{\frac{1}{2}}t + t^2  \right)^{\frac{(1-\beta_p)}{2}} t^{1-\frac{p}{\beta_p}} \Vert \D^2 u \Vert_{\Ld^{p}(0,t;\Ld^{p}(\R^3_+))}^{\beta_p},
\end{align*}
where we have used the Hölder's inequality in the last line. Furthermore, thanks to the maximal $\Ld^p \Ld^p$ regularity for the Stokes system, we have
\begin{align*}
\Vert \D^2 u \Vert_{\Ld^{p}(0,t;\Ld^{p}(\R^3_+))} \lesssim \Vert u_0 \Vert_{D^{1-\frac{1}{p},p}_p(\R^3_+)} + \Vert j_f-\rho_f u \Vert_{\Ld^{p}(0,t;\Ld^{p}(\R^3_+))} +\Vert (u \cdot \nabla) u\Vert_{\Ld^{p}(0,t;\Ld^{p}(\R^3_+))}<\infty,
\end{align*}
thanks to the assumption \eqref{data:hyp}, Corollary \ref{BrinkLpLp} and Proposition \ref{convecLPLP}. This allows to conclude the proof.
\end{proof}

\begin{rem}
By Corollary \ref{coro:uL1linfini} and Corollary \ref{grad:nonunif}, we get $u \in \Ld^1(0,t; \W^{1,\infty}_0(\R^3_+))$ when $t>0$ is a finite strong existence time, so that the characteristic curves for the Vlasov equation are classically defined on $(0,t)$ and the representation formula \eqref{formule-rep} can be applied.
\end{rem}
\section{Exit geometric condition and absorption}\label{Section:EGCabs}
The main goal of this section is to describe precisely the effect of the absorption boundary condition \eqref{bcond-f} satisfied by the distribution function $f$ solution to the Vlasov equation \eqref{eq:Vlasov}. In short, we will study the time of absorption when one starts from a compact support for the initial distribution function $f_0$.
The simple geometry of the flat boundary will allow us to base our study upon the characterics curves for the Vlasov equation.
%will thanks to the presence of the gravitation term and the particular geometry of the half-space. 
%This will virtually allow to provide any polynomial decay in time of the moments of $f$, modulo a decay of $f_0$: this important consequence will be written down in Section \ref{Section:Bootstrap}.

As explained in the introduction, we rely on different ideas mainly taken from the work of Glass, Han-Kwan and Moussa in \cite{GHKM} (but which lead to different types of results). In some sense, the adaptation to the half-space case is less tedious because we only deal with a flat boundary. Here, the set $\B_v(R)$ refers to $\B(0,R) \subset \R^3$ with $R \in [0,+\infty]$. 
We first introduce the so-called \textit{exit geometric condition}.
\begin{defi}\label{DEF:EGC}
Let $L,R>0$. We say that a vector field $\U:\R^+ \times \R^3_+ \rightarrow \R^3$ satisfies the exit geometric condition (EGC) in time $T$ with respect to $\big(\R^2 \times (0,L)\big) \times \B_v(R)$ if
\begin{align}\label{EGC-T}
\underset{(x,v) \in \,  (\R^2 \times (0,L)) \times \B_v(R) }{\sup} \ \tau^{+}_{\U}(0,x,v)<T,
\end{align}
where $\tau^{+}_{\U}$ refers to Definition \eqref{def:tau+} of the characteristic curves $(\mathrm{X_U},\mathrm{V_U})$ of the Vlasov equation associated to a velocity field $\U$ in \eqref{EDO-charac}, that is the solution to \begin{equation}\label{EDO-charac-grav}
\left\{
      \begin{aligned}
        \frac{\mathrm{d}}{\mathrm{d}s}\mathrm{X}_{\U}(s;t,x,v) &=\mathrm{V}_{\U}(s;t,x,v),\\
\frac{\mathrm{d}}{\mathrm{d}s}\mathrm{V}_{\U}(s;t,x,v)&= P\U(s,\mathrm{X}_{\U}(s;t,x,v)) -\mathrm{V}_{\U}(s;t,x,v)+G,\\
	\mathrm{X}_{\U}(t;t,x,v)&=x,\\
	\mathrm{V}_{\U}(t;t,x,v)&=v.
      \end{aligned}
    \right.
\end{equation}
\end{defi}

%The effect of the absorption is merely the following and we will mainly base our analysis upon this one in the next Section.
A direct consequence of an EGC satisfied by a velocity field is the following.
\begin{propo}\label{PropoABS-form}
Suppose that a velocity field $\mathrm{U} \in \Ld^1_{\mathrm{loc}}(\R^+;\mathrm{W}^{1,\infty}_0(\R^3_+))$ satisfies an EGC in time $T>0$ with respect to $\big(\R^2 \times (0,L)\big) \times \B_v(R)$ for some fixed $L, R>0$. Then, if $f$ is the solution to the Vlasov equation
\begin{equation*}
\left\{
      \begin{aligned}
      \partial_t f +v\cdot \nabla_x f + \mathrm{div}_v \, (f(\mathrm{U}-v)+fG)&=0,\\
f_{\mid t=0}&=f_0,\\
f &=0, \ \mathrm{on} \ \Sigma^{-},
      \end{aligned}
    \right.
\end{equation*}
with initial data $f_0$, we have for almost every $(x,v) \in \R^3_+ \times \R^3$ and any $t > T$
\begin{align}\label{formuleRep:split3}
\begin{split}
f(t,x,v)&= e^{3t} \mathbf{1}_{\mathcal{O}^t_{\U}}(x,v) \mathbf{1}_{\vert \V^0_t(x,v) \vert >R}   \, f_0(\X^0_t(x,v),\V^0_t(x,v)) \\[2mm]
& \quad + e^{3t} \mathbf{1}_{\mathcal{O}^t_{\U}}(x,v) \mathbf{1}_{\vert \V^0_t(x,v) \vert \leq R} \,  \mathbf{1}_{ \X^0_t(x,v)_3 >L}  \, f_0(\X^0_t(x,v),\V^0_t(x,v)),
\end{split}
\end{align}
where $(\X,\V)=(\X_{\mathrm{U}}, \V_{\mathrm{U}})$ and 
$$\mathcal{O}^t_{\U}=\left\lbrace (x,v) \in \R^3_+ \times \R^3 \mid \forall \sigma \in [0,t], \ \  \mathrm{X}_{\U}(\sigma;t,x,v) \in \R^3_+ \right\rbrace,
$$
\end{propo}
\begin{proof}
We drop the dependency un $\U$. Let $t>T$. Observe that the regularity of $u$ allows us to manipulate the characteristic curves for the Vlasov equation in a classic sense on $[0,t]$. From the representation formula \eqref{formule-rep}, we have
\begin{align*}
f(t,x,v)= e^{3t} \mathbf{1}_{\mathcal{O}^t}(x,v) \,  \mathbf{1}_{\vert \V^0_t(x,v) \vert \leq R} \,  \mathbf{1}_{ \X^0_t(x,v)_3 \leq L}  \, f_0(\X^0_t(x,v),\V^0_t(x,v)) + \widetilde{f}(t,x,v),
\end{align*}
where $\widetilde{f}(t,x,v)$ denotes the expression in the right-hand side of \eqref{formuleRep:split3}. We thus have to prove that the first term of the previous equality vanishes. Using the change of variable $(x,v) \mapsto (\X(0;t,x,v),\V(0;t,x,v))$ (see Proposition \ref{Propo:diffeoZ}) together with the fact that
\begin{align*}
(\X,\V)(0;t,\mathcal{O}^t)=\left\lbrace  (x,v) \in \R^3_+ \times \R^3  \mid  \tau^{+}(0,x,v) >t  \right\rbrace,
\end{align*}
%where
%\begin{align*}
%\mathcal{O}^t=\left\lbrace (x,v) \in \R^3_+ \times \R^3 \mid \tau^-(t,x,v)<0\right\rbrace,
%\end{align*}
we get
\begin{align*}
& \int_{\R^3_+ \times \R^3}  e^{3t} \mathbf{1}_{\mathcal{O}^t}(x,v)  \,  \mathbf{1}_{\vert \V^0_t(x,v) \vert \leq R} \,  \mathbf{1}_{ \X^0_t(x,v)_3 \leq L}  \,  f_0(\X^0_t(x,v),\V^0_t(x,v)) \, \mathrm{d}x \, \mathrm{d}v \\
& = \int_{\R^3_+ \times \R^3}  \mathbf{1}_{\tau^+(0,x,v)>t} \,  \mathbf{1}_{\vert v \vert \leq R} \, \mathbf{1}_{ x_3  \leq L} \, f_0(x,v) \, \mathrm{d}x \, \mathrm{d}v \\
& =\int_{ x_3  \leq L, \ \vert v \vert \leq R }  \mathbf{1}_{\tau^+(0,x,v)>t} \,  f_0(x,v) \, \mathrm{d}x \, \mathrm{d}v.
\end{align*}
By definition of the EGC in time $T>0$ with respect to $\big(\R^2 \times (0,L)\big) \times \B_v(R)$, we have $\tau^+(0,x,v)<T$ for all $(x,v) \in \big(\R^2 \times (0,L)\big) \times \B_v(R)$ therefore the last integral is actually zero because $t>T$. Since it is true for all $t > T$ and since $f_0$ is nonnegative, this concludes the proof.
\end{proof}

\medskip
The main idea that we follow now is to compare the Vlasov equation with velocity field $u$ (solution to the Navier-Stokes equations) to the “free" Vlasov equation without coupling. We thus consider the following characteristic curves $(\mathrm{X}^g,\mathrm{V}^g)$ for the Vlasov equation associated with the vector field $(x,v) \mapsto (v,G-v)$, namely
\begin{equation}\label{EDO-charac-grav}
\left\{
      \begin{aligned}
        \frac{\mathrm{d}}{\mathrm{d}s}\mathrm{X}^g(s;t,x,v) &=\mathrm{V}^g(s;t,x,v),\\
\frac{\mathrm{d}}{\mathrm{d}s}\mathrm{V}^g (s;t,x,v)&= G-\mathrm{V}^g(s;t,x,v),\\
	\mathrm{X}^g(t;t,x,v)&=x,\\
	\mathrm{V}^g(t;t,x,v)&=v.
      \end{aligned}
    \right.
\end{equation}
These are the equations satisfied by the characteristic curves associated to a trivial velocity field (i.e $\U =0$) in the Vlasov equation and where the particles only undergo the effect of the gravity force $G$, without being coupled to a surrounding fluid. We have the formulas
\begin{equation}\label{expr:Zt-grav}
\left\{
      \begin{aligned}
        \mathrm{X}^g(t;0,x,v)&=x+(1-e^{-t})v+(t+e^{-t}-1)G,  \\
        \mathrm{V}^g(t;0,x,v)&= e^{-t}v+(1-e^{-t})G,
      \end{aligned}
    \right.
\end{equation}
and in particular, because $G=(0,0,-g)$, we have
\begin{align*}
\left\{
      \begin{aligned}
\X^g(t;0,x,v)_3&=x_3+(1-e^{-t})v_3-(t+e^{-t}-1)g,\\
\V^g(t;0,x,v)_3&=-g+e^{-t}(v_3+g).
  \end{aligned}
    \right.
\end{align*}
\begin{rem}\label{rmq:charac:gravity}
%If $(x,v) \in (\R^2 \times (0,L)) \times \R^3$ is fixed, we observe that $\V_3^g(t;0,x,v) < 0$ for all times $t \geq 0$ when $v_3+g \geq 0$, and that $\V_3^g(t;0,x,v) < 0$ if $e^{-t} <g(g+v_3)^{-1}$ when $v_3+g >0$. In addition, it shows that if $\V_3^g(t;0,x,v) < 0$ then $\V_3^g(s;0,x,v) < 0$ for all $s \geq t$.
%We will see later that the initial EGC related to $(\mathrm{X}^g,\mathrm{V}^g)$ is always satisfied in some time $T$, for compactly supported initial data, and depends only on this initial support.
In view of the property \eqref{traj:exit}, we observe that for the characteristic curves $\mathrm{Z}^g:=(\mathrm{X}^g,\mathrm{V}^g)$, we have 
$$ \Big\lbrace (x,v) \in \big(\R^2 \times (0,L)\big) \times \R^3  \mid \tau^{+}_g(0,x,v) = \infty \ \text{or} \ \mathrm{Z}^g(\tau^{+}(0,x,v);0,x,v) \in \Sigma^0 \Big\rbrace= \emptyset,$$
where $\tau^+_g$ refers to the forward exit time associated to the curves $(\X^g,\V^g)$. This means that the curve $\mathrm{X}^g$ leaves the domain with a transversal exit. In addition, if $(x,v) \in (\R^2 \times (0,L)) \times \R^3$ is fixed and if $\V^g(t;0,x,v)_3 < 0$ then $\V^g(s;0,x,v)_3 < 0$ for all $s \geq t \geq 0$.
\end{rem}

%REMARQUE: on n'a pas mis la condition de sortie transversale dans la dèf de l'EGC. Mais pour les charac libres, c'est toujours vérifié. Cette propriété de transversalité donne de la stabilité de l'EGC sous petite perturbation: l'EGC est vérifié pour les charac totale si u petit, (et la prop de tranversalité aussi). Mais on ne se sert pas de la transversalité sur les charac totales ! Juste du premier moment où ça touche le bord (c'est parce que la théorie des solutions faibles pour Vlasov ne voit pas ce genre de chose, via les fonctions test).

\medskip

As we shall see in the end of the current section, quantitative information about the EGC are easily available for the trivial velocity field. The following stability result, which is directly inspired by \cite{GHKM}, will thus enable us to show that any velocity field solution to the Navier-Stokes equations satisfies the EGC in some finite time, provided that its $\Ld^1 \Ld^{\infty}$ norm is small enough.

\begin{lem}\label{LM:perturbEGC}
Let $\alpha>0$. There exists a constant $\kappa_{\alpha}>0$ such that the following holds. Suppose that the trivial vector field (related to $(\mathrm{X}^g,\mathrm{V}^g)$) satisfies the EGC with respect to $\left(\R^2 \times (0,L) \right) \times \B_v(0,R)$ in time $T>0$, where $L,R>0$ are given.  Then, any vector field $\U \in \Ld^1_{\mathrm{loc}}(\R^+;\mathrm{W}^{1,\infty}_0(\R^3_+))$ such that
\begin{align}\label{smallnessU-L1Linfty}
\int_0^{T+\alpha} \Vert \U(s) \Vert_{\Ld^{\infty}(\R^3_+)} \, \mathrm{d}s \leq \kappa_{\alpha},
\end{align}
satisfies the EGC in time $T+\alpha$ with respect to $\left(\R^2 \times (0,L) \right) \times \B_v(0,R)$.
\end{lem}
\begin{proof}
For all $(x,v) \in \left(\R^2 \times (0,L) \right) \times \B_v(0,R)$ and $(s,t) \in \R^+ \times \R^+$, we consider $(\X_{\U}(s,t,x,v), \V_{\U}(s,t,x,v))$ (resp. $(\X^g(s,t,x,v), \V^g(s,t,x,v))$ ) the characteristic curves associated to $\U$ (resp. to the trivial velocity field). We first set
\begin{align*}
(\mathrm{Y},\mathrm{W}):=(\mathrm{X}_{\U}-\X^g,\V_{\U}-\V^g),
\end{align*}
which satisfy the following equations
\begin{equation*}\
\left\{
      \begin{aligned}
        \frac{\mathrm{d}}{\mathrm{d}s}\mathrm{Y}(s;0,x,v) &=\mathrm{W}(s;0,x,v),\\
\frac{\mathrm{d}}{\mathrm{d}s}\mathrm{W} (s;0,x,v)&= (P \, \U)(s,\mathrm{X}_{\U}(s;0,x,v))-\mathrm{W}(s;0,x,v),\\
	\mathrm{Y}(0;0,x,v)&=0,\\
	\mathrm{W}(0;0,x,v)&=0.
      \end{aligned}
    \right.
\end{equation*}
We now fix $(x,v) \in \left(\R^2 \times (0,L) \right) \times \B_v(0,R)$. We observe that we have
for all $t \in \R^+$
\begin{equation*}
\mathrm{Y}(t;0,x,v)=\int_0^t (1-e^{\tau-t}) (P \U)(\tau, \X_{\U}(\tau,0,x,v) \, \mathrm{d}\tau,
\end{equation*}
so that for all $t \in [0,T+\alpha]$
\begin{equation}\label{ineq:stabEGC}
\vert \mathrm{Y}(t;0,x,v ) \vert \leq \int_0^t (1-e^{\tau-t}) \Vert  (P \U)(\tau) \Vert_{\Ld^{\infty}(\R^3)} \, \mathrm{d}\tau \leq \int_0^{T+\alpha}  \Vert  \U(\tau) \Vert_{\Ld^{\infty}(\R^3_+)} \, \mathrm{d}\tau,
 \end{equation}
because of the property \eqref{opP:Linfini}.

Furthermore, thanks to the EGC satisfied by the trivial velocity field in time $T$ and Remark \ref{rmq:charac:gravity}, we have $\X^g(T;0,x,v)_3<0$ and $\V^g(T;0,x,v)_3<0$. Hence, 
%\begin{align}\label{ineqX^g_3:stabEGC}
%\begin{split}
%\X^g(T+1;0,x,v)_3&=\X^g(T+1;T,\X^g(T;0,x,v)_3,\V^g(T;0,x,v)_3)_3\\
%&=\X^g(T;0,x,v)_3+(1-e^{-1})(\V^g(T;0,x,v)_3+g)-g \\
%&<-e^{-1}g<0.
%\end{split}
%\end{align}
\begin{align}\label{ineqX^g_3:stabEGC}
\begin{split}
\X^g(T+\alpha;0,x,v)_3&=\X^g(T+\alpha;T,\X^g(T;0,x,v)_3,\V^g(T;0,x,v)_3)_3\\
&=\X^g(T;0,x,v)_3+(1-e^{-\alpha})(\V^g(T;0,x,v)_3+g)-\alpha g \\
&<(1-\alpha-e^{-\alpha})g<0,
\end{split}
\end{align}
because $\alpha>0$. If we set $\eta_{\alpha}:=(e^{-\alpha}+\alpha-1)g>0$ and $\kappa_{\alpha}:=\eta_{\alpha}/2$, we see that if $\U$ satisfies the condition
\begin{align*}
\int_0^{T+\alpha} \Vert \U(s) \Vert_{\Ld^{\infty}(\R^3_+)} \, \mathrm{d}s \leq \kappa_{\alpha},
\end{align*}
then the estimate \eqref{ineq:stabEGC} turns into
\begin{align}\label{majYW}
\sup_{t \in [0 , T+\alpha ]} \vert \Y(t;0,x,v) \vert \leq \frac{\eta_{\alpha}}{2}. 
\end{align}
In view of \eqref{ineqX^g_3:stabEGC} and Remark \ref{rmq:charac:gravity}, we get the existence of $t_0=t_0(x,v) \in (0,T+\alpha]$ such that $\X^g(t_0;0,x,v)\notin \left\lbrace h \in \R^3 \mid h_3>-\eta_{\alpha} \right\rbrace$. From \eqref{majYW}, we deduce that $\X_{\U} (t_0;0,x,v)\notin \R^3_+$ and therefore $\tau^+_{\U}(0,x,v)<T+\alpha$. As it is true for any $(x,v) \in (\R^2 \times (0,L)) \times \B_v(R)$, this means by definition that the EGC is satisfied for $\U$ in time $T+\alpha$ with respect to $\left(\R^2 \times (0,L) \right) \times \B_v(0,R)$.
\end{proof}
%\begin{rem}\label{Rmk:perturbEGC}
%We observe that the result of Lemma \ref{LM:perturbEGC} remains true if one replaces $T+1$ by $T+\alpha$ for any $\alpha>0$. 
%\end{rem}

\medskip

Then, in view of the simple and explicit form of the characteristic curves $(\mathrm{X}^g,\mathrm{V}^g)$ for the Vlasov equation associated to a trivial velocity field, we can easily obtain precise information on the EGC satisfied for this velocity field. Indeed, for all $(x,v) \in \left(\R^2 \times (0,L) \right) \times \B_v(0,R)$ and for all $t \geq 0$, we have
\begin{align*}
 \X^g(t;0,x,v)_3&=x_3 + (1-e^{-t})(v_3+g)-tg \\
&<L+\vert v \vert +g -tg \\
&<L+R+g-tg.
\end{align*}
It naturally leads to the following definition and properties about the EGC for the trivial vector field and the characteristic curves $(\mathrm{X}^g,\mathrm{V}^g)$.
\begin{defi}\label{def:t_0(R)}
We set
\begin{align}
t_0(L,R)&:= \dfrac{L+R+g}{g}, \\
t_0&:=t_0(1,1).
\end{align}
\end{defi}
\begin{lem}
If $L,R$ are given, the trivial vector field $\mathrm{U} \equiv 0$ (associated to $(\mathrm{X}^g,\mathrm{V}^g)$) satisfies the EGC in time $t_0(L,R)$ with respect to $\left(\R^2 \times (0,L) \right) \times \B_v(R)$.
\end{lem}
\begin{proof}
Assume that  $\X^g(t,0,x,v)_3>0$ for all $t \in [0,t_0(L,R)]$ and $(x,v) \in \left(\R^2 \times (0,L) \right) \times \B_v(0,R)$. Then
\begin{align*}
0&<x_3 + (1-e^{-t})(v_3+g)-tg \\
&<L+R+g-tg \\
&=g(t_0(L,R)-t),
\end{align*}
so that we get a contradiction by taking $t=t_0(L,R)$.
\end{proof}

\begin{rem}
%We note that $t_0=t_0(1,1)$ only depends on the spatial support of $f_0$ (namely, on $L$). This means that once $t_0$ is fixed, we can't assume that $L$ is taken small. However, we can use a smallness assumption as \eqref{smallness:condition:THM} about the initial data which doesn't take into account the support of $f_0$ (see Notation \ref{nota:lesssim_0}).
We note that $t_0=t_0(1,1)$ does not depend on the initial data. This time only has to be seen as a reference time after which we will use the absorption phenomenon. The subsequent analysis could have been performed by replacing $t_0(1,1)$ by $t_0(\kappa,\kappa)$ for any $\kappa>0$.
\end{rem}

\medskip

The following result is, in some sense, of reverse nature: given a time $t>t_0$, we describe which proportion of the initial velocities will lead to absorption before time $t$. More precisely, we state the following lemma.
\begin{lem}\label{EGCt:reverse}
There exist some continuous increasing functions $\Ld_g:[t_0, + \infty) \rightarrow \R^+$ and $\mathrm{R}_g: [t_0, + \infty) \rightarrow \R^+$ such that for all $t>t_0$ the trivial velocity vector field (associated to ($\X^g,\V^g)$) satisfies the EGC in time $t$ with respect to $\left( \R^2 \times (0,1+\Ld_g(t)) \right) \times \B_v(0,1+\mathrm{R}_g(t))$. Furthermore, there exist $C_{g},\underline{C}_g, \overline{C}_g>0$ such that for all $s>t_0$
\begin{align*}
\dfrac{1}{1+\Ld_{g}(s)} \leq \dfrac{C_{g}}{1+s}, \\ 
\dfrac{\underline{C}_g}{1+s} \leq \dfrac{1}{1+\mathrm{R}_{g}(s)} \leq \dfrac{\overline{C}_g}{1+s}.
\end{align*}
\end{lem}
\begin{proof}
%As an initial EGC condition holds for the trivial vector field with respect to $\left( \R^2 \times (0,L) \right) \times \B(0,1)$ in time $t_0$, this initial EGC remains valid with respect to $\left( \R^2 \times (0,L) \right) \times \B_v(0,1)$ in time $t>t_0$. Hence, we only have to prove that the trivial velocity vector field satisfies the initial EGC in time $t$ with respect to $\left( \R^2 \times (0,L) \right) \times \mathscr{C}_v(0;1,1+\phi_{L,g}(t))$, where $\mathscr{C}_v(0;1,1+\phi_{L,g}(t))$ stands for the annulus in velocity
%\begin{align*}
%\mathscr{C}_v(0;a,b):= \left\lbrace  v \in \R^3 \mid a <\vert v \vert <b) \right\rbrace.
%\end{align*}
Since the characteristic curve $\X^g$ is a vertical line, a necessary and sufficient condition ensuring an EGC in time $t$ with respect to $\left( \R^2 \times (0,1+\Ld_g(t)) \right) \times \B_v(0,1+\mathrm{R}_g(t))$ (for some positive functions $\Ld_g$ and $\mathrm{R}_g$ to be determined) is that for all $(x,v) \in \left( \R^2 \times (0,1+\Ld_g(t)) \right) \times \B_v(0,1+\mathrm{R}_g(t))$, the following inequality holds:
\begin{align}
\X^g(t;0,x,v)_3=x_3+(1-e^{-t})(v_3+g)-tg<0.
\end{align}
%Indeed, one can easily show that the function $s \mapsto \X^g(s;0,x,v)_3$ is strictly decreasing on $[\overline{t},t]$, for some time $\overline{t} \geq t_0$ which is independent of $(x,v)$.
Indeed, one can easily show that the real function $s \mapsto \X^g(s;0,x,v)_3$ becomes strictly decreasing after its first cancellation. We then set
\begin{align*}
\Ld_g(s):=\dfrac{1}{2}(sg-(1-e^{-s})g)-1, \ \ \mathrm{R}_g(s):=\dfrac{1}{2}\left(\dfrac{sg}{1-e^{-s}} -g\right)-1.
\end{align*}
Now, if $0<x_3<1+\Ld_g(t)$ and $\vert v \vert <1+\mathrm{R}_g(t)$, we observe that 
\begin{align*}
x_3+(1-e^{-t})v_3 & \leq x_3+(1-e^{-t})\vert v \vert < tg-(1-e^{-t})g,
\end{align*}
therefore we have $\X^g(t;0,x,v)_3<0$. It remains to show that $\Ld_g$ and $\mathrm{R}_g$ are positive on $[t_0, + \infty[$. As they are nondecreasing functions, we only have to prove that $\Ld_g(t_0) > 0$ and 
$\mathrm{R}_g(t_0)>0$: recalling the explicit Definition \ref{def:t_0(R)} of the time $t_0=1+2g^{-1}$, we have
\begin{align*}
\Ld_g(t_0)&=\dfrac{1}{2}\left((1+\frac{2}{g})g-(1-e^{-(1+\frac{2}{g})})g\right)-1=\frac{1}{2}e^{-(1+\frac{2}{g})})g>0,\\
\mathrm{R}_g(t_0)&=\dfrac{1}{2}\left(\dfrac{(1+\frac{2}{g})g}{1-e^{-(1+\frac{2}{g})}} -g\right)-1=\dfrac{2}{1-e^{-(1+\frac{2}{g})}}-1+ g \left( \dfrac{1}{1-e^{-(1+\frac{2}{g})}}-\dfrac{1}{2} \right)>1+ \frac{g}{2}>0,
\end{align*}
which is the desired claim. Concerning the last part, a direct computation shows that for all $s \geq t_0$, we have
\begin{align*}
\dfrac{1+s}{1+\Ld_g(s)} \leq \dfrac{1+t_0}{1+\Ld_g(t_0)}, \ \ \dfrac{1+s}{1+\mathrm{R}_g(s)} \leq \dfrac{1+t_0}{1+\mathrm{R}_g(t_0)},
\end{align*}
and that $s \mapsto \dfrac{1+\mathrm{R}_g(s)}{1+s}$ is bounded from above on $[t_0,+\infty)$. The proof is then complete.
\end{proof}

\section{The bootstrap argument}\label{Section:Bootstrap}
In this section, we provide a proof of Theorem \ref{thm1}, relying on the absorption effect highlighted in Section \ref{Section:EGCabs}. Our strategy is based on a bootstrap argument reminiscent of the ideas of \cite{HKMM,HK}.
Roughly speaking, we will prove that one can propagate the controls \eqref{borne-gradient0} and \eqref{smallnessU-L1Linfty} for the velocity field $u$.
%so that, at each time $t$, an EGC is satisfied with respect to some compact in space and velocity (evolving linearly with $t$) and the change of variable $v \mapsto \V(0;t,x,v)$ is admissible. This will implies a sufficient (polynomial) decay of the moments of $f$ and then of the Brinkman force.
\subsection{Initialization of the bootstrap procedure}
In order to set up a boostrap procedure, we introduce the following quantities.
\begin{defi}
We set $T_0:=t_0+1$ where  $t_0$ is given in Definition \ref{def:t_0(R)}.
\end{defi}
\begin{defi}
We consider $\delta_0>0$ which satisfies $\delta_0 e^{\delta_0}<1/9$ and $\delta_0<\kappa_{1/2}$, where $\kappa_{1/2}$ is given in Lemma \ref{LM:perturbEGC}.
\end{defi}

Let $(u,f)$ be a global weak solution to the Vlasov-Navier-Stokes system in the sense of Definition \ref{sol-faible} and associated to an admissible initial data $(u_0,f_0)$ satisfying \eqref{data:hyp}. We start with the following lemma.
\begin{lem}\label{T=T_0strong}
We have
\begin{align}\label{ineq:BrinkL2L2:loc}
 \int_0^{T_0} \Vert j_f(s) -\rho_f u(s) \Vert_{\Ld^2(\R^3_+)}^2  \, \mathrm{d}s \lesssim N_q(f_0) \eta(T_0)\left[ \E(0) + T_0^2 \right],
\end{align}
for some nondecreasing positive continuous function $\eta$ and the exponent $q$ appearing in Assumption \eqref{data:hyp}. Furthermore, under the smallness assumption \eqref{smallness:condition:THM}, the time $T_0$ is a strong existence time, in the sense of Definition \ref{strongtime}.
\end{lem}
\begin{proof}
We first set $F:= j_f -\rho_f u$. From Lemma \ref{ineg-Brinkman} and \eqref{borneLinftyLinfty-rhoetj} in Proposition \ref{controlLinfiniLOC}, we infer that
\begin{align*}
 \int_0^{T_0} \Vert F(s) \Vert_{\Ld^2(\R^3_+)}^2  \, \mathrm{d}s & \leq \underset{s \in [0,T_0]}{\sup} \Vert \rho_f (s) \Vert_{\Ld^{\infty}(\R^3_+)}  \left[ \E(0) + \E(0)^{1/2}T_0 + T_0^2 \right] \\
 & \lesssim N_q(f_0) \eta(T_0)\left[ \E(0) +  T_0^2 \right].
\end{align*}
Thanks to the Cauchy-Schwarz inequality, we also have
\begin{align*}
\int_0^{T_0} \Vert F(s) \Vert_{\Ld^2(\R^3_+)} \, \mathrm{d}s \leq {T_0}^{1/2} \left( \int_0^{T_0} \Vert F(s) \Vert_{\Ld^2(\R^3_+)}^2 \, \mathrm{d}s \right)^{1/2},
\end{align*}
therefore the inequality \eqref{ineq:BrinkL2L2:loc} now entails
\begin{align*}
\Vert  u_0 \Vert_{\H^1(\R^3_+)}^2 +  \int_0^{T_0} \left[ \Vert F(s) \Vert_{\Ld^2(\R^3_+)}^2 + \Vert F(s) \Vert_{\Ld^2(\R^3_+)} \right] \, \mathrm{d}s 
%&\lesssim \Vert  u_0 \Vert_{\H^1(\R^3_+)}^2 +  \underset{s \in [0,T_0]}{\sup} \Vert \rho_f (s) \Vert_{\Ld^{\infty}(\R^3_+)}  \left[ \E(0) + \E(0)T_0 + T_0^2 \right] \\
%& \quad +  \underset{s \in [0,t]}{\sup} \Vert \rho_f (s) \Vert_{\Ld^{\infty}(\R^3_+)}^{1/2}  \left[ \E(0) + \E(0)T_0 + T_0^2 \right]^{1/2} \\
&\lesssim \Vert  u_0 \Vert_{\H^1(\R^3_+)}^2+ N_q(f_0)\eta(T_0) \left[ \E(0) + \E(0)T_0 + T_0^2 \right]\\
& \quad +  N_q(f_0)^{1/2} {T_0}^{1/2}\eta(T_0)^{1/2} \left[ \E(0) + T_0^2 \right]^{1/2} .
\end{align*}
If $\mathrm{C}_{\star}$ refers to the universal constant given in Proposition \ref{propdatasmall:VNSreg}, we can thus use the smallness assumption \eqref{smallness:condition:THM} to  ensure that
\begin{align*}
\Vert  u_0 \Vert_{\H^1(\R^3_+)}^2 +  \int_0^{T_0} \left[ \Vert F(s) \Vert_{\Ld^2(\R^3_+)}^2 + \Vert F(s) \Vert_{\Ld^2(\R^3_+)} \right] \, \mathrm{d}s< \mathrm{C}_{\star},
\end{align*}
which means that $T_0$ is a strong existence time.
\end{proof}

\begin{coro}\label{coro:uL1linfini}
Under the smallness assumption \eqref{smallness:condition:THM}, we have
\begin{align}\label{uL1linfini}
\int_{0}^{T_0} \Vert u(s) \Vert_{\Ld^{\infty}(\R^3_+)}  \, \mathrm{d}s < \dfrac{\delta_0}{2}.
\end{align}
\end{coro}
\begin{proof}
By Lemma \ref{T=T_0strong}, we know that $T_0$ is a strong existence time therefore the parabolic regularization of the Navier-Stokes equations stated in Proposition \ref{propdatasmall:VNSreg} holds for $u$ on $[0,T_0]$. Namely, we get $$u \in \Ld^{\infty}(0,T_0;\H^1_{\mathrm{div}}(\R^3_+)) \cap \Ld^{2}(0,T_0;\H^2(\R^3_+)),$$ and there exists $\widetilde{C}>0$ such that for all $t \in [0,T_0]$
\begin{align}\label{ineq:strongNS-uloc}
\Vert \nabla u(t) \Vert_{\Ld^2(\R^3_+)}^2 + \int_0^t \Vert \D^2 u(s) \Vert_{\Ld^2(\R^3_+)}^2 \, \mathrm{d}s \leq  \widetilde{C} \left( \Vert u_0 \Vert_{\H^1(\R^3_+)}^2 +  \Vert j_f - \rho_f u \Vert_{\Ld^2(0,T_0;\Ld^2(\R^3_+))}^2\right).
\end{align}

We then use the Gagliardo-Nirenberg-Sobolev inequality (see Theorem \ref{gagliardo-nirenberg} in the Appendix) with $(p,r,q,j,m)=(\infty,2,6,0,2)$ so that there exists a universal constant $C_0>0$ such that 
\begin{align*}
\Vert u(s) \Vert_{\Ld^{\infty}(\R^3_+)} \leq C_0 \Vert \mathrm{D}^2(s) \Vert_{\Ld^{2}(\R^3_+)}^{1/2} \Vert u(s) \Vert_{\Ld^6(\R^3_+)}^{1/2}.
\end{align*}
By Sobolev's embedding, we infer from \eqref{ineq:strongNS-uloc} that
\begin{align*}
\int_0^{T_0} \Vert u(s) \Vert_{\Ld^{\infty}(\R^3_+)} \, \mathrm{d}s & \leq C_0 \int_0^{T_0} \Vert \D^2 u(s) \Vert_{\Ld^2(\R^3_+)}^{1/2}   \Vert \nabla u(s) \Vert_{\Ld^2(\R^3_+)}^{1/2}\, \mathrm{d}s \\
& \leq  C_0 T_0^{3/4} \Vert \nabla u\Vert_{\Ld^{\infty}(0,T_0;\Ld^2(\R^3_+))}^{1/2} \left( \int_0^{T_0} \Vert \D^2 u(s) \Vert_{\Ld^2(\R^3_+)}^2 \, \mathrm{d}s \right)^{1/4}  \\
& \leq C_0 T_0^{3/4} \widetilde{C}^{1/16} \left( \Vert u_0 \Vert_{\H^1(\R^3_+)}^2 +  \Vert j_f - \rho_f u \Vert_{\Ld^2(0,T_0;\Ld^2(\R^3_+))}^2\right)^{1/16} \\
& \lesssim C_0 T_0^{3/4} \widetilde{C}^{1/16} \left( \Vert u_0 \Vert_{\H^1(\R^3_+)}^2 +  N_q(f_0) \eta(T_0)\left[ \E(0) + \E(0)T_0 + T_0^2 \right]\right)^{1/16},
\end{align*}
thanks to \eqref{ineq:BrinkL2L2:loc}. We can now proceed exactly in the same way as in the proof of Lemma \ref{T=T_0strong}: since $T_0$ is fixed, we can use the smallness condition \eqref{smallness:condition:THM} to ensure that
\begin{align*}
\int_0^{T_0} \Vert u(s) \Vert_{\Ld^{\infty}(\R^3_+)} \, \mathrm{d}s< \dfrac{\delta_0}{2}.
\end{align*}
This concludes the proof of the corollary.
\end{proof}

\bigskip

In order to set up the bootstrap argument, we introduce
\begin{align}\label{def:tstar}
t^{\star} :=\sup \left\lbrace \text{strong existence times } t \geq T_0  \text{ such that } \int_{T_0} ^t \Vert \nabla u(s) \Vert_{\Ld^{\infty}(\R^3_+)}  \, \mathrm{d}s < \delta_0, \ \ \int_{T_0} ^t \Vert u(s) \Vert_{\Ld^{\infty}(\R^3_+)}  \, \mathrm{d}s < \dfrac{\delta_0}{2} \right\rbrace.
\end{align}
Our main goal is now to show that $t^{\star} = + \infty$.

\begin{lem}
We have $t^{\star}>T_0$.
\end{lem}
\begin{proof}
In view of Lemma \ref{T=T_0strong}, we see that using one more time the smallness assumption \eqref{smallness:condition:THM}, the time $t=T_0+\varepsilon_0$ is a strong existence time, if $\varepsilon_0>0$ is fixed. Now, we observe that a continuity argument, Proposition \ref{grad:nonunif} and \eqref{borneL1Linfty-u} ensure that
\begin{align*}
\int_{T_0} ^{T_0+\varepsilon_0} \Vert \nabla u(s) \Vert_{\Ld^{\infty}(\R^3_+)} \, \mathrm{d}s < \delta_0, \ \ \int_{T_0} ^{T_0+\varepsilon_0} \Vert  u(s) \Vert_{\Ld^{\infty}(\R^3_+)} \, \mathrm{d}s < \dfrac{\delta_0}{2},
\end{align*}
if $\varepsilon_0$ is small enough. By definition of $t^{\star}$, this concludes the proof.
\end{proof}
We then argue by contradiction and shall assume from now on that $t^{\star} <\infty$.

\subsection{Absorption and decay in time of the moments}
We eventually explain how one can take advantage of the absorption effect at the boundary. Recall first the concept of EGC introduced in Definition \ref{DEF:EGC}. In view of the bounds satisfied before $t^{\star}$, the vector field $u$ will satisfy an EGC at time $t<t^{\star}$ with respect to some compact depending on $t$.
%(evolving linearly with $t$) while the change of variable $v \mapsto \V(0;t,x,v)$ is admissible. At time $t$, we can thus deduce the support of $f(t)$ from the support of $f_0$, through the image of the previous compact set by the flow $(\X,\V)$. 
This will imply a (polynomial) decay in time of the moments of $f$ and then of the Brinkman force $j_f-\rho_f u$.
%
%\begin{nota}
%We will denote by $\mathrm{B}_{\mathrm{decay}}$ a generic constant which refers to an expression of the type $$\varphi \Big(  1+N_q f_0 + K_{p,r}(f_0) + \mathrm{E}(0) \Big)$$ where a finite number of exponents $q,p,r$ are involved and where  $\varphi : \R^+ \rightarrow \R^+$ is onto, continuous and nondecreasing, which may vary from line to line.
%\end{nota} 
%

Recalling the notations introduced in Definition \ref{notation:moments} and Definition \ref{decay_f0}, we state the following fundamental lemma which highlights the precise link between absorption and decay.
\begin{lem}\label{LM:decaymom-abs}
%Let $t \in (T_0,t^{\star})$. If
%\begin{align}
%\int_{0}^T \Vert \nabla u(s) \Vert_{\Ld^{\infty}(\R^3_+)} \, \mathrm{d}s < \delta_1 \ \ \ \ \int_{0}^T \Vert  u(s) \Vert_{\Ld^{\infty}(\R^3_+)} \, \mathrm{d}s < \delta,
%\end{align}
%with $\delta$ independent of $T$ and
Let $t \in (T_0,t^{\star})$ and $r \in [1, \infty)$. Suppose that $k_1,k_2,\ell \in \R^+$ and $q>3$ satisfy
\begin{align*}
q>k_1+\ell +3,
\end{align*}
Assume that the velocity field $u$ satisfies an EGC in time $t$ with respect to $\big( \R^2 \times (0,1+\mathrm{L}(t) ) \big)\times \B_v(0,1+\mathrm{R}(t))$ for some nondecreasing functions $\mathrm{L}: \R^+ \rightarrow \R^+$ and $\mathrm{R}: \R^+ \rightarrow \R^+$. Then we have for almost every $x \in \R^3_+$
%\begin{align}
%m_{\ell} f(t,x) &  \lesssim  \dfrac{N_q(f_0)}{(1+\mathrm{R}(t))^{k_1}}+ H_{k_2}(f_0) \left[\dfrac{(1+\mathrm{R}(t))^{\ell+3}}{(1+\mathrm{L}(t))^{k_2} }+ \dfrac{(1+\mathrm{R}(t))^{3}}{{(1+\mathrm{L}(t))^{k_2} }} \right],\label{decay:pointmom-GEN}\\[3mm]
%%M_{\ell}f(t) & \lesssim \dfrac{K_{q,1}(f_0)}{(1+\mathrm{R}(t))^k}, \label{decay:pointMOM} \\[2mm]
%\Vert m_{\ell} f(t) \Vert_{\Ld^r(\R^3_+)} & \lesssim \dfrac{K_{q,r}(f_0)}{(1+\mathrm{R}(t))^{k_1}}+ F_{k_2,r}(f_0) \left[\dfrac{(1+\mathrm{R}(t))^{\ell+3}}{(1+\mathrm{L}(t))^{k_2} }+ \dfrac{(1+\mathrm{R}(t))^{3}}{{(1+\mathrm{L}(t))^{k_2} }} \right], \label{decay:Lpmom-GEN}
%\end{align}
\begin{align}
m_{\ell} f(t,x) &  \lesssim  \dfrac{N_q(f_0)}{(1+\mathrm{R}(t))^{k_1}}+ \dfrac{H_{\ell,k_2}(f_0)}{{(1+\mathrm{L}(t))^{k_2} }},\label{decay:pointmom-GEN}\\[3mm]
%M_{\ell}f(t) & \lesssim \dfrac{K_{q,1}(f_0)}{(1+\mathrm{R}(t))^k}, \label{decay:pointMOM} \\[2mm]
\Vert m_{\ell} f(t) \Vert_{\Ld^r(\R^3_+)} & \lesssim \dfrac{K_{q,r}(f_0)}{(1+\mathrm{R}(t))^{k_1}}+ \dfrac{F_{\ell,k_2,r}(f_0)}{{(1+\mathrm{L}(t))^{k_2} }}, \label{decay:Lpmom-GEN}
\end{align}
where $\lesssim$ only depends on $k_1, k_2, q,\ell,g,\delta_0,r$.
\end{lem}
\begin{proof}
%The whole proof is based on the fact that
%\begin{align}
%f(t,x,v)= e^{3t} \mathbf{1}_{\tau^-(t,x,v)<0} \,  \mathbf{1}_{\vert \V^0_t(x,v) \vert >1+\mathrm{R}(t)} f_0(\X^0_t(x,v),\V^0_t(x,v)).
%\end{align}
%for almost every $(x,v) \in \R^3_+ \times \R^3$. Indeed, one can write
%\begin{multline*}
%f(t,x,v)= e^{3t} \mathbf{1}_{\tau^-(t,x,v)<0} \,  \mathbf{1}_{\vert \V^0_t(x,v) \vert \leq 1+\mathrm{R}(t)} f_0(\X^0_t(x,v),\V^0_t(x,v)) \\ + e^{3t} \mathbf{1}_{\tau^-(t,x,v)<0} \,  \mathbf{1}_{\vert \V^0_t(x,v) \vert >1+\mathrm{R}(t)} f_0(\X^0_t(x,v),\V^0_t(x,v)).
%\end{multline*}
%and perform the same argument as in the proof of Theorem \ref{thm:absorb:cpct}, that is 
%\begin{multline*}
%\int_{\R^3_+ \times \R^3} e^{3t} \mathbf{1}_{\tau^-(t,x,v)<0} \,  \mathbf{1}_{\vert \V^0_t(x,v) \vert \leq 1+\mathrm{R}(t)} \,  f_0(\X^0_t(x,v),\V^0_t(x,v)) \, \mathrm{d}x \, \mathrm{d}v \\
%=\int_{\R^3_+ \times \R^3}  \mathbf{1}_{\tau^+(0,x,v)>t} \,  \mathbf{1}_{\vert v \vert \leq 1+\mathrm{R}(t)} \, f_0(x,v) \, \mathrm{d}x \, \mathrm{d}v \\
%=\int_{(\mathrm{supp} \, f_0)_x \times \B_v(0,1+\mathrm{R}(t))}  \mathbf{1}_{\tau^+(0,x,v)>t} \,  f_0(x,v) \, \mathrm{d}x \, \mathrm{d}v =0,
%\end{multline*}
The whole proof is based on the fact that
\begin{align}\label{frep:separate}
f(t,x,v)=f^{\natural}(t,x,v)+f^{\flat}(t,x,v),
\end{align}
where 
\begin{align}
\label{eq:f-beca}f^{\natural}(t,x,v)&:=e^{3t} \mathbf{1}_{\mathcal{O}^t}(x,v)\,  \mathbf{1}_{\vert \V^0_t(x,v) \vert >1+\mathrm{R}(t)}  \,  f_0(\X^0_t(x,v),\V^0_t(x,v)), \\
\label{eq:f-bemo}f^{\flat}(t,x,v)&:=e^{3t}\mathbf{1}_{\mathcal{O}^t}(x,v)\,  \mathbf{1}_{\vert \V^0_t(x,v) \vert \leq 1+\mathrm{R}(t)} \,  \mathbf{1}_{ \X^0_t(x,v)_3 >1+\mathrm{L}(t)}  \,  f_0(\X^0_t(x,v),\V^0_t(x,v)).
\end{align}
Indeed, because of the EGC satisfied by $u$ in time $t$ with respect to $\big( \R^2 \times (0,1+\mathrm{L}(t) ) \big)\times \B_v(0,1+\mathrm{R}(t))$, we can apply Proposition \ref{PropoABS-form} so that \eqref{frep:separate} holds. The decay of the moments of $f^{\natural}$ and $f^{\flat}$ is then studied separetely.

$\bullet$ For the term $f^{\natural}$ defined in \eqref{eq:f-beca}, we use the admissible change of variable $v \mapsto \V^0_t(x,v)=\Gamma_{t,x}(v)$ given by Lemma \ref{chgmt-var-prepa}, whose Jacobian inverse is bounded by $2 e^{-3t}$. We get
\begin{align*}
m_{\ell} f^{\natural}(t,x)&= \int_{\R^3} \vert v \vert^{\ell} f^{\natural}(t,x,v) \, \mathrm{d}v \\
&=\int_{\R^3} \vert v \vert^{\ell} e^{3t} \mathbf{1}_{\mathcal{O}^t}(x,v) \,  \mathbf{1}_{\vert \V^0_t(x,v) \vert >1+\mathrm{R}(t)} \,  f_0(\X^0_t(x,v),\V^0_t(x,v)) \,  \mathrm{d}v \\[2mm]
&=\int_{\R^3} \vert \Gamma_{t,x}^{-1}(w) \vert^{\ell} e^{3t}  \mathbf{1}_{\mathcal{O}^t}(x,\Gamma_{t,x}^{-1}(w)) \mathbf{1}_{\vert w \vert >1+\mathrm{R}(t)}  \, f_0(\Lambda_{t,w}(x),w) \vert \det \, \mathrm{D}_w \Gamma_{t,x}^{-1} (w)  \vert \, \mathrm{d}w \\
& \lesssim \dfrac{1}{(1+\mathrm{R}(t))^{k_1} }\int_{\R^3} \vert w \vert ^{k_1} \vert \Gamma_{t,x}^{-1}(w) \vert^{\ell}  \mathbf{1}_{\mathcal{O}^t}(x,\Gamma_{t,x}^{-1}(w))   \, f_0(\Lambda_{t,w}(x),w) \, \mathrm{d}w,
\end{align*}
where we have used the notations of Lemma \ref{chgmt:varSPACE}, namely $\Lambda_{t,w}(x)=\X^{0}_t(x, \Gamma_{t,x}^{-1}(w))$. Thanks to the inequality \eqref{ineqGamma-1}, we also have
\begin{align*}
\vert \Gamma_{t,x}^{-1}(w) \vert &\leq  \vert w\vert +(1-e^{-t}) g + \int_0^t e^{\tau-t} \Vert Pu(\tau) \Vert_{\Ld^{\infty}(\R^3)} \, \mathrm{d}\tau, \\
& \leq  \vert w \vert +g+ \delta_0,
\end{align*}
and because of the Definition \eqref{def:Ot} of $\mathcal{O}^t$, we see that if $(x, \Gamma_{t,x}^{-1}(w)) \in \mathcal{O}^t$ then $  \Lambda_{t,w}(x) \in \R^3_+$. Therefore we obtain 
\begin{align*}
m_{\ell} f^{\natural} (t,x)& \lesssim \dfrac{1}{(1+\mathrm{R}(t))^{k_1} }\int_{\R^3} \left[ \vert w \vert^{k_1+\ell}+\vert w \vert ^{k_1} \right]  \mathbf{1}_{\mathcal{O}^t}(x,\Gamma_{t,x}^{-1}(w))   \, f_0(\Lambda_{t,w}(x),w) \, \mathrm{d}w \\
&=\dfrac{1}{(1+\mathrm{R}(t))^{k_1} }\int_{\R^3} \dfrac{ \vert w \vert^{k_1+\ell}+\vert w \vert ^{k_1} }{1+\vert w \vert^q}  \mathbf{1}_{\mathcal{O}^t}(x,\Gamma_{t,x}^{-1}(w))  (1+\vert w \vert^q)  \, f_0(\Lambda_{t,w}(x),w) \, \mathrm{d}w \\
& \leq \dfrac{N_q(f_0)}{(1+\mathrm{R}(t))^{k_1} }\int_{\R^3} \dfrac{ \vert w \vert^{k_1+\ell}+\vert w \vert ^{k_1} }{1+\vert w \vert^q}  \, \mathrm{d}w \\
&\lesssim \dfrac{N_q(f_0)}{(1+\mathrm{R}(t))^{k_1} },
\end{align*}
for $q>k_1+\ell+3$. This entails the contribution from $f^{\natural}$ in the inequality \eqref{decay:pointmom-GEN}. Concerning its contribution to the inequality \eqref{decay:Lpmom-GEN}, we combine the previous change of variable in velocity with the Minkowski's integral inequality (see \cite[Theorem 202]{HLP}) and the bound \eqref{ineqGamma-1} to write with the notations of Lemma \ref{chgmt:varSPACE}
\begin{align*}
\Vert m_{\ell} f^{\natural}(t) \Vert_{\Ld^r(\R^3_+)} &= \left[ \int_{\R^3_+} \left(\int_{\R^3} \vert v \vert^{\ell} f^{\natural}(t,x,v) \, \mathrm{d}v \right)^r \hspace{-2mm} \mathrm{d}x \right]^{1/r} \\
&\lesssim \dfrac{1}{(1+\mathrm{R}(t))^{k_1} }  \left[ \int_{\R^3_+} \left(\int_{\R^3} \vert w \vert ^{k_1} \vert \Gamma_{t,x}^{-1}(w) \vert^{\ell}  \mathbf{1}_{\mathcal{O}^t}(x,\Gamma_{t,x}^{-1}(w))   \, f_0(\Lambda_{t,w}(x),w) \, \mathrm{d}w \right)^r \hspace{-2mm} \mathrm{d}x \right]^{1/r} \\
&\leq  \dfrac{1}{(1+\mathrm{R}(t))^{k_1} } \int_{\R^3}  \vert w \vert^{k_1}  \left(\int_{\R^3_+} \vert \Gamma_{t,x}^{-1}(w)) \vert^{r \ell} \mathbf{1}_{\mathcal{O}^t}(x,\Gamma_{t,x}^{-1}(w))    \, f_0(\Lambda_{t,w}(x),w)^r  \, \mathrm{d}x \right)^{1/r} \hspace{-4mm} \mathrm{d}w \\
&  \leq \dfrac{1}{(1+\mathrm{R}(t))^{k_1} } \int_{\R^3}  \left[ \vert w \vert^{k_1+\ell}+\vert w \vert ^{k_1} \right]
  \left(\int_{\R^3} \mathbf{1}_{x \in \R^3_+} \mathbf{1}_{\mathcal{O}^t}(x,\Gamma_{t,x}^{-1}(w))    \, f_0(\Lambda_{t,w}(x),w)^r  \, \mathrm{d}x \right)^{1/r} \hspace{-4mm} \mathrm{d}w.
\end{align*}
Performing the admissible change of variable $x \mapsto \Lambda_{t,w}(x)$ in the interior integral, with a Jacobian inverse bounded by $2$ (see Lemma \ref{chgmt:varSPACE}), we get
\begin{align*}
\Vert m_{\ell} f^{\natural}(t) \Vert_{\Ld^r(\R^3_+)} & \lesssim  \dfrac{1}{(1+\mathrm{R}(t))^{k_1} } \int_{\R^3}  \left[ \vert w \vert^{k_1+\ell}+\vert w \vert ^{k_1} \right] \\
 &  \qquad \qquad \qquad \qquad \qquad \left(\int_{\R^3} \mathbf{1}_{\Lambda_{t,w}^{-1}(x) \in \R^3_+}    \mathbf{1}_{\mathcal{O}^t}\left(\Lambda_{t,w}^{-1}(x),\Gamma_{t,\Lambda_{t,w}^{-1}(x)}^{-1}(w)\right)    \, f_0(x,w)^r  \, \mathrm{d}x \right)^{1/r} \hspace{-4mm} \mathrm{d}w \\
 &\lesssim \dfrac{1}{(1+\mathrm{R}(t))^{k_1} } \int_{\R^3}  \dfrac{ \vert w \vert^{k_1+\ell}+\vert w \vert ^{k_1} }{1+\vert w \vert^q} (1+\vert w \vert^q) \Vert f_0(\cdot,w) \Vert_{\Ld^r(\R^3_+)}       \, \mathrm{d}w \\
& \leq \dfrac{K_{q,r}(f_0)}{(1+\mathrm{R}(t))^{k_1} }\int_{\R^3} \dfrac{ \vert w \vert^{k_1+\ell}+\vert w \vert ^{k_1} }{1+\vert w \vert^q}  \, \mathrm{d}w \\
& \lesssim \dfrac{K_{q,r}(f_0)}{(1+\mathrm{R}(t))^{k_1} },
\end{align*}
for $q>k_1+\ell+3$.

$\bullet$ For the term $f^{\flat}$ defined in \eqref{eq:f-bemo}, we also perform the change of variable $v \mapsto \V^0_t(x,v)=\Gamma_{t,x}(v)$ which entails
\begin{align*}
m_{\ell} f^{\flat}(t,x)&= \int_{\R^3} \vert v \vert^{\ell} f^{\flat}(t,x,v) \, \mathrm{d}v \\
&=\int_{\R^3} \vert v \vert^{\ell} e^{3t} \mathbf{1}_{\mathcal{O}^t}(x,v) \,  \mathbf{1}_{\vert \V^0_t(x,v) \vert \leq 1+\mathrm{R}(t)} \,  \mathbf{1}_{ \X^0_t(x,v)_3 >1+\mathrm{L}(t)}  \,  f_0(\X^0_t(x,v),\V^0_t(x,v)) \,  \mathrm{d}v \\[2mm]
&=\int_{\R^3} \vert \Gamma_{t,x}^{-1}(w) \vert^{\ell} e^{3t}  \mathbf{1}_{\mathcal{O}^t}(x,\Gamma_{t,x}^{-1}(w)) \, \mathbf{1}_{\vert w \vert \leq 1+\mathrm{R}(t)}  \, \mathbf{1}_{ \Lambda_{t,w}(x)_3 >1+\mathrm{L}(t)}  \ f_0(\Lambda_{t,w}(x),w) \vert \det \, \mathrm{D}_w \Gamma_{t,x}^{-1} (w)  \vert \, \mathrm{d}w \\
& \lesssim \dfrac{1}{(1+\mathrm{L}(t))^{k_2} }\int_{\R^3}  (\vert w \vert^{\ell}+1) \mathbf{1}_{\vert w \vert \leq 1+\mathrm{R}(t)}  \,  \mathbf{1}_{\mathcal{O}^t}(x,\Gamma_{t,x}^{-1}(w)) \vert  \Lambda_{t,w}(x)_3 \vert^{k_2}  \, f_0(\Lambda_{t,w}(x),w) \, \mathrm{d}w,
\end{align*}
so that
\begin{align*}
m_{\ell} f^{\flat}(t,x)  \lesssim \dfrac{H_{\ell, k_2}(f_0)}{(1+\mathrm{L}(t))^{k_2}}.
\end{align*}
%\begin{align*}
%m_{\ell} f^{\flat}(t,x)&= \int_{\R^3} \vert v \vert^{\ell} f^{\flat}(t,x,v) \, \mathrm{d}v \\
%&=\int_{\R^3} \vert v \vert^{\ell} e^{3t} \mathbf{1}_{\tau^-(t,x,v)<0} \,  \mathbf{1}_{\vert \V^0_t(x,v) \vert \leq 1+\mathrm{R}(t)} \,  \mathbf{1}_{ \X^0_t(x,v)_3 >1+\mathrm{L}(t)}  \,  f_0(\X^0_t(x,v),\V^0_t(x,v))  \ddv \\[2mm]
%&=\int_{\R^3} \vert \Gamma_{t,x}^{-1}(w) \vert^{\ell} e^{3t}  \mathbf{1}_{\tau^-(t,x,\Gamma_{t,x}^{-1}(w))<0} \, \mathbf{1}_{\vert w \vert \leq 1+\mathrm{R}(t)}  \, \mathbf{1}_{ \Lambda_{t,w}(x)_3 >1+\mathrm{L}(t)}  \ f_0(\Lambda_{t,w}(x),w) \vert \det \ \mathrm{D}_w \Gamma_{t,x}^{-1} (w)  \vert \, \mathrm{d}w \\
%& \lesssim \dfrac{1}{(1+\mathrm{L}(t))^{k_2} }\int_{\R^3}  (\vert w \vert^{\ell}+1) \mathbf{1}_{\vert w \vert \leq 1+\mathrm{R}(t)}  \,  \mathbf{1}_{\tau^-(t,x,\Gamma_{t,x}^{-1}(w))<0} \vert  \Lambda_{t,w}(x)_3 \vert^{k_2}  \, f_0(\Lambda_{t,w}(x),w) \, \mathrm{d}w \\[2mm]
%& \lesssim \dfrac{H_{\ell, k_2}(f_0)}{(1+\mathrm{L}(t))^{k_2}}.
%%&  \lesssim \dfrac{H_{k_2}(f_0)}{(1+\mathrm{L}(t))^{k_2} }\int_{\vert w \vert \leq 1+\mathrm{R}(t)}  (\vert w \vert^{\ell}+1)  \, \mathrm{d}w \\[2mm]
%%&  \lesssim H_{k_2}(f_0) \left[\dfrac{(1+\mathrm{R}(t))^{\ell+3}}{(1+\mathrm{L}(t))^{k_2} }+ \dfrac{(1+\mathrm{R}(t))^{3}}{{(1+\mathrm{L}(t))^{k_2} }} \right].
%\end{align*}
Here, we have again used the fact that if $(x, \Gamma_{t,x}^{-1}(w)) \in \mathcal{O}^t$ then $  \Lambda_{t,w}(x) \in \R^3_+$. This provides the contribution of $f^{\flat}$ to the estimate \eqref{decay:pointmom-GEN}. To get its contribution to the inequality \eqref{decay:Lpmom-GEN}, we perform exactly the same computations as for $f^{\natural}$ to write
\begin{align*}
\Vert m_{\ell} f^{\flat}(t) \Vert_{\Ld^r(\R^3_+)} &= \left[ \int_{\R^3_+} \left(\int_{\R^3} \vert v \vert^{\ell} f^{\flat}(t,x,v) \, \mathrm{d}v \right)^r \, \mathrm{d}x \right]^{1/r} \\
& \lesssim \dfrac{1}{(1+\mathrm{L}(t))^{k_2} } \int_{\vert w \vert \leq 1+\mathrm{R}(t)} \hspace{-2mm} \left[ \vert w \vert^{\ell}+1 \right]  \\
& \qquad \qquad \qquad \qquad \qquad \left(\int_{\R^3} \mathbf{1}_{\Lambda_{t,w}^{-1}(x) \in \R^3_+}  \mathbf{1}_{\mathcal{O}^t}\left(\Lambda_{t,w}^{-1}(x),\Gamma_{t,\Lambda_{t,w}^{-1}(x)}^{-1}(w)\right) \hspace{-1mm} \vert x_3 \vert^{k_2 r}  \, f_0(x,w)^r  \, \mathrm{d}x \right)^{1/r} \hspace{-5mm} \mathrm{d}w \\
&  \lesssim \dfrac{1}{(1+\mathrm{L}(t))^{k_2} } \int_{\vert w \vert \leq 1+\mathrm{R}(t)}  \left[ \vert w \vert^{\ell}+1 \right]   \Vert x_3^{k_2} f_0(\cdot,w) \Vert_{\Ld^{r}(\R^3_+)} \, \mathrm{d}w,
\end{align*}
so that
\begin{align*}
\Vert m_{\ell} f^{\flat}(t) \Vert_{\Ld^r(\R^3_+)} \lesssim \dfrac{F_{\ell,k_2,r}(f_0)}{{(1+\mathrm{L}(t))^{k_2} }}.
\end{align*}
Gathering all the pieces together, we have proven inequalities \eqref{decay:pointmom-GEN} and \eqref{decay:Lpmom-GEN}. This achieves the proof.
\end{proof}

Of course, the previous lemma is only available if we could ensure that a certain EGC is satisfied by the velocity field $u$ at time $t$. 
%In view of Lemma \ref{EGCt:reverse} in Section \ref{Section:EGCabs}, we indeed have the following result.
Thanks to Lemma \ref{EGCt:reverse} stated in Section \ref{Section:EGCabs}, we are able to obtain such a condition.
\begin{propo}\label{PropEGC_after}
For any $t \in (T_0,t^{\star})$, the velocity field $u$ satisfies an EGC in time $t$ with respect to $(\R^2 \times (0,1+\Ld(t)) \times \B_v(0,1+\mathrm{R}(t))$ for some nondecreasing continuous functions $\Ld: [T_0,+\infty) \rightarrow \R^+$ and  $\mathrm{R}: [T_0, + \infty) \rightarrow \R^+$ satisfying for all $s>T_0$
%\begin{align*}
%\dfrac{1}{1+\Ld(s)} \lesssim \dfrac{1}{1+s}, \\ 
%\dfrac{1}{1+s} \lesssim \dfrac{1}{1+\mathrm{R}(s)} \lesssim \dfrac{1}{1+s},
%\end{align*}
\begin{align*}
\dfrac{1}{1+\Ld(s)} \lesssim \dfrac{1}{1+s}, \ \  \dfrac{1}{1+\mathrm{R}(s)} \lesssim \dfrac{1}{1+s},
\end{align*}
where $\lesssim$ only depends on $T_0$ and $g$.

In particular, if $r \in [1, \infty)$, $k_1,k_2,\ell \in \R^+$ and $q>3$ satisfy
\begin{align*}
q>k_1+\ell +3,
\end{align*} then for almost every $x \in \R^3_+$
%\begin{align}
%\label{decay:pointmom}m_{\ell} f(t,x) &  \lesssim  \dfrac{N_q(f_0)}{(1+t)^{k_1}}+ H_{k_2}(f_0)  \left[\dfrac{1}{(1+t)^{k_2-\ell-3} }+ \dfrac{1}{{(1+t)^{k_2-3} }} \right], \\[3mm]
%\label{decay:Lpmom}\Vert m_{\ell} f(t) \Vert_{\Ld^r(\R^3_+)} & \lesssim \dfrac{K_{q,r}(f_0)}{(1+t)^{k_1}}+ F_{k_2,r}(f_0) \left[\dfrac{1}{(1+t)^{k_2-\ell-3} }+ \dfrac{1}{{(1+t)^{k_2-3} }} \right], 
%\end{align}
\begin{align}
\label{decay:pointmom}m_{\ell} f(t,x) &  \lesssim  \dfrac{N_q(f_0)}{(1+t)^{k_1}}+  \dfrac{H_{\ell,k_2}(f_0)}{(1+t)^{k_2}} , \\[3mm]
\label{decay:Lpmom}\Vert m_{\ell} f(t) \Vert_{\Ld^r(\R^3_+)} & \lesssim \dfrac{K_{q,r}(f_0)}{(1+t)^{k_1}}+ \dfrac{F_{\ell,k_2,r}(f_0)}{(1+t)^{k_2}},
\end{align}
where $\lesssim$ only depends on $k_1, k_2, q,\ell,g,\delta_0$ and $T_0$.
\end{propo}
\begin{proof}
Thanks to Lemma \ref{EGCt:reverse}, we know that if $t-1/2>t_0$
%(for some $\alpha>0$ to be determined later)
, the trivial velocity field satisfies an EGC in time $t-1/2$ with respect to $(\R^2 \times (0,1+\widetilde{\Ld}(t-1/2)) \times \B_v(0,1+\widetilde{\mathrm{R}}(t-1/2))$  for some continuous nondecreasing functions $\widetilde{\Ld}: [t_0,+\infty) \rightarrow \R^+$ and  $\widetilde{\mathrm{R}}: [t_0, + \infty) \rightarrow \R^+$ such that for any $s>t_0$
%\begin{align}\label{Maj:LR}
%\dfrac{1}{1+\widetilde{\Ld}(s)} \lesssim \dfrac{1}{1+s}, \ \ \dfrac{1}{1+s} \lesssim \dfrac{1}{1+\widetilde{\mathrm{R}}(s)} \lesssim \dfrac{1}{1+s},
%\end{align}
\begin{align}\label{Maj:LR}
\dfrac{1}{1+\widetilde{\Ld}(s)} \lesssim \dfrac{1}{1+s}, \ \  \dfrac{1}{1+\widetilde{\mathrm{R}}(s)} \lesssim \dfrac{1}{1+s},
\end{align}
where $\lesssim$ depends on $t_0$ and $g$. Furthermore, by definition of $\delta_0$, we have
\begin{align*}
\int_0^t \Vert u(s) \Vert_{\Ld^{\infty}(\R^3_+)} \, \mathrm{d}s <\delta_0 <\kappa_{1/2},
\end{align*}
so that, owing to Lemma \ref{LM:perturbEGC}, we get the fact that the velocity field $u$ satisfies an EGC in time $t$ with respect to $(\R^2 \times (0,1+\widetilde{\Ld}(t-1/2)) \times \B_v(0,1+\widetilde{\mathrm{R}}(t-1/2))$. We have $t>T_0=t_0+1>t_0+1/2$ so that if we set $\Ld(s):=\widetilde{\Ld}(s-1/2)$ and $\mathrm{R}(s):=\widetilde{\mathrm{R}}(s-1/2)$, we observe that the velocity field $u$ satisfies an EGC in time $t$ with respect to  $(\R^2 \times (0,1+\Ld(t)) \times \B_v(0,1+\mathrm{R}(t))$. Furthermore, by \eqref{Maj:LR}, for all $s>t_0+1/2$
\begin{align*}
\dfrac{1}{1+\Ld(s)} \lesssim \dfrac{1}{1+s}, \ \ \dfrac{1}{1+\mathrm{R}(s)} \lesssim \dfrac{1}{1+s}. 
%\dfrac{1}{1+s} \lesssim \dfrac{1}{1-\alpha+s} \lesssim \dfrac{1}{1+\mathrm{R}(s)} &\lesssim \dfrac{1}{1-\alpha+s} \lesssim \dfrac{1}{1+s}.
\end{align*}
The last part of the statement readily comes from the previous estimates and the inequalities \eqref{decay:pointmom-GEN}-\eqref{decay:Lpmom-GEN} in Lemma \ref{LM:decaymom-abs}.
\end{proof}

\bigskip

We now fix $T \in (T_0,t^{\star})$. Since the decay provided by the absorption can only be considered after time $T_0$, we will split the study of the decay estimates between $[0,T_0]$ and $[T_0,T]$. In particular, the local in time estimates of Section \ref{Section:Prelim} will help us to treat the part on $[0,T_0]$.

In view of the assumptions \eqref{data:hyp} and \eqref{smallness:condition:THM} on the initial data, we will not be very precise about the exact range of exponents which are used when the quantities $N_q(f_0), H_{q,m}(f_0), K_{q,r}(f_0), F_{q,m,r}(f_0)$ of Definition \ref{decay_f0} may appear. Note that they are nondecreasing when the parameters $q$ and $m$ increase. In the last part of the current section, we will consider the maximum of the finite family of indices we have used previously. Moreover, these indices may sometimes depend on the exponent $\gamma$ that we use to quantify the decay in time but this one will be somehow fixed in the last step of the bootstrap.
Thus, all the following results must be understood as if we take the involved exponents large enough. Let us also recall the symbol $\lesssim_0$ introduced in Notation \ref{nota:lesssim_0}.

\begin{coro}\label{OKdecay:BrinkPoint}
There exists $\varepsilon>0$ small enough such that the following holds. Let $p \in (3,3+\varepsilon)$. There exist nondecreasing positive functions $\varphi_2$ and $\varphi_p$ vanishing at $0$ such that for any $k>0$ satisfying $q>k+3$, we have
\begin{align}\label{decayBrinkL2point}
\forall t \in [0,T], \ \ \Vert j_f(t)-\rho_f u(t) \Vert_{\Ld^2(\R^3_+)} \lesssim_{0} \dfrac{\varphi_2 \big(N_q(f_0)+H_{0,k}(f_0) \big)}{(1+t)^{\frac{k}{2}}},
\end{align}
and if $\Vert u \Vert_{\Ld^{\infty}(0,T;\Ld^6(\R^3_+))} \lesssim_0 1$, we have
\begin{align}\label{decayBrinkLppoint}
\forall t \in [0,T], \ \ \Vert j_f(t)-\rho_f u(t) \Vert_{\Ld^p(\R^3_+)} \lesssim_{0} \dfrac{\varphi_p \big(N_q(f_0)+H_{0,k}(f_0) \big)}{(1+t)^{k\frac{p-1}{p}}}.
\end{align}
\end{coro}
\begin{proof}
We separate the study between $[0,T_0]$ and $[T_0,T]$.

$\bullet$ If $t \in [T_0,T]$, we use Lemma \ref{lem-Dp} and get for all $r \geq 2$
\begin{align*}
\|( j_f - \rho_f u )(t)\|_{\mathrm{L}^r(\R^3_+)} &\leq   \| \rho_f(t) \|_{\Ld^\infty(\R^3_+))}^{\frac{r-1}{r}} \left(\int_{\R^3_+ \times \R^3} f(t,x,v) |v-u(t,x)|^r \, \mathrm{d} x \, \mathrm{d} v\right)^{1/r} \\
& \lesssim  \| \rho_f (t) \|_{\Ld^\infty(\R^3_+))}^{\frac{r-1}{r}} \left(\int_{\R^3_+ \times \R^3} f(t,x,v) \vert v \vert^r \, \mathrm{d} x \, \mathrm{d} v + \int_{\R^3_+ } \rho_f(t,x) \vert u(t,x)\vert^r \, \mathrm{d} x \right)^{1/r}.
\end{align*}
We first focus on the estimate \eqref{decayBrinkL2point}. We have
\begin{align*}
&\|( j_f - \rho_f u )(t)\|_{\mathrm{L}^2(\R^3_+)} \\
&\lesssim  \Vert m_0 f(t) \Vert_{\Ld^{\infty}(\R^3_+)}^{\frac{1}{2}} \left(\Vert m_2 f(t) \Vert_{\Ld^1(\R^3_+)} + \Vert m_0 f(t) \Vert_{\Ld^{\infty}(\R^3_+)} \E(t) \right)^{1/2}\\
& \lesssim  \Vert m_0 f(t) \Vert_{\Ld^{\infty}(\R^3_+)}^{\frac{1}{2}} \Bigg(\Vert m_2 f(t) \Vert_{\Ld^1(\R^3_+)}  \\ 
& \qquad \qquad \qquad \qquad \quad+ \Vert m_0 f(t) \Vert_{\Ld^{\infty}(\R^3_+)} \left[ \E(0) + \E(0)^{1/2}T_0+T_0^2 + \int_{T_0}^{T} \Vert j_f(\tau) \Vert_{\Ld^1(\R^3_+)} \, \mathrm{d}\tau   \right] \Bigg)^{1/2},
\end{align*}
where we have used the energy inequality \eqref{ineq-energy} and Lemma \ref{j_f:L^1L^1}. The estimates \eqref{decay:pointmom} and \eqref{decay:Lpmom} of Proposition \ref{PropEGC_after} then imply
\begin{align*}
\|( j_f - \rho_f u )(t)\|_{\mathrm{L}^2(\R^3_+)} &\leq  \left( \dfrac{N_{q}(f_0)}{(1+t)^{k}}+\dfrac{H_{0,k}(f_0)}{(1+t)^{k}} \right)^{1/2} \\
& \quad \times \Bigg(\dfrac{K_{q_1,1}(f_0)}{(1+t)^{k_1}}+\dfrac{F_{2,k_1,1}(f_0)}{(1+t)^{k_1}}+ \left(\dfrac{N_{q_2}(f_0)}{(1+t)^{k_2}}+\dfrac{H_{0,k_2}(f_0)}{(1+t)^{k_2}} \right) \\
& \qquad \qquad \qquad \qquad \times \left[ \E(0) + \E(0)^{1/2}T_0+T_0^2 + \int_{T_0}^{T} \left( \dfrac{K_{q_3,1}(f_0)}{(1+\tau)^{k_3}} + \dfrac{F_{1,k_3,1}(f_0)}{(1+\tau)^{k_3}} \right) \, \mathrm{d}\tau   \right] \Bigg)^{1/2},
\end{align*}
provided that $q>k+3$, $q_1>k_1+5$, $q_2>k_2+3$ and $q_3>k_3+4$. Thus, we obtain for all $t \in [T_0,T]$
\begin{align*}
\|( j_f - \rho_f u )(t)\|_{\mathrm{L}^2(\R^3_+)}  \lesssim_0  \dfrac{N_{q}(f_0)^{1/2}}{(1+t)^{k/2}}+\dfrac{H_{0,k}(f_0)^{1/2}}{(1+t)^{k/2}}.
\end{align*}
if we can choose $k_1, k_2>0$ and $k_3>1$, which is possible in view of the assumptions \eqref{data:hyp} of Theorem \ref{thm1}.

We now treat the estimate \eqref{decayBrinkLppoint}, assuming $\Vert u \Vert_{\Ld^{\infty}(0,T;\Ld^6(\R^3_+))} \lesssim_0 1$. By Hölder's inequality, we now have for $p \in (3,3+\eps)$
\begin{align*}
\| (j_f - \rho_f u )(t)\|_{\mathrm{L}^p(\R^3_+)}  & \lesssim  \| \rho_f (t) \|_{\Ld^\infty(\R^3_+))}^{\frac{p-1}{p}} \left(\int_{\R^3_+ \times \R^3} f(t,x,v) \vert v \vert^p \, \mathrm{d} x \, \mathrm{d} v + \int_{\R^3_+ } \rho_f(t,x) \vert u(t,x)\vert^p \, \mathrm{d} x \, \mathrm{d} v\right)^{1/p} \\
& \leq  \| m_0 f (t) \|_{\Ld^\infty(\R^3_+))}^{\frac{p-1}{p}} \left(\Vert m_p f(t) \Vert_{\Ld^1(\R^3_+)} + \Vert u \Vert_{\Ld^{\infty}(0,T;\Ld^6(\R^3_+))}^{p} \| m_0 f (t) \|_{\Ld^{r_{\varepsilon}(p)}(\R^3_+)} \right)^{1/p},
\end{align*}
where $2 <r_{\varepsilon}(p) <6/(3-\varepsilon)$. Thanks to  $\Vert u \Vert_{\Ld^{\infty}(0,T;\Ld^6(\R^3_+))} \lesssim_0 1$, the estimates \eqref{decay:pointmom}-\eqref{decay:Lpmom} of Proposition \ref{PropEGC_after} and \eqref{data:hyp}, we obtain in the same way as before that
\begin{align*}
\|( j_f - \rho_f u )(t)\|_{\mathrm{L}^p(\R^3_+)} \lesssim_0   \dfrac{N_{q}(f_0)^{\frac{p-1}{p}}}{(1+t)^{k\frac{p-1}{p}}}+\dfrac{H_{0,k}(f_0)^{\frac{p-1}{p}}}{(1+t)^{k\frac{p-1}{p}}},
\end{align*}
for $q>k+3$ and $t \in [T_0,T]$. Note that the the treament of the term $\Vert m_p f(t) \Vert_{\Ld^1(\R^3_+)}$ via \eqref{decay:Lpmom} requires $q \geq7$, which is allowed by \eqref{data:hyp}.
%
%-----
%We use the triangular inequality and the Hölder inequality to write
%\begin{align*}
%\Vert j_f(t)-\rho_f u(t) \Vert_{\Ld^2(\R^3_+)} &\leq \Vert j_f(t) \Vert_{\Ld^2(\R^3_+)} +\Vert\rho_f u(t) \Vert_{\Ld^2(\R^3_+)} \\
%&\leq \Vert j_f(t) \Vert_{\Ld^2(\R^3_+)} +\Vert\rho_f(t) \Vert_{\Ld^{\infty}(\R^3_+)}\Vert u(t) \Vert_{\Ld^2(\R^3_+)},
%\end{align*}
%We use the energy inequality \eqref{ineq-energy} and get 
%\begin{align*}
%\Vert j_f(t)-\rho_f u(t) \Vert_{\Ld^2(\R^3_+)} \leq \Vert j_f(t) \Vert_{\Ld^2(\R^3_+)} +\Vert\rho_f(t) \Vert_{\Ld^{\infty}(\R^3_+)}\left[ \E(0) + g \int_{0}^{T} \Vert j_f(\tau) \Vert_{\Ld^1(\R^3_+)} \, \mathrm{d}\tau   \right]^{1/2},
%\end{align*}
%We separate the study of the last r.h.s between the interval $[0,T_0]$ and $[T_0,T]$. If $t \in [T_0,T]$, the estimates \eqref{decay:pointmom} and \eqref{decay:Lpmom} of Proposition \ref{PropEGC_after} then imply
%\begin{align*}
%\Vert j_f(t)-\rho_f u(t) \Vert_{\Ld^2(\R^3_+)} &\lesssim \dfrac{K_{q,2}(f_0)}{(1+t)^{k}}+\dfrac{F_{1,k,2}(f_0)}{(1+t)^{k}} \\
%& \quad + \left( \dfrac{N_q(f_0)}{(1+t)^{k}}+  \dfrac{H_{0,k}(f_0)}{(1+t)^{k}}\right) \Bigg[ \E(0) +\E(0)^{1/2}T_0+T_0^2 \\ 
%& \qquad \qquad \qquad \qquad \qquad \qquad \qquad  + g \int_{T_0}^T \left(\dfrac{K_{q,1}(f_0)}{(1+\tau)^{k_2}} + \dfrac{F_{1,k_2,1}(f_0)}{(1+\tau)^{k_2}} \right)\, \mathrm{d}\tau   \Bigg]^{1/2},
%\end{align*}
%provided that $q>k+1+3$ and $q>k_2+1+3$. Taking $k_2$ (and then $q$) large enough, we obtain the desired estimate on $[T_0,T]$.
%
%

$\bullet$ If $t \in [0,T_0]$, and for the proof of \eqref{decayBrinkL2point}, we use the triangular inequality and the Hölder's inequality to write
\begin{align*}
\Vert j_f(t)-\rho_f u(t) \Vert_{\Ld^2(\R^3_+)} &\leq \Vert j_f(t) \Vert_{\Ld^2(\R^3_+)} +\Vert\rho_f u(t) \Vert_{\Ld^2(\R^3_+)} \\
&\leq \Vert j_f(t) \Vert_{\Ld^2(\R^3_+)} +\Vert\rho_f(t) \Vert_{\Ld^{\infty}(\R^3_+)}\Vert u(t) \Vert_{\Ld^2(\R^3_+)}.
\end{align*}
Hence, by interpolation, Lemma \ref{j_f:L^1L^1} and Lemma \ref{rough:bounds}, we get
\begin{align*}
\Vert j_f(t)-\rho_f u(t) \Vert_{\Ld^2(\R^3_+)} &\leq  \Vert j_f(t) \Vert_{\Ld^{3/2}(\R^3_+)}^{3/4} \Vert j_f(t) \Vert_{\Ld^{\infty}(\R^3_+)}^{1/4} + \Vert\rho_f(t) \Vert_{\Ld^{\infty}(\R^3_+)}[\E(0) +\E(0)^{1/2}T_0+T_0^2] \\
&\leq   \phi(T_0) \Vert j_f \Vert_{\Ld^{\infty}(0,T_0;\Ld^{\infty}(\R^3_+))}^{1/4} + \Vert\rho_f \Vert_{\Ld^{\infty}(0,T_0;\Ld^{\infty}(\R^3_+))}[\E(0) +\E(0)^{1/2}T_0+T_0^2].
\end{align*}
The local estimates \eqref{borneLinftyLinfty-rhoetj} then provides
\begin{align*}
\Vert j_f(t)-\rho_f u(t) \Vert_{\Ld^2(\R^3_+)} & \lesssim  \phi(T_0) (N_q(f_0) \eta(T_0))^{1/4} + (N_q(f_0) \eta(T_0))[\E(0) +\E(0)^{1/2}T_0+T_0^2] \\
& \lesssim_0 N_q(f_0)^{1/4}+N_q(f_0) \\
& \lesssim_0 \dfrac{\varphi(N_q(f_0))}{(1+T_0)^{k/2}},
\end{align*}
for some nonincreasing function $\varphi$ vanishing at $0$, provided that $N_q(f_0)$ is taken small enough by the smallness assumption \eqref{smallness:condition:THM}. Concerning the estimate \eqref{decayBrinkLppoint} under the assumption $\Vert u \Vert_{\Ld^{\infty}(0,T;\Ld^6(\R^3_+))} \lesssim_0 1$, we proceed as before by writing
\begin{align*}
\Vert  j_f(t)-\rho_f u(t) \Vert_{\Ld^p(\R^3_+)} &\leq \Vert  j_f(t) \Vert_{\Ld^p(\R^3_+)}+ \Vert  \rho_f u(t) \Vert_{\Ld^p(\R^3_+)} \\
&\leq \Vert  j_f(t) \Vert_{\Ld^p(\R^3_+)}+ \Vert  \rho_f(t) \Vert_{\Ld^{\widetilde{r_{\varepsilon}}}(\R^3_+)} \Vert  u(t) \Vert_{\Ld^{6}(\R^3_+)},
\end{align*}
where $
6 < \widetilde{r_{\varepsilon}} < \dfrac{18+6\varepsilon}{3-\varepsilon}$. By interpolation, this turns into
\begin{align*}
\Vert  (j_f-\rho_f u) \Vert_{\Ld^p(0,T_0;\Ld^p(\R^3_+))}^p & \leq  \Vert j_f \Vert_{\Ld^{\infty}(0,T_0;\Ld^{3/2}(\R^3_+))}^{3/2} \Vert  j_f \Vert_{\Ld^{\infty}(0,T_0;\Ld^{\infty}(\R^3_+))}^{p-3/2} \\
& \quad   +    \Vert  \rho_f \Vert_{\Ld^{\infty}(0,T_0;\Ld^1(\R^3_+))}^{p/\widetilde{r_{\varepsilon}}}  \Vert \rho_f \Vert_{\Ld^{\infty}(0,T_0;\Ld^{\infty}(\R^3_+))}^{p-p/\widetilde{r_{\varepsilon}}} \Vert  u \Vert_{\Ld^{\infty}(0,T_0;\Ld^6(\R^3_+))}^p  \\
& \lesssim_0   N_q(f_0)^{p-3/2}+N_q(f_0)^{p-p/\widetilde{r_{\varepsilon}}},
\end{align*}
thanks to the Lemma \ref{rough:bounds}, Lemma \ref{loose:mass}, the bound \eqref{borneLinftyLinfty-rhoetj} and $\Vert u \Vert_{\Ld^{\infty}(0,T;\Ld^6(\R^3_+))} \lesssim_0 1$.
By taking $N_q(f_0)$ small enough, we get for all $t \in [0,T_0]$
\begin{align*}
\Vert j_f(t)-\rho_f u(t) \Vert_{\Ld^p(\R^3_+)}  \lesssim_0 \dfrac{\varphi(N_q(f_0))}{(1+T_0)^{k\frac{p-1}{p}}}.
\end{align*}
Combining the estimates on $[0,T_0]$ and on $[T_0,T]$, we finally obtain the result.
\end{proof}

\begin{coro}
For all $0 \leq t \leq T$ and for any $k>1$ such that $q>k+3$, we have
\begin{align}
\Vert \nabla u(t) \Vert_{\Ld^2(\R^3_+)}^2 + \dfrac{1}{2}\int_0^t \Vert \mathrm{D}^2 u(s) \Vert_{\Ld^2(\R^3_+)}^2 \, \mathrm{d}s   \lesssim_0 \Vert  u_0 \Vert_{\H^1(\R^3_+)}^2 + \varphi_2 \big( N_q(f_0)+H_{0,k}(f_0)\big).
\end{align}
Furthermore, the following estimate holds for all $0 \leq t \leq T$
\begin{align}
\Vert u \Vert_{\Ld^{\infty}(0,t; \Ld^6(\R^3_+))}^2 \lesssim \Vert \nabla u \Vert_{\Ld^{\infty}(0,t; \Ld^2(\R^3_+))}^2 \lesssim_0 \Vert  u_0 \Vert_{\H^1(\R^3_+)}^2 + \varphi_2 \big( N_q(f_0)+H_{0,k}(f_0) \big), \label{controlparab1}
\end{align}
and in particular, the estimate \eqref{decayBrinkLppoint} holds true.
\end{coro}
\begin{proof}
As we have $T \in (T_0,t^{\star})$ and in view of \eqref{ineq:BrinkL2L2:loc}, we only treat the case where $t \in [T_0,T]$. Since the time $t$ is a strong existence time (see Definition \ref{strongtime}), Proposition \ref{propdatasmall:VNSreg} entails the following estimate 
\begin{align*}
\Vert \nabla u(t) \Vert_{\Ld^2(\R^3_+)}^2 + \dfrac{1}{2}\int_0^t \Vert \mathrm{D}^2 u(s) \Vert_{\Ld^2(\R^3_+)}^2 &\lesssim  \Vert \nabla u_0 \Vert_{\Ld^2(\R^3_+)}^2 + \Vert j_f-\rho_f u \Vert_{\Ld^2(0,t;\Ld^2(\R^3_+))}^2.
\end{align*}
Applying Corollary \ref{OKdecay:BrinkPoint} with $k>1$ leads to
\begin{align*}
\Vert \nabla u(t) \Vert_{\Ld^2(\R^3_+)}^2 + \dfrac{1}{2}\int_0^t \Vert \mathrm{D}^2 u(s) \Vert_{\Ld^2(\R^3_+)}^2 &\lesssim_0 \Vert \nabla u_0 \Vert_{\Ld^2(\R^3_+)}^2 +  \varphi_2 \big( N_q(f_0)+\H_{0,k}(f_0) \big),
\end{align*}
which is the first part of the result. The other statement is then a consequence of the Sobolev embedding.
\end{proof}

\subsection{Estimates with a polynomial weight in time}
We recall that we have fixed $T \in (T_0,t^{\star})$. As in \cite{HK}, the guiding line in order to obtain polynomial weighted estimates for the fluid velocity $u$ is to derive a Stokes system satisfied by a weighted version of $u$, in such a way that the maximal regularity result of Section \ref{AnnexeMaxregStokes} can be applied. Note that since we do not have enough information about the regularity of $u(T_0)$ with respect to the interpolation spaces defined in \eqref{domaineDqs}, we rely on this result on the whole interval $[0,T]$.

In what follows, the function $\Psi$ is a generic continuous positive and increasing function that may change from one line to another. Let us also recall the symbol $\lesssim_0$ introduced in Notation \ref{nota:lesssim_0}.
%which may also depend on quantities involving the initial data but increasing when these quantities increase: for simplicity, this dependency will not be explicitly written down.

\medskip

We first state the following $\Ld^2 \H^2$ regularity result which is reminiscent of \cite{HK}.
\begin{lem}
For all $\gamma \in (0,3/4)$, we have
%\begin{align}\label{maxregL2L2}
%\Vert (1+t)^{\gamma} \partial_t u \Vert_{\Ld^2(0,T;\Ld^2(\R^3_+))} + \Vert (1+t)^{\gamma} \D^2u \Vert_{\Ld^2(0,T;\Ld^2(\R^3_+))} \lesssim_{T_0, \gamma} \mathrm{A}_0+ \Psi \left( \Vert u_0 \Vert_{\Ld^2(\R^3_+)}^2+\Vert u_0 \Vert_{\Ld^1(\R^3_+)}^2 \right) .
%\end{align}
\begin{align}\label{maxregL2L2}
\Vert (1+t)^{\gamma} \D^2u \Vert_{\Ld^2(0,T;\Ld^2(\R^3_+))} \lesssim_{0} \Psi \left( \Vert u_0 \Vert_{\Ld^1 \cap \Ld^2(\R^3_+)} +N_q(f_0)+H_{0,k}(f_0) \right)+  \Vert  u_0 \Vert_{\H^1(\R^3_+)}.
\end{align}
\end{lem}
\begin{proof}
Let $\gamma \in (0,3/4)$ and set $U:=(1+t)^{\gamma} u$. We observe that $U$ satisfies the following Stokes system on $[0,T]$
\begin{equation}\label{Stokes:Uweighted}
\left\{
      \begin{aligned}
         \partial_t U +A_2 U &=S(u,f),\\
         U_{\mid x_3=0}&=0,\\
        U_{\mid t=0}&=u_0, 
    \end{aligned}
    \right.
\end{equation}
where $A_2$ stands for the Stokes operator on $\Ld^2_{\mathrm{div}}(\R^3_+)$ and where
\begin{align}\label{def:S(u,f)}
S(u,f):=(1+t)^{\gamma} \P(j_f-\rho_f u) -(1+t)^{\gamma}\mathbb{P} \left( u \cdot \nabla \right) u+\gamma (1+t)^{\gamma-1} u.
\end{align}
According to the maximal regularity result in $\Ld^2\Ld^2$ for the Stokes system (see Section \ref{AnnexeMaxregStokes} in the Appendix), we get
\begin{align*}
\Vert  \partial_t U \Vert_{\Ld^2(0,T;\Ld^2(\R^3_+))} + \Vert  \D^2 \, U \Vert_{\Ld^2(0,T;\Ld^2(\R^3_+))} \lesssim \Vert  S(u,f) \Vert_{\Ld^2(0,T;\Ld^2(\R^3_+))} + \Vert  u_0 \Vert_{\H^1(\R^3_+)},
\end{align*}
where $\lesssim$ is independent of $T$, $u$ and $f$.

We now estimate each term in the modified source term $S(u,f)$. Owing to Corollary \ref{OKdecay:BrinkPoint} with $k=7/2$ and Theorem \ref{cond:decay}, we have
\begin{align}\label{eq:ufirst}
(1+t)^{\gamma-1} \Vert u(t) \Vert_{\Ld^2(\R^3_+)} \leq \dfrac{\Psi \left( \Vert u_0 \Vert_{\Ld^1 \cap \Ld^2(\R^3_+)} +N_q(f_0)+H_{0,k}(f_0) \right)}{(1+t)^{3/4-\gamma+1}},
\end{align}
which is bounded independently of $T$ in $\Ld^2(0,T)$ because $\gamma \in (0,3/4)$.
For the term coming from the Brinkman force $j_f-\rho_f u$, we use the continuity of the Leray projection on $\Ld^2(\R^3_+)$ and the estimate \eqref{decayBrinkL2point} of Corollary \ref{OKdecay:BrinkPoint} which entail
%\begin{align}\label{eq:Brink1}
%\Vert (1+t)^{\gamma}(j_f-\rho_f u) \Vert_{\Ld^2(0,T;\Ld^2(\R^3))} \lesssim_{T_0,\gamma} \mathrm{A}_0,
%\end{align}
\begin{align}\label{eq:Brink1}
\begin{split}
\Vert (1+t)^{\gamma}\P(j_f-\rho_f u) \Vert_{\Ld^2(0,T;\Ld^2(\R^3_+))}^2 &\lesssim_{0} \int_0^T (1+s)^{2 \gamma} \Vert j_f(s)-\rho_f u(s) \Vert_{\Ld^2(\R^3_+)}^2  \, \mathrm{d}s \\
& \lesssim_{0} \varphi_2\big(N_q(f_0)+H_{0,k}(f_0) \big) \int_0^T \dfrac{1}{(1+s)^{k-2\gamma}} \, \mathrm{d}s \\
& \lesssim_0 \Psi \left( \Vert u_0 \Vert_{\Ld^1 \cap \Ld^2(\R^3_+)} +N_q(f_0)+H_{0,k}(f_0) \right),
\end{split}
\end{align}
for $q>k+3>4+3/2$. For the convection term $(u \cdot \nabla )u$, we proceed as in the proof of Theorem \ref{RegParabNS} in the Appendix (more precisely, see \eqref{majConv1}-\eqref{majConv2}) to obtain, with the Young inequality, that
\begin{align*}
\Vert (u \cdot \nabla ) u(t) \Vert_{\Ld^2(\R^3_+)}^2 \leq C{\eta^{-4}} \Vert u(t) \Vert_{\Ld^2(\R^3_+)}^2  \Vert \nabla u(t) \Vert_{\Ld^2(\R^3_+)}^4 \Vert \nabla u(t) \Vert_{\Ld^2(\R^3_+)}^4 + \frac{3\eta^{4/3}}{4} \Vert \D^2 u(t) \Vert_{\Ld^2(\R^3_+)}^2,
\end{align*}
where $\eta>0$ can be taken as small as we want. We then use the Gagliardo-Nirenberg-Sobolev inequality to write
\begin{align*}
\Vert \nabla u(t) \Vert_{\Ld^2(\R^3_+)} \lesssim \Vert u(t) \Vert_{\Ld^2(\R^3_+)}^{1/2} \Vert \D^2 u(t) \Vert_{\Ld^2(\R^3_+)}^{1/2},
\end{align*}
therefore this yields 
\begin{align*}
\Vert (1+t)^{\gamma} \P (u \cdot \nabla ) u \Vert_{\Ld^2(0,T;\Ld^2(\R^3_+))}^2 \leq \int_{0}^T (1+t)^{2\gamma} \left[ C{\eta^{-4}} \Vert u(t) \Vert_{\Ld^2(\R^3_+)}^4  \Vert \nabla u(t) \Vert_{\Ld^2(\R^3_+)}^4 + \frac{3 \eta^{4/3}}{4}  \right] \Vert \D^2 u(t) \Vert_{\Ld^2(\R^3_+)}^2 \, \mathrm{d}t.
\end{align*}
Furthermore, owing to Theorem \ref{cond:decay} and the estimate \eqref{controlparab1}, we can write
\begin{align*}
\Vert u(t) \Vert_{\Ld^2(\R^3_+)}^4 \Vert  \nabla u(t) \Vert_{\Ld^2(\R^3_+)}^4 &\leq \Psi \left( \Vert u_0 \Vert_{\Ld^1 \cap \Ld^2(\R^3_+)} +N_q(f_0)+H_{0,k}(f_0) \right)^4 \Vert \nabla u \Vert_{\Ld^{\infty}(0,t; \Ld^2(\R^3_+))}^4 \\
%&\lesssim  \Psi(N_q(f_0), K_{m,1}(f_0), \E(0), \Vert  u_0 \Vert_{\Ld^1(\R^3_+)}, \Vert \nabla u_0 \Vert_{\Ld^2(\R^3_+)} )
&\lesssim_0 \Psi \left( \Vert u_0 \Vert_{\Ld^1 \cap \Ld^2(\R^3_+)} +N_q(f_0)+H_{0,k}(f_0) \right)^4 .
\end{align*}
Using the smallness assumption \eqref{smallness:condition:THM}, we can first choose $\eta$ and then $\Vert u_0 \Vert_{\Ld^1 \cap \Ld^2(\R^3_+)}^2 +N_q(f_0)+H_{0,k}(f_0) $ small enough so that
\begin{align}\label{toabsorb:convec}
\Vert (1+t)^{\gamma} \P ( u \cdot \nabla ) u \Vert_{\Ld^2(0,T;\Ld^2(\R^3_+))} \leq \frac{1}{2} \Vert (1+t)^{\gamma} \D^2u \Vert_{\Ld^2(0,T;\Ld^2(\R^3_+))}.
\end{align}
We eventually combine \eqref{eq:ufirst}, \eqref{eq:Brink1} together with \eqref{toabsorb:convec} so that we get
\begin{align*}
\Vert (1+t)^{\gamma} \D^2u \Vert_{\Ld^2(0,T;\Ld^2(\R^3_+))} & \lesssim_0 \Psi \left( \Vert u_0 \Vert_{\Ld^1 \cap \Ld^2(\R^3_+)} +N_q(f_0)+H_{0,k}(f_0) \right)+  \Vert  u_0 \Vert_{\H^1(\R^3_+)},
%& \lesssim \mathrm{A}_0 + \Psi \left( \Vert u_0 \Vert_{\Ld^2(\R^3_+)}^2+\Vert u_0 \Vert_{\Ld^1(\R^3_+)}^2 \right),
\end{align*}
which concludes the proof.
\end{proof}
By interpolation, the previous weighted $\Ld^2$ maximal parabolic estimates allow us to obtain improved $\Ld^p$ estimates of the weighted source term $S(u,f)$ defined in \eqref{def:S(u,f)}, for $p>3$ (close enough to $3$). We refer to \cite{HK} where the proof of the two following results can be found in the whole space case and apply \textit{mutatis mutandis} to the half-space case thanks to Theorem \ref{cond:decay}, Corollary \ref{OKdecay:BrinkPoint} and the estimate \eqref{maxregL2L2}.
%We refer to \cite[Lemma 5.5]{HK} and \cite[Lemma 5.6]{HK}, the proof of which applies \textit{mutadis mutandis} to the half-space case thanks to Theorem \ref{cond:decay} and Proposition \ref{OKdecay:BrinkPoint}. 
\begin{coro}\label{coro:source1}
Let $p>3$. For all $\gamma \in (0, \frac{17}{8}-\frac{7}{4p})$, we have
\begin{align*}
\Vert (1+t)^{\gamma-1}u \Vert_{\Ld^p(0,T;\Ld^p(\R^3_+))} \lesssim_0 \Psi \left( \Vert u_0 \Vert_{\Ld^1 \cap \Ld^2(\R^3_+)} +N_q(f_0)+H_{0,k}(f_0) \right)+  \Vert  u_0 \Vert_{\H^1(\R^3_+)}.
\end{align*}
\end{coro}

\begin{coro}\label{coro:source2}
There exists $\varepsilon>0$ such that for all $p \in (3,3+\varepsilon)$, the following holds. For all $\gamma > 0$, we have
\begin{align*}
\Vert (1+t)^{\gamma} (u \cdot \nabla)u \Vert_{\Ld^p(0,T;\Ld^p(\R^3_+))} \lesssim_0 \Psi \left( \Vert u_0 \Vert_{\Ld^1 \cap \Ld^2(\R^3_+)} +N_q(f_0)+H_{0,k}(f_0) \right) \Vert (1+t)^{\gamma} \D^2 u \Vert_{\Ld^p(0,T;\Ld^p(\R^3_+))}.
\end{align*}
\end{coro}

\bigskip

We are now in position to provide some $\Ld^p \Ld^p$ estimates for the weighted source term 
\begin{align*}
S(u,f)=(1+t)^{\gamma} \P(j_f-\rho_f u) -(1+t)^{\gamma}\mathbb{P} \left( u \cdot \nabla \right) u+\gamma(1+t)^{\gamma-1} u,
\end{align*}
with $p>3$ (close enough to $3$) and $\gamma>0$ large enough, thanks to the next lemma.
\begin{lem}\label{LMsourceLplast}
There exists $\varepsilon>0$ such that for all $p \in (3,3+\varepsilon)$, the following holds. For all $\gamma \in (0, \frac{17}{8}-\frac{7}{4p})$, we have
\begin{align*}
\Vert S(u,f) \Vert_{\Ld^p(0,T;\Ld^p(\R^3_+))} \lesssim_{0}  \left[\Psi \left( \Vert u_0 \Vert_{\Ld^1 \cap \Ld^2(\R^3_+)} +N_q(f_0)+H_{0,k}(f_0) \right)+  \Vert  u_0 \Vert_{\H^1(\R^3_+)} \right] \\
 \times \left[ 1+ \Vert (1+t)^{\gamma} \D^2 u \Vert_{\Ld^p(0,T;\Ld^p(\R^3_+))} \right].
\end{align*}
\end{lem}
\begin{proof}
We first use Corollary \ref{coro:source1} and Corollary \ref{coro:source2} to estimate the contribution of $(1+t)^{\gamma}\mathbb{P} \left( u \cdot \nabla \right) u-\gamma(1+t)^{\gamma-1} u$, that is
\begin{multline*}
\Vert (1+t)^{\gamma}\mathbb{P} \left( u \cdot \nabla \right) u+\gamma(1+t)^{\gamma-1} u \Vert_{\Ld^p(0,T;\Ld^p(\R^3_+))} \\
\lesssim_{0} \left[\Psi \left( \Vert u_0 \Vert_{\Ld^1 \cap \Ld^2(\R^3_+)} +N_q(f_0)+H_{0,k}(f_0) \right)+  \Vert  u_0 \Vert_{\H^1(\R^3_+)} \right] \\
 \times \left[ 1+ \Vert (1+t)^{\gamma} \D^2 u \Vert_{\Ld^p(0,T;\Ld^p(\R^3_+))} \right],
\end{multline*}
for $p \in (3,3+\varepsilon)$, where we have also used the continuity of the Leray projection of $\Ld^p(\R^3_+)$.

In order to treat the term coming from the Brinkman force, we use the estimate \eqref{decayBrinkLppoint} of Corollary \ref{OKdecay:BrinkPoint} (which is valid in view of \eqref{controlparab1}). We get
\begin{align}\label{eq:Brink1}
\begin{split}
\Vert (1+t)^{\gamma}\P(j_f-\rho_f u) \Vert_{\Ld^p(0,T;\Ld^p(\R^3_+))}^p &\lesssim_{0} \int_0^T (1+s)^{p \gamma} \Vert j_f(s)-\rho_f u(s) \Vert_{\Ld^p(\R^3_+)}^p  \, \mathrm{d}s \\
& \lesssim_{0} \varphi_p\big(N_q(f_0)+H_{0,k}(f_0) \big) \int_0^T \dfrac{1}{(1+s)^{k(p-1)-p\gamma}} \, \mathrm{d}s \\
& \lesssim_0 \Psi \left( \Vert u_0 \Vert_{\Ld^1 \cap \Ld^2(\R^3_+)} +N_q(f_0)+H_{0,k}(f_0) \right),
\end{split}
\end{align}
for $k$ large enough with respect to $\gamma \in (0, \frac{17}{8}-\frac{7}{4p})$. More precisely, we have to choose $k$ such that
\begin{align*}
q>k+3>\frac{1+p\gamma}{p-1}+3.
\end{align*}
Since $\frac{1+p\gamma}{p-1} <\frac{17}{16}p-\frac{3}{8}<\frac{17}{16}\frac{3(2+s)}{4}-\frac{3}{8}$ with $s \in (2,3)$ (see \eqref{data:hyp}), this procedure is allowed by the assumption \eqref{data:hyp} of Theorem \ref{thm1} on the exponents $q$ and $m$. 
% We have $\frac{17}{16}\frac{3(2+s)}{4}-\frac{3}{8}}<231/64$
Gathering all the pieces together, we end up with
\begin{align*}
\Vert S(u,f) \Vert_{\Ld^p(0,T;\Ld^p(\R^3_+))}   \lesssim_{0} \left[ \Psi \left( \Vert u_0 \Vert_{\Ld^1 \cap \Ld^2(\R^3_+)} +N_q(f_0)+H_{0,k}(f_0) \right)+  \Vert  u_0 \Vert_{\H^1(\R^3_+)} \right] \\
 \times \left[ 1+ \Vert (1+t)^{\gamma} \D^2 u \Vert_{\Ld^p(0,T;\Ld^p(\R^3_+))} \right],
\end{align*}
and this concludes the proof.
\end{proof}

\begin{coro}\label{coroMaxpEND}
There exists $\varepsilon>0$ such that for all $p \in (3,3+\varepsilon)$, the following holds. For all $\gamma \in (0, \frac{17}{8}-\frac{7}{4p})$, we have
%\begin{align*}
%\Vert (1+t)^{\gamma} \D^2 u \Vert_{\Ld^p(0,T;\Ld^p(\R^3_+))} \lesssim_{0,T_0,\gamma} \Psi \Big(\Vert u_0 \Vert_{\Ld^2(\R^3_+)}^2 + \Vert u_0 \Vert_{\Ld^1(\R^3_+)}^2 \Big)+1.
%\end{align*}
\begin{align*}
\Vert (1+t)^{\gamma} \D^2 u \Vert_{\Ld^p(0,T;\Ld^p(\R^3_+))} \lesssim_{0} 1.
\end{align*}
\end{coro}
\begin{proof}
We apply the maximal $\Ld^p \Ld^p$ regularity for the Stokes system (see Section \ref{AnnexeMaxregStokes} in the Appendix) satisfied by
$\U:=(1+t)^{\gamma} u$:
\begin{equation}
\left\{
      \begin{aligned}
         \partial_t U +A_p U &=S(u,f),\\
         U_{\mid x_3=0}&=0,\\
        U_{\mid t=0}&=u_0, 
    \end{aligned}
    \right.
\end{equation}
where $A_p$ stands for the Stokes operator in $\Ld^p_{\mathrm{div}}(\R^3_+)$. This entails
\begin{align*}
\Vert (1+t)^{\gamma} \D^2 u \Vert_{\Ld^p(0,T;\Ld^p(\R^3_+))} \lesssim \Vert u_0 \Vert_{\mathrm{D}_p^{1-\frac{1}{p},p}(\R^3_+)} +  \Vert S(u,f) \Vert_{\Ld^p(0,T;\Ld^p(\R^3_+))},
\end{align*}
where $\lesssim$ involves a universal constant and where $\mathrm{D}_p^{1-\frac{1}{p},p}(\R^3_+)$ has been defined in \eqref{domaineDqs}. Recall the meaning of the notation $\lesssim_0$ from Definition \ref{nota:lesssim_0}. Using Lemma \ref{LMsourceLplast}, we then have 
%\begin{align*}
%\Vert (1+t)^{\gamma} \D^2 u \Vert_{\Ld^p(0,T;\Ld^p(\R^3_+))} &\lesssim_{T_0,\gamma} \Vert u_0 \Vert_{\mathrm{D}_p^{1-\frac{1}{p},p}(\R^3_+)} + \mathrm{A}_0   \\
%& \qquad +\mathrm{A}_0 \Psi \Big(\Vert u_0 \Vert_{\Ld^2(\R^3_+)}^2 + \Vert u_0 \Vert_{\Ld^1(\R^3_+)}^2 \Big) \left[ 1+ \Vert (1+t)^{\gamma} \D^2 u \Vert_{\Ld^p(0,T;\Ld^p(\R^3))} \right].
%\end{align*}
\begin{align*}
\Vert (1+t)^{\gamma} \D^2 u \Vert_{\Ld^p(0,T;\Ld^p(\R^3_+))} &\leq \Vert u_0 \Vert_{\mathrm{D}_p^{1-\frac{1}{p},p}(\R^3_+)} \\
& \quad + \psi \Big( 1+\mathscr{E}(0) + \mathscr{N}_{q,m}(f_0)  \Big) \\
& \qquad \times  \left[ \Psi \left( \Vert u_0 \Vert_{\Ld^1 \cap \Ld^2(\R^3_+)} +N_q(f_0)+H_{0,k}(f_0) \right)+  \Vert  u_0 \Vert_{\H^1(\R^3_+)} \right] \\
 & \qquad \qquad \qquad \times \left[ 1+ \Vert (1+t)^{\gamma} \D^2 u \Vert_{\Ld^p(0,T;\Ld^p(\R^3_+))} \right].
\end{align*}
Using the smallness assumption \eqref{smallness:condition:THM}, we can ensure that
%\begin{align*}
%\mathrm{A}_0 \Psi \Big(\Vert u_0 \Vert_{\Ld^2(\R^3_+)}^2 + \Vert u_0 \Vert_{\Ld^1(\R^3_+)}^2 \Big)<\frac{1}{2}
%\end{align*}
\begin{align*}
 \psi \Big( 1+\mathscr{E}(0) + \mathscr{N}_{q,m}(f_0)  \Big)  \times \left[ \Psi \left( \Vert u_0 \Vert_{\Ld^1 \cap \Ld^2(\R^3_+)} +N_q(f_0)+H_{0,k}(f_0) \right)+  \Vert  u_0 \Vert_{\H^1(\R^3_+)} \right]  <\frac{1}{2},
\end{align*}
and this allows to absorb the term $\Vert (1+t)^{\gamma} \D^2 u \Vert_{\Ld^p(0,T;\Ld^p(\R^3_+))}$ in the l.h.s. This concludes the proof.
\end{proof}

\subsection{End of the proof of Theorem \ref{thm1}}
We are now able to close the bootstrap argument. Recall that $T \in (T_0,t^{\star})$, where $t^{\star}$ is given in Definition \ref{def:tstar}.
\begin{coro}\label{coro:CloseBoot}
The following inequalities hold
\begin{align*}
\int_{T_0}^T \Vert  u(s) \Vert_{\Ld^{\infty}(\R^3_+)} \, \mathrm{d}s &\lesssim_{0}  \Psi \left( \Vert u_0 \Vert_{\Ld^1 \cap \Ld^2(\R^3_+)} +N_q(f_0)+H_{0,k}(f_0) \right), \\
\int_{T_0}^T \Vert \nabla u(s) \Vert_{\Ld^{\infty}(\R^3_+)} \, \mathrm{d}s &\lesssim_{0}  \Psi \left( \Vert u_0 \Vert_{\Ld^1 \cap \Ld^2(\R^3_+)} +N_q(f_0)+H_{0,k}(f_0) \right).
\end{align*}
\end{coro}
\begin{proof}
We use the Gagliardo-Nirenberg-Sobolev inequality twice (see Section \ref{gagliardo-nirenberg} in the Appendix) to write that for $p>3$ (to be determined later)
\begin{align*}
\Vert u \Vert_{\Ld^{\infty}(\R^3_+)} &\lesssim \Vert \D^2 u \Vert_{\Ld^{p}(\R^3_+)}^{\alpha_p} \Vert u \Vert_{\Ld^{2}(\R^3_+)}^{1-\alpha_p}, \ \ \alpha_p:=\dfrac{3p}{7p-6}, \\[2mm]
\Vert \nabla u \Vert_{\Ld^{\infty}(\R^3_+)} &\lesssim \Vert \D^2 u \Vert_{\Ld^{p}(\R^3_+)}^{\beta_p} \Vert u \Vert_{\Ld^{2}(\R^3_+)}^{1-\beta_p}, \ \ \beta_p:=\dfrac{5p}{7p-6},
\end{align*}
where $\lesssim$ hides a universal constant. Thanks to the Hölder's inequality in time, we then get for all $\gamma >0$
\begin{align}\label{Holdu}
\begin{split}
\int_{T_0}^T \Vert  u(s) \Vert_{\Ld^{\infty}(\R^3_+)} \, \mathrm{d}s  &\lesssim \int_0^T \Vert \D^2 u(s) \Vert_{\Ld^{p}(\R^3_+)}^{\alpha_p} \Vert u(s) \Vert_{\Ld^{2}(\R^3_+)}^{1-\alpha_p} \, \mathrm{d}s \\
& \leq  \Vert (1+t)^{\gamma} \D^2 u \Vert_{\Ld^p(0,T;\Ld^p(\R^3_+))}^{\alpha_p} \left(  \int_0^T (1+t)^{-\gamma \frac{p \alpha_p }{p-\alpha_p}} \Vert u(s) \Vert_{\Ld^2(\R^3_+)}^{(1-\alpha_p)\frac{p}{p-\alpha_p}} \, \mathrm{d}s \right)^{\frac{p-\alpha_p}{p}},
\end{split}
\end{align}
as well as
\begin{align}\label{HoldNau}
\int_{T_0}^T \Vert  \nabla u(s) \Vert_{\Ld^{\infty}(\R^3_+)} \, \mathrm{d}s  \lesssim  \Vert (1+t)^{\gamma} \D^2 u \Vert_{\Ld^p(0,T;\Ld^p(\R^3_+))}^{\beta_p} \left(  \int_0^T (1+t)^{-\gamma \frac{p \beta_p }{p-\beta_p}} \Vert u(s) \Vert_{\Ld^2(\R^3_+)}^{(1-\beta_p)\frac{p}{p-\beta_p}} \, \mathrm{d}s \right)^{\frac{p-\beta_p}{p}}.
\end{align}
We also note that we have
\begin{align*}
\dfrac{p-\alpha_p}{p \alpha_p }=\dfrac{7}{3}-\dfrac{3}{p}, \ \ \dfrac{p-\beta_p}{p \beta_p }=\dfrac{7}{5}-\dfrac{11}{5p}.
\end{align*}
In view of Corollary \ref{coroMaxpEND}, we thus choose $p>3$ close to $3$ (remembering the assumption \eqref{data:hyp}) and take $\gamma$ close to $\frac{17}{8}-\frac{7}{4p}$ (which is strictly greater than $\frac{7}{3}-\frac{3}{p}$ and $\frac{7}{5}-\frac{11}{p}$ for $p$ close to $3$) so that 
\begin{align*}
\gamma \dfrac{p \alpha_p }{p-\alpha_p}>1, \ \ \gamma \dfrac{p \beta_p }{p-\beta_p}>1.
\end{align*}
Observe that we can first choose such $p$ and $\gamma$ and then perform the whole analysis of the two previous subsections 
%(using a finite number of exponents $q,m,r,\widetilde{q},\widetilde{m},\widetilde{r}$ involved in the different quantities $N_q(f_0), K_{m,r}(f_0), H_{\widetilde{q}}(f_0)$ and $F_{\widetilde{m},\widetilde{r}}(f_0)$, as well as a finite number of times the smallness assumption \eqref{smallness:condition:THM}, whether the estimates depend on $\gamma$ or not).
(using a finite number of quantities defined in Definition \ref{decay_f0}, as well as a finite number of times the smallness assumption \eqref{smallness:condition:THM}, whether the estimates depend on $\gamma$ or not). We then take the largest exponent involved in these quantities. We can now come back to \eqref{Holdu} and \eqref{HoldNau}, and use the uniform (in time) inequality $\Vert u(s) \Vert_{\Ld^2(\R^3_+)} \lesssim  \Psi \left( \Vert u_0 \Vert_{\Ld^1 \cap \Ld^2(\R^3_+)} +N_q(f_0)+H_{0,m}(f_0) \right)$ on $[0,T]$ with the previous exponents (ensuring uniform in time bound), therefore we get the desired results thanks to Corollary \ref{coroMaxpEND}.
\end{proof}
Recall that we have assumed that $t^{\star}<\infty$.
 \begin{proof}[End of the proof of Theorem \ref{thm1}]

As before, we use the notation $F:=j_f-\rho_f u$. Applying Corollary \ref{OKdecay:BrinkPoint} with $k>1$, we get
\begin{align*}
\int_0^{t^{\star}} \Vert F(s) \Vert_{\Ld^2(\R^3_+)}^2 \, \mathrm{d}s \lesssim_0 \varphi_2\big( N_q(f_0)+H_{0,k}(f_0) \big),
\end{align*}
while combining the Cauchy-Schwarz inequality with Corollary \ref{OKdecay:BrinkPoint} with $k>3$ (since $q \geq7$) yields
\begin{align*}
\int_0^{t^{\star}} \Vert F(s) \Vert_{\Ld^2(\R^3_+)}\, \mathrm{d}s & \leq \left( \int_0^{t^{\star}} \dfrac{\mathrm{d}s}{(1+s)^{2}} \right)^{1/2} \left( \int_0^{t^{\star}} (1+s)^{2} \Vert F(s) \Vert_{\Ld^2(\R^3_+)}^2 \, \mathrm{d}s   \right)^{1/2} \\
& \lesssim_0 \varphi_2\big( N_q(f_0)+H_{0,k}(f_0) \big).
\end{align*}
By the smallness assumption \eqref{smallness:condition:THM}, we can thus ensure that
\begin{align}
\Vert u_0 \Vert_{\H^1(\R^3_+)}^2 +  \int_0^{t^{\star}} \left[ \Vert F(s) \Vert_{\Ld^2(\R^3_+)}^2 + \Vert F(s) \Vert_{\Ld^2(\R^3_+)} \right] \, \mathrm{d}s < \dfrac{\mathrm{C}_{\star}}{2},
\end{align}
where $\mathrm{C}_{\star}$ refers to the universal constant from Proposition \ref{propdatasmall:VNSreg}. Since $F \in \Ld^2_{\mathrm{loc}}(\R^+;\Ld^2(\R^3_+)) \hookrightarrow \Ld^1_{\mathrm{loc}}(\R^+;\Ld^2(\R^3_+))$ by Proposition \ref{integBrinkman-1}, a continuity argument shows that there exists a strong existence time strictly larger than $t^{\star}$. Furthermore, by an appropriate choice of the data in \eqref{smallness:condition:THM}, we can use Corollary \ref{coro:CloseBoot} to ensure
\begin{align*}
\int_{T_0}^{t^{\star}} \Vert  u(s) \Vert_{\Ld^{\infty}(\R^3_+)} \, \mathrm{d}s < \dfrac{\delta}{4}, \ \ \int_{T_0}^{t^{\star}} \Vert \nabla u(s) \Vert_{\Ld^{\infty}(\R^3_+)} \, \mathrm{d}s < \dfrac{\delta}{2}.
\end{align*}
According to Proposition \ref{controlLinfiniLOC} and Proposition \ref{grad:nonunif}, this means that there exists a strong existence time $t>t^{\star}$ such that 
\begin{align*}
\int_{T_0}^{t} \Vert  u(s) \Vert_{\Ld^{\infty}(\R^3_+)} \, \mathrm{d}s <\dfrac{\delta}{2}, \ \ \int_{T_0}^{t} \Vert \nabla u(s) \Vert_{\Ld^{\infty}(\R^3_+)} \, \mathrm{d}s < \delta.
\end{align*}
This is a contradiction with the very definition of $t^{\star}$ therefore we necessarily have $t^{\star}=+\infty$. The proof of Theorem \ref{thm1} is now complete.
\end{proof}

\begin{appendix}

%\section*{\begin{center} ANNEXE \end{center} }

\section{Appendix}
\subsection{DiPerna-Lions theory in $\R^3_+\times \R^3$}\label{DiPernaLions}
\begin{thm}
Let $\chi \in \mathscr{C}^{\infty}(\R)$  such that $\vert \chi(z) \vert \leq \vert z \vert$ and $\chi' \in \Ld^{\infty}(\R)$. Take $f_0 \in \Ld^1 \cap \Ld^{\infty}(\R^3_+ \times \R^3)$, $G \in \R^3$  and a vector field $u \in \Ld^1_{\mathrm{loc}}(\R^+; \W^{1,1}(\R^3_+))$. Consider the following boundary value problem on $\R^3_+ \times \R^3$:
\begin{equation}\label{Vlasov:DPL}
\left\{
      \begin{aligned}
      \partial_t f +v\cdot \nabla_x f + {\rm div}_v(f\chi(u-v)+fG)&=0,\\
f_{\mid t=0}&=f_0,\\
f &=0, \ \mathrm{on} \ \Sigma^{-},
\end{aligned}
    \right.
\end{equation}
where $\chi(Z)=(\chi(Z_1), \chi(Z_2),\chi(Z_3))$ for all $Z \in \R^3$. Then we have, for all fixed $T>0$ 

\medskip

\textbullet \ \ \underline{Well-posedness}: There exists a unique 
$f \in \Ld^{\infty}_{\mathrm{loc}}(\R^+;\Ld^1 \cap\Ld^{\infty}(\R^3_+ \times \R^3)) $
which is a weak solution of the Cauchy problem \eqref{Vlasov:DPL}. Furthermore, $$ f \in \mathscr{C}(\R^+;\Ld^p_{\mathrm{loc}}(\overline{\R^3_+} \times \R^3)),$$ for all $p \in [1,\infty)$ and the function $f$ has a trace on $\partial \R^3_+ \times \R^3$ defined in the following sense: there exists a unique element $\gamma f \in \Ld^{\infty}_{\mathrm{loc}}([0,T] \times \partial \R^3_+ \times \R^3 )$ such that for any test function $\psi \in \mathscr{C}^{\infty}([0,T] \times \overline{\R^3_+ } \times \R^3)$ with compact support in space and velocity, and for all $0 \leq t_1 \leq t_2 \leq T$
\begin{multline*}
\int_{t_1}^{t_2} \int_{\R^3_+ \times \R^3}  f(t,x,v) \left[  \partial_t \psi + v \cdot \nabla_x \psi +(\chi(u(t,x)-v)+G ) \cdot \nabla_v \psi \right](t,x,v) \, \mathrm{d}x \, \mathrm{d} v \, \mathrm{d}t \\ 
=\int_{\R^3_+ \times \R^3} f(t_2,x,v) \psi(t_2,x,v) \, \mathrm{d}v \,\mathrm{d}x-\int_{\R^3_+ \times \R^3} f(t_1,x,v) \psi(t_1,x,v) \, \mathrm{d}x \,\mathrm{d}v \\
+ \int_{t_1}^{t_2}\int_{\partial \R^3_+ \times \R^3} \left[ (\gamma f) \psi(t,x,v) \right] v \cdot n(x) \,  \mathrm{d\sigma}(x) \, \mathrm{d}v  \, \mathrm{d}t.
\end{multline*}

%\medskip
%
%
%\textbullet \ \ \underline{Renormalization}: For every $\beta \in \mathscr{C}^1(\R)$, for all test function $\psi \in \mathscr{C}^{\infty}([0,T] \times \overline{\R^3_+ } \times \R^3)$ with compact support in space velocity, and for all $0 \leq t_1 \leq t_2 \leq T$, we have
%\begin{multline*}
%\int_{t_1}^{t_2} \int_{\R^3_+ \times \R^3}  \beta(f(t,x,v)) \left[  \partial_t \psi + v \cdot \nabla_x \psi +(u-v) \cdot \nabla_v \psi \right](t,x,v) \, \mathrm{d}v \, \mathrm{d} x \, \mathrm{d}t \\ 
%=\int_{\R^3_+ \times \R^3} \beta(f(t_2,x,v)) \psi(t_2,x,v) \, \mathrm{d}v \,\mathrm{d}x-\int_{\R^3_+ \times \R^3} \beta(f(t_1,x,v)) \psi(t_1,x,v) \, \mathrm{d}v \,\mathrm{d}x \\
%+ \int_{t_1}^{t_2}\int_{\partial \R^3_+ \times \R^3} \left[ \beta(\gamma f) \psi(t,x,v) \right] v \cdot n(x) \, \mathrm{d}v \, \mathrm{d\sigma}(x) \, \mathrm{d}t \\
%-3 \int_{t_1}^{t_2} \int_{\R^3_+ \times \R^3} \psi \left[f\beta'(f)-\beta(f) \right](t,x,v)\, \mathrm{d}v \, \mathrm{d} x \, \mathrm{d}t .
%\end{multline*}
\medskip

\textbullet \ \ \underline{Stability}: If
\begin{align*}
u_n \longrightarrow u  \ \ \text{in} \ \ \Ld^1_{\mathrm{loc}}(\R^+;\Ld^1(\R^3_+)) \ \text{ and } \
f_{0,n} \longrightarrow f_0  \ \ \text{in} \ \ \Ld^1_{\mathrm{loc}}(\R^3_+ \times \R^3),
\end{align*}
the corresponding sequence of solutions $(f_n)$ satisfies for all $p \in[1,\infty)$,
\begin{align*}
f_n \longrightarrow f  \ \ \text{in} \ \Ld^{\infty}_{\mathrm{loc}}(\R^+;\Ld^p(\R^3_+ \times \R^3)).
\end{align*}
\end{thm}
Such a result can be found in \cite[Theorem 3.2 - Proposition 3.2]{BGM} for the well-posedness and renormalization properties and in \cite[Theorem VI.1.9]{BF} for the stability property.

\subsection{The Cauchy problem for the Vlasov-Navier-Stokes system}\label{Section:PBCauchy}
Since the existence of global weak solutions to the Vlasov-Navier-Stokes system on a half-space with a gravity force has not been explicitly written in the literature, we provide some rather short elements of proof. We especially focus on the difficulties which may appear when obtaining the strong energy inequalities satisfied by these solutions (namely, the inequalities \eqref{ineq-energyNS} and \eqref{ineq-energy}).
\begin{thm}[Existence of weak solutions]\label{thm:existence}
Consider $(u_0,f_0)$ a pair of admissible initial conditions in the sense of Definition \ref{CIadmissible}. Then there exists a global weak solution $(u,f)$ in the sense of Definition \ref{sol-faible} to the system (\ref{eq:Vlasov})-(\ref{eq:NS2}) with boundary conditions (\ref{bcond-fluid})-(\ref{bcond-f}) and with initial data $(u_0,f_0)$.
\end{thm}
The proof of this result follows a now classic method for those non linear coupled problems: first, we introduce an approximated problem of the whole system. For each fixed $n \in \N\setminus \lbrace 0 \rbrace$ an appropriate fixed point procedure gives the existence of $(u_n,f_n)$ which are solutions to the following regularized system on each interval $[0,T]$
\begin{align}
\label{approxV:limitn}\partial_t f_n+v \cdot \nabla_x f_n    + {\rm div}_v   (f_n \chi_n( u_n-  v)  + f_n G ) &=0, \\
\label{approxNS:limitn}\partial_t u_n +( \J_n u_n \cdot \nabla) u_n- \Delta u_n + \nabla p_n &=  \int_{\R^3} f_n\chi_n(v- u_n) \, \mathrm{d}v, \\ 
\mathrm{div} \, u_n &=0,\\
(f_n,u_n)_{\mid t=0}&=(f_0^n,\J_n u_0),
\end{align}
where $T>0$ is arbitrary and where we consider the following operator, regularizations and data, and their respective convergences when $n$ goes to infinity: 
\begin{itemize}
\item $\chi_n \underset{n \rightarrow + \infty}{\longrightarrow} \mathrm{Id}_{\R}$, where for all $n \in \N\setminus \lbrace 0 \rbrace$, $ \chi_n \in \mathscr{D}(\R)$ is an odd, increasing and bounded function satisfying $\vert \chi(z) \vert \leq \vert  z \vert$, $\Vert \chi_n' \Vert_{\Ld^{\infty}(\R)} \leq 1$ and $z \chi_n(z) \geq 0$. Above, $\chi_n$ is applied componentwise.
\item $f_0^n \underset{n \rightarrow + \infty}{\longrightarrow} f_0$ strongly in $\Ld^p(\R^3_+\times \R^3)$ for all $1 \leq p < \infty$ and weakly-$\ast$ in $\Ld^{\infty}(\R^3_+\times \R^3)$, and where $f_0^n=\eta^n f_0$ with $(\eta^n)_n$ a family of  positive functions compactly supported in velocity such that $0 \leq \eta^n \leq 1$ and increasing to $1$.
\item for all $n \in \N\setminus \lbrace 0 \rbrace$, the operator $ \mathrm{J}_n$ refers to the Yosida approximation of the identity defined by
\begin{align*}
 \mathrm{J}_n:=(\mathrm{I}+n^{-1}A_2 )^{-1},
 \end{align*}
where $A_2$ is the Stokes operator on $\Ld^2(\R^3_+)$ (see Section \ref{AnnexeMaxregStokes} in the Appendix). The operator $\mathrm{J}_n$ acts on $\Ld^2_{\mathrm{div}}(\R^3_+)$ and satisfies: for all $v \in \Ld^2_{\mathrm{div}}(\R^3_+)$, $ \mathrm{J}_n u \in D(A_2),\Vert \J_n v \Vert_{\Ld^2(\R^3_+)} \leq \Vert v \Vert_{\Ld^2(\R^3_+)} $ and $\J_n v \underset{n\rightarrow + \infty}{\longrightarrow} v$ in $\Ld^2(\R^3_+)$.
\end{itemize}
In the previous construction, $u_n$ is a strong solution to the Navier-Stokes equations while $f_n$ is a weak solution to the Vlasov equation which is compactly supported in velocity. Note that the presence of the gravity force is harmless in this procedure.

Then, we pass to the limit in the previous approximated problem by compactness arguments, in order to recover solution to the whole system of equations.
To do so, we shall obtain energy estimates which are independent of $n$: all in all, it can be proven that $(u_n)_n$ is weakly compact in $\Ld^{\infty}_{\mathrm{loc}}(\R^+;\Ld^{2}(\R^3_+)) \cap \Ld_{\mathrm{loc}}^2(\R^+;\H^1(\R^3_+))$ while $(f_n)$ is weakly-$\ast$ compact in $ \Ld^{\infty}_{\mathrm{loc}}(\R^+ ;\Ld^{\infty}(\R^3_+ \times \R^3))$. Furthermore, a standard application of the Aubin-Lions lemma ensures that $(u_n)_n$ is strongly compact in $\Ld^{2}_{\mathrm{loc}}(\R^+;\Ld^{2}_{\mathrm{loc}}(\R^3_+))$ and this also implies that, up to an additional extraction, $(u_n)_n$ converges almost everywhere. Note that the treament of the gravity term $G$ requires an application of Grönwall's lemma for the quantity $M_2 f_n$, which turns out to be bounded itself in $\Ld^{\infty}_{\mathrm{loc}}(\R^+)$. Theses convergences are enough to pass to the limit in the weak formulation of the Vlasov equation and of the Navier-Stokes equations because each term converges in the sense of distributions (see e.g. \cite{BMM} for a proof, in a more involved context). An additional diagonal extraction along increasing intervals of the type $[0,T_j]$ with $T_j \rightarrow + \infty$ allows to obtain global solutions.

\bigskip

We then explain how one can recover the strong energy inequality \eqref{ineq-energyNS} for the Navier-Stokes equations, as well as the strong energy inequality \eqref{ineq-energy} for the Vlasov-Navier-Stokes system: since the case $s=0$ is classic (at least for the Navier-Stokes equations with a given source term, see e.g \cite{BF}), we restrict ourselves to the case $s>0$ and we pay attention to the particular treatment of the coupling terms.

We first write the strong energy inequality \eqref{ineq-energyNS} for the Navier-Stokes equations. It seems that the construction of solutions satisfying this strong energy inequality (and not only the same inequality with $s=0$ and $t \geq 0$) has long remained an open problem for general unbounded domains, due to the interaction of the pressure with the distant boundaries \cite{HeywoodOPEN} (see however a general modern treatment in \cite{Far}). Nevertheless, the strategy applied in \cite{MS} for exterior domains actually extends for the half-space case.
%(note that, to our knowledge, this fact is only mentioned in \cite{MS} and not in the subsequent literature).
%, and that the reference for the half-space case given \cite{GaldINTRO} seems uncorrect).
We thus refer to \cite{MS} where the proof of \eqref{ineq-energyNS} for the Navier-Stokes equations with a forcing term is written down explicitly: namely, it consists in obtaining a localised energy inequality between almost any times $s<t$, by multiplying the regularized version of the Navier-Stokes equations by $2 u_n  \varphi_k$ (where $\varphi_k:=\varphi(\cdot/k)$ with $\varphi \in \mathscr{D}(\R^3)$ satisfying $0 \leq \varphi \leq 1$, $\varphi(x)=1$ for $\vert x \vert \leq 1$ and $\varphi(x)=0$ for $ \vert x \vert \geq 2$), integrating on $(s,t) \times \R^3_+$, performing integration by parts, and then understanding the behavior of the remaining terms, especially for the pressure ones. All in all, this reads as 
\begin{align}\label{locIneqNS}
\begin{split}
\int_{\R^3_+} \vert u_n(t) \vert^2 \varphi_k \, \mathrm{d}x   + 2\int_s^t \int_{\R^3_+}\vert \nabla u_n \vert^2 \varphi_k \, \mathrm{d}x \, \mathrm{d}\tau\leq \mathrm{R}_{s,t}^{k,n} +\int_{\R^3_+} \vert u_n(s) \vert^2 \varphi_k \, \mathrm{d}x+ 2\int_s^t \int_{\R^3_+} \varphi_k u_n \cdot S_{\chi_n}(u_n,f_n) \, \mathrm{d}x \, \mathrm{d}\tau,
\end{split}
\end{align}
with
\begin{align*}
S_{\chi_n}(u_n,f_n):=\int_{\R^3} f_n \chi_n(v-u_n) \, \mathrm{d}v,
\end{align*}
and where $\mathrm{R}_{s,t}^{k,n}$ is a remainder term such that
%\begin{align*}
$\underset{k\rightarrow + \infty}{\lim} \underset{n\rightarrow + \infty}{\limsup} \, \mathrm{R}_{s,t}^{k,n} =0$
%\end{align*}
(see \cite{MS} for a proof). 

Thus, the main difficulty essentially lies in the treatment of the last term in the right-hand side of the inequality \eqref{locIneqNS}. First, since $(u_n)_n$ is bounded in $\Ld^{\infty}_{\mathrm{loc}}(\R^+;\Ld^{2}(\R^3_+)) \cap \Ld_{\mathrm{loc}}^2(\R^+;\Ld^6(\R^3_+))$, an interpolation argument shows that $(u_n)_n$ is bounded in $\Ld^{r}_{\mathrm{loc}}(\R^+;\Ld^{s}(\R^3_+))$ if 
\begin{align*}
\frac{2}{r}+\frac{3}{s}=\frac{3}{2}, \ \ 2 \leq s \leq 6, \ \ 2 \leq r \leq \infty.
\end{align*}
For every $n$, we have in particular $u_n \in \Ld^{20/9}_{\mathrm{loc}}(\R^+;\Ld^{5}(\R^3_+)) \hookrightarrow \Ld^{19/9 }_{\mathrm{loc}}(\R^+;\Ld^{5}(\R^3_+))$. Furthermore, there exist $r_1>19/9$ and $s_1>5$ such that $(u_n)_n$ is bounded in $\Ld^{r_{1}}_{\mathrm{loc}}(\R^+;\Ld^{s_{1}}(\R^3_+))$. Since $(u_n)$ converges almost everywhere towards $u$, a corollary of Vitali convergence theorem shows that $u_n \underset{n \rightarrow + \infty}{\longrightarrow} u$ in $\Ld^{19/9 }_{\mathrm{loc}}(\R^+;\Ld^{5}_{\mathrm{loc}}(\R^3_+))$.
%pour $s=5+\varepsilon$, on a $r_{\varepsilon}=\frac{4\varepsilon+20}{3 \varepsilon+9}$ qui décroît en $\varepsilon$ en partant de $20/9$ en $\varepsilon=0$.

We also know that  $(M_2 f_n)_n$ is bounded in $\Ld^{\infty}_{\mathrm{loc}}(\R^+)$ so that Proposition \ref{interpo-moment} and the maximum principle \eqref{max:principle} ensure that $(j_{f_n})_n$ in bounded in $\Ld^{\infty}_{\mathrm{loc}}(\R^+; \Ld^{5/4}(\R^3_+))$ while $(\rho_{f_n})_n$ is bounded in $\Ld^{\infty}_{\mathrm{loc}}(\R^+; \Ld^{5/3}(\R^3_+)) \cap \Ld^{\infty}_{\mathrm{loc}}(\R^+; \Ld^{1}(\R^3_+))$, and thus bounded in $\Ld^{\infty}_{\mathrm{loc}}(\R^+; \Ld^{p_{\theta}}(\R^3_+))$ where $p_{\theta}=\frac{5}{3+2\theta}$ with $\theta \in [0,1]$. As there exists $s_2<6$ close enough to $6$ and $\theta \in [0,1]$ such that $p_{\theta}^{-1}+s_2^{-1}=4/5$ and $(u_n)_n$ is bounded in $\Ld^{r_2}_{\mathrm{loc}}(\R^+;\Ld^{s_{2}}(\R^3_+))$ for some $r_2>2>19/10=(19/9)'$, Hölder's inequality shows that $(\rho_{f_n} u_n)_n$ is bounded in $\Ld^{19/10}_{\mathrm{loc}}(\R^+;\Ld^{5/4}(\R^3_+))$, as well as $(j_{f_n})_n$.
%Avec $s_2=6_\varepsilon$, on trouve $p_{\theta}$ associé pour $\theta=(\varepsilon-1)/(2(\varepsilon-6))$ and $r_2=(24-4\varepsilon)/(12-3\varepsilon)>2$ croissant en $\varepsilon$.
By the properties of $\chi_n$, we observe that for all $n$
\begin{align*}
\vert S_{\chi_n}(u_n,f_n) \vert \leq j_{f_n} + \vert u_n \vert \rho_{f_n},
\end{align*}
therefore the previous bound entails that $(S_{\chi_n}(u_n,f_n))_n$ converges weakly in $\Ld^{19/10}_{\mathrm{loc}}(\R^+;\Ld^{5/4}(\R^3_+))$ to some limit $S$, up to extraction. Since $(S_{\chi_n}(u_n,f_n))_n$  also converges to $j_f - \rho_f u$ in the sense of distributions, we can identify $S=j_f-\rho_f u$.

All in all, combining the strong convergence of $(u_n)_n$ and the weak convergence of $(S_{\chi_n}(u_n,f_n))_n$ we have just obtained, we get for all fixed $k$
\begin{align*}
\int_s^t \int_{\R^3_+} \varphi_k u_n \cdot S_{\chi_n}(u_n,f_n) \, \mathrm{d}x \, \mathrm{d}\tau \underset{n \rightarrow + \infty}{\longrightarrow} \int_s^t \int_{\R^3_+} \varphi_k u \cdot (j_f-\rho_f u) \, \mathrm{d}x \, \mathrm{d}\tau,
\end{align*}
and by the dominated convergence theorem, we can then pass to the limit when $k \rightarrow + \infty$ in the inequality \eqref{locIneqNS}.

In this procedure, one has to first pick one $s>0$ such that $u_n(s) \underset{n \rightarrow + \infty}{\longrightarrow} u(s)$ in $\Ld^2(K)$ for any compact $K \subset \R^3$ and the inequality \eqref{ineq-energyNS} eventually holds for almost any $s>0$ (including $s=0$) and almost any $t \geq s$. Standard lower-semicontinuity arguments allow one to extend this inequality to any $t \geq s$.

\medskip

Finally, we come back to the proof of the energy inequality \eqref{ineq-energy} for the Vlasov-Navier-Stokes system: in view of \eqref{ineq-energyNS}, we only need to prove that for all $0 \leq s \leq t$
\begin{align}\label{2ndMomineqVlasov}
M_2 f(t)+ 2 \int_s^t M_2 f(\tau) \, \mathrm{d}\tau \leq M_2 f(s) + 2 \int_s^t \int_{\R^3_+ \times \R^3} f(\tau,x,v) v \cdot [u(t,x)+G] \, \mathrm{d}x \, \mathrm{d}v \, \mathrm{d}\tau. 
\end{align}
We proceed as follows: mimicking the beginning of the proof of Lemma \ref{interpo-estimate} for $\alpha=2$ (which does not use the energy inequality \eqref{ineq-energy} but only the regularity of $u$), we obtain the fact that $(t,x,v) \mapsto \vert v \vert^2 f \in \Ld^{\infty}_{\mathrm{loc}}(\R^+;\Ld^1(\R^3_+ \times \R^3))$. We now consider $\psi_n(v):=\theta_n(\vert v \vert^2)$ where $\theta_n(z):= n\theta(z/n)$ is defined thanks to a positive function $\theta \in \mathscr{D}(\R^+)$ equal to the identity function on $[0,1]$, less than this function on $(1,+\infty)$ and with a derivative uniformly bounded by $1$. After a suitable localization in space that we do not detail here, the function $\psi_n$ is an admissible test function in the weak formulation of the Vlasov equation with a trace (see Section \ref{DiPernaLions}) and this implies
% la fonction test en x est dans D(adhérence du demi-plan), ie fonction test de tout l'espace donc le support déborde un peu en dehors du demi plan: on peut ndonc prendre =1 sur le demi plan et c'est tout.
\begin{align*}
\int_{s}^{t} \int_{\R^3_+ \times \R^3}  f(\tau,x,v) 2 \theta'\left(\frac{\vert v \vert^2}{n} \right)& v \cdot \left[ (u(t,x)-v+G) \right] \, \mathrm{d}x \, \mathrm{d}v \, \mathrm{d}\tau \\ 
&=\int_{\R^3_+ \times \R^3} f(t,x,v) \theta_n(\vert v \vert^2) \, \mathrm{d}v \,\mathrm{d}x-\int_{\R^3_+ \times \R^3} f(s,x,v) \theta_n(\vert v \vert^2) \, \mathrm{d}x \, \mathrm{d}v \\
& \qquad \qquad+ \int_{s}^{t}\int_{\Sigma^+} \left[ (\gamma f) \theta_n(\vert v \vert^2) \right] v \cdot n(x) \, \mathrm{d\sigma}(x) \, \mathrm{d}v \,  \mathrm{d}\tau \\
&\geq \int_{\R^3_+ \times \R^3} f(t,x,v) \theta_n(\vert v \vert^2) \, \mathrm{d}v \,\mathrm{d}x-\int_{\R^3_+ \times \R^3} f(s,x,v) \theta_n(\vert v \vert^2) \, \mathrm{d}x \, \mathrm{d}v,
\end{align*}
because of the boundary condition \eqref{bcond-f}. Since $\theta_n(\vert v \vert^2) \underset{n \rightarrow +\infty}{\longrightarrow} \vert v \vert^2$ with $\theta_n(\vert v \vert^2) \leq \vert v \vert^2$, and $\theta'(\vert v \vert^2/n) \underset{n \rightarrow +\infty}{\longrightarrow} 1$ with $\theta'(\vert v \vert^2/n) \leq 1$, we can use the dominated convergence theorem to pass to the limit when $n \rightarrow +\infty$ in the previous inequality: indeed, $M_2 f \in \Ld^1_{\mathrm{loc}}(\R^+)$ while $u \in \Ld^2_{\mathrm{loc}}(\R^+;\Ld^5(\R^3_+))$ and $m_1 f \in \Ld^2_{\mathrm{loc}}(\R^+;\Ld^{5/4}(\R^3_+))$ by interpolation.
%\textcolor{red}{\textbf{A cacher :}} On majore par $f \vert v \vert^2$ dans le r.h.s et OK pour le TCD. Dans le l.h.s, on majore en valeur absolue par
%$$ \vert f v \cdot (u-v+g) \vert \leq f \vert v \vert \vert u \vert + \vert v \vert^2 f + \vert v \vert g f=\mathrm{I}+\mathrm{II}+\mathrm{III}$$
%Le terme $\mathrm{II}$ est $\Ld^1((s,t)\times \R^3_+ \times \R^3)$ car $M_2 f$ est localement intégrable sur $\R^+$. On peut traiter $\mathrm{III}$ de la même manière en écrivant $\vert v \vert  f \leq (1+\vert v \vert^2)f$ qui est lui aussi dans $\Ld^1((s,t)\times \R^3_+ \times \R^3)$. Pour $\mathrm{I}$, on montre qu'il est dans $L^1((s,t)\times \R^3_+ \times \R^3)$ car
%$$\int_s^t \int_{\R^3_+ \times \R^3}  \vert u \vert \vert v  \vert f \leq \int_s^t \int_{\R^3_+} \underbrace{\vert u \vert}_{\in \Ld^2_t \Ld^{5}_x} \underbrace{m_1 f}_{\in \Ld^2_t \Ld^{5/4}_x}<\infty, \ \ \text{by Hölder}$$
%car en effet $\Vert m_1 f(\tau) \Vert_{\Ld^{5/4}(\R^3_+)}\lesssim \Vert f \Vert_{\Ld^{\infty}_{t,x,v}}^{1/5} M_2 f(\tau)^{4/5}$.
Adding \eqref{ineq-energyNS} to \eqref{2ndMomineqVlasov}, we eventually get the energy inequality \eqref{ineq-energy}.

\subsection{Gagliardo-Nirenberg-Sobolev inequality on $\R^3_+$}
\begin{thm}\label{gagliardo-nirenberg}
Let $1 \leq p,q,r \leq \infty$ and $m \in \N$. Suppose $j \in \N$ and $\alpha \in [0,1]$ satisfy the relations
\begin{align*}
&\dfrac{1}{p}=\dfrac{j}{3}+\left( \dfrac{1}{r}-\dfrac{m}{3} \right)\alpha+\dfrac{1-\alpha}{q},\\
&\dfrac{j}{m} \leq \alpha \leq 1,
\end{align*}
with the exception $\alpha<1$ if $m-j-d/r \in \N$. 

Then for all $g \in \Ld^q(\R^3_+)$, if $\mathrm{D}^m g \in \Ld^r(\R^3_+)$, we have $\mathrm{D}^j g \in \Ld^p(\R^3_+)$ with the estimate 
$$ \Vert \mathrm{D}^j g  \Vert_{\Ld^p(\R^3_+)} \lesssim \Vert \mathrm{D}^m g \Vert_{\Ld^r(\R^3_+)} ^{\alpha} \Vert g \Vert_{\Ld^q(\R^3_+)}^{1-\alpha},$$
where $\lesssim$ refers to a constant only depending  on the dimension $d$.
\end{thm}
This result can be found in \cite[Thm 1.5.2]{CherMil}, at least for the case of the whole space $\R^3$. In the case of the half-space $\R^3_+$, we can rely on the existence of an extension operator mapping functions defined on $\R^3_+$ to functions defined on $\R^3$ and which is continuous for the topology of Sobolev spaces. More precisely, for all $m \in \N \setminus \left\lbrace 0 \right\rbrace$ and $ \ell \in [1,\infty)$, there exists an operator $E : \W^{m,\ell}(\R^3_+) \longrightarrow \W^{m,\ell}(\R^3)$ (which is independent of $\ell$) such that for all $v \in \W^{m,\ell}(\R^3_+)$, we have
\begin{itemize}
\item $(E v)_{\mid \R^3_+}=v$,
\item for all $i \in \left\lbrace 0, \cdots, m \right \rbrace$, $\Vert \mathrm{D}^i(E v)  \Vert_{\Ld^{\ell}(\R^3)} \lesssim \Vert \mathrm{D}^i v  \Vert_{\Ld^{\ell}(\R^3_+)}$.
\end{itemize}
This result can be found in (the proofs of) \cite[Section 2.3.3]{Demen} and be combined with the Gagliardo-Nirenberg-Sobolev inequality on $\R^3$ to deduce Theorem \ref{gagliardo-nirenberg}.

\subsection{Maximal $\Ld^p \Ld^q$ regularity for the Stokes system on $\R^3_+$}\label{AnnexeMaxregStokes}

Let $1 < q <\infty$ be fixed. Given a vector field $u \in \Ld^q(\R^3_+)$, this can be uniquely decomposed as
\begin{align*}
& u=\widetilde{u} + \nabla p, \\
& \widetilde{u}  \in \Ld^q_{\mathrm{div}}(\R^3_+), \ p \in \Ld^q(\R^3_+), \ \nabla p \in \Ld^q(\R^3_+),
\end{align*}
where $\Ld^q_{\mathrm{div}}(\R^3_+)$ stands for the closure in $\Ld^q(\R^3_+)$ of $\mathscr{D}_{\mathrm{div}}(\R^3_+)$.
We recall that the projection $\mathbb{P}_q : u \mapsto \widetilde{u} $ is continuous from $\Ld^q(\R^3_+)$ to $\Ld^q_{\mathrm{div}}(\R^3_+)$.

For $1<q<\infty$, we consider the following Stokes operator
\begin{align*}
A_q := -\mathbb{P}_q \Delta , \ \  D(A_q):=\Ld^q_{\mathrm{div}}(\R^3_+) \cap \W^{1,q}_0(\R^3_+) \cap \W^{2,q}(\R^3_+).
\end{align*}
We also set
\begin{align}\label{domaineDqs}
\mathrm{D}_q^{1-\frac{1}{s},s}(\R^3_+):=\left( D(A_q),\Ld^q_{\mathrm{div}}(\R^3_+) \right)_{1/s,s},
\end{align}
where $( \ , \ )_{1/s,s}$ refers to the real interpolation space of exponents $(1/s,s)$. In the case of the Stokes operator $A_q$, which generates an analytic semigroup $e^{-tA_q}$, the quantity 
\begin{align}
\Vert u \Vert_{\Ld^q(\R^3_+)}+ \left( \int_0^{\infty} \Vert A_q e^{-t A_q} u \Vert_{\Ld^q(\R^3_+)}^s \, \mathrm{d}t  \right)^{1/s}
\end{align}
defines an equivalent norm on $\mathrm{D}_q^{1-\frac{1}{s},s}(\R^3_+)$ (see \cite[Chapter 5]{Lunardi}).

The main result is the following and can be found with further references in \cite{Giga}.

\begin{thm}
Consider $0 < T \leq \infty$ and $1 <q,s <\infty$. Then, for every $u_0 \in \mathrm{D}_q^{1-\frac{1}{s},s}(\R^3_+)$ which is divergence free and $f \in \Ld^s(0,T;\Ld^q_{\mathrm{div}}(\R^3_+))$, there exists a unique solution $u$ of the Stokes system
\begin{equation*}
\left\{
      \begin{aligned}
\partial_t u +A_q u &=f, \\
u_{\mid x_3=0}&=0, \\
u(0,x)&=u_0(x),
      \end{aligned}
    \right.
\end{equation*}
satisfying 
\begin{align*}
& u \in \Ld^s(0,T';D(A_q)) \ \text{for all finite } T' \leq T,
%& \partial_t u \in \Ld^s(0,T;\Ld^q(\R^3_+)),
\end{align*}
and
\begin{align*}
\Vert \partial_t u \Vert_{\Ld^s(0,T;\Ld^q(\R^3_+))}  + \Vert \mathrm{D}^2 u\Vert_{\Ld^s(0,T;\Ld^q(\R^3_+))}\leq C\left(\Vert u_0 \Vert_{\mathrm{D}_q^{1-\frac{1}{s},s}(\R^3_+)} + \Vert f \Vert_{\Ld^s(0,T;\Ld^q(\R^3_+))} \right),
\end{align*}
where $C=C(q,s)>0$.

Furthermore, if $u_0  \in \W^{1,q}_0(\R^3_+)\cap \Ld^q_{\mathrm{div}}(\R^3_+)$ and if $s \in (1,2]$, the statement holds and we can replace $\Vert u_0 \Vert_{\mathrm{D}_q^{1-1/s,s}(\R^3_+)}$ by $\Vert  u_0 \Vert_{\W^{1,q}_0(\R^3_+)}$ in the right hand side of the previous inequality.

\end{thm}

\medskip
For the sake of completeness, we bring some precisions concerning the last statement of the theorem (even if a related fact can be found in \cite[Remark 2.5]{Giga}). Suppose that $u \in D(A_q^{\frac{1}{2}}) \cap \Ld^q_{\mathrm{div}}(\R^3_+)$.  If $s \in (1,2)$, we write $\frac{1}{2}=1-\frac{1}{s}+\varepsilon$ with $\varepsilon=\frac{2-s}{2s}>0$ (so that $1-s\varepsilon<1$). Since $A_q=A$ generates an analytic semigroup, we have by standard functional calculus manipulations (see e.g. \cite[Chapter 3]{Haase})
\begin{align*}
 \int_0^{\infty} \Vert A e^{-t A} u \Vert_{\Ld^q(\R^3_+)}^s \, \mathrm{d}t  &=  \int_0^{1} \Vert A e^{-t A} u \Vert_{\Ld^q(\R^3_+)}^s \, \mathrm{d}t  +  \int_1^{\infty} \Vert A e^{-t A} u \Vert_{\Ld^q(\R^3_+)}^s \, \mathrm{d}t \\
&=  \int_0^{1} \Vert A^{\frac{1}{s}-\varepsilon} e^{-tA} A^{\frac{1}{2}}u  \Vert_{\Ld^q(\R^3_+)}^s \, \mathrm{d}t  +  \int_1^{\infty} \Vert A e^{-t A} u \Vert_{\Ld^q(\R^3_+)}^s \, \mathrm{d}t \\
& \lesssim \Vert A^{\frac{1}{2}} u  \Vert_{\Ld^q(\R^3_+)}^s \int_0^{1} \frac{1}{t^{1-s\varepsilon}} \, \mathrm{d}t  +  \Vert u \Vert_{\Ld^q(\R^3_+)}^s \int_1^{\infty}   \frac{ 1}{t^s} \, \mathrm{d}t<\infty,
\end{align*}
where we have used the fact that the function $x \mapsto x^{\frac{1}{s}-\varepsilon}e^{-tx}$ and $x \mapsto xe^{-tx}$ are respectively bounded on $\R^+$ by (a constant times) $t^{\varepsilon-\frac{1}{s}}$ and $t^{-1}$.

In the case where $s=2$, we proceed in a similar way and write 
\begin{align*}
 \int_0^{\infty} \Vert A e^{-t A} u \Vert_{\Ld^q(\R^3_+)}^2 \, \mathrm{d}t  & \leq \int_0^{1} \Vert A^{\frac{1}{2}} e^{-tA} A^{\frac{1}{2}}u  \Vert_{\Ld^q(\R^3_+)}^s \, \mathrm{d}t  +  \int_1^{\infty} \Vert A e^{-t A} u \Vert_{\Ld^q(\R^3_+)}^2 \, \mathrm{d}t \\
 & \lesssim \Vert A^{\frac{1}{2}} u  \Vert_{\Ld^q(\R^3_+)}^2 \int_0^{1} \frac{e^{-\frac{1}{t}}}{t} \, \mathrm{d}t  +  \Vert u \Vert_{\Ld^q(\R^3_+)}^2 \int_1^{\infty}   \frac{ 1}{t^2} \, \mathrm{d}t<\infty,
\end{align*}
where we have used the fact that the function $x \mapsto \sqrt{x}e^{-tx}$ is bounded by $t^{-1/2}e^{-\frac{1}{2t}}$.

Thus, the conclusion follows because $D(A_q^{\frac{1}{2}})=\W^{1,q}_0(\R^3_+)\cap \Ld^q_{\mathrm{div}}(\R^3_+)$ (see \cite{BM}).

\subsection{Conditional decay of the energy: proof of Theorem \ref{cond:decay}}\label{Section:CondDecayAPPENDIX}
Let $u$ be a Leray solution (with strong energy inequality \eqref{ineq-energyNS-decay}) to the Navier-Stokes system \eqref{systNS:decay} with source term $F$ and initial data $u_0$. We provide some elements of proof concerning the conditional decay of the $\Ld^2$ norm of $u$ on $[0,T]$, provided that the forcing term $F$ satisfies
\begin{align}\label{HYPdecay:BRINK2}
\forall t \in [0,T], \ \ \Vert F(t) \Vert_{\Ld^2(\R^3_+)} \leq \dfrac{\mathrm{C}}{(1+t)^{7/4}},
\end{align}
for some constant $\mathrm{C}>0$. As explained in Section \ref{Section:CondDecay}, this enters within the scope of the situation treated in \cite{Wieg}, with the adaptation to the half-space case coming from \cite{BM}.

\medskip

We first recall the following $\Ld^r \Ld^q$ estimates for the semigroup generated by the Stokes operator $A_q$ on $\Ld^q_{\mathrm{div}}(\R^3_+)$ (see Section \ref{AnnexeMaxregStokes} in the Appendix). A proof can be found in \cite{BM}.
\begin{lem}\label{decayStokes:LpLq}
Let $a \in \Ld^2_{\mathrm{div}}(\R^3_+) \cap \Ld^r(\R^3_+)$ for some $r \in [1,\infty]$. Then for all $t \geq 0$
\begin{align}
\Vert e^{-t A_q} a \Vert_{\Ld^q(\R^3_+)} \leq  \dfrac{C}{t^{\frac{3}{2}(\frac{1}{r}-\frac{1}{q})}}\Vert  a \Vert_{\Ld^r(\R^3_+)},
\end{align}
where $C$ in independent of $a$ and $t$, provided either $1<r \leq q<\infty$, or $1 \leq r <q  \leq \infty$.
\end{lem}

\begin{proof}[Proof of Theorem \ref{cond:decay}]
Consider a  given time-dependent smooth positive cut-off function $g : \R^+ \rightarrow (0,+\infty)$. Since the Stokes operator $A=A_2$ on $\Ld^2_{\mathrm{div}}(\R^3_+)$ is a self-adjoint and positive operator (see Section \ref{AnnexeMaxregStokes} in the Appendix), we can use a spectral decomposition of this operator thanks to a resolution of the identity $(E_{\lambda})_{\lambda \geq 0}$ (see \cite[Section 3.2]{Sohr}), namely
\begin{align*}
A=\int_0^{+\infty} \lambda \, \mathrm{d}E_{\lambda}.
\end{align*}
By standard functional calculus rules and the Plancherel-Parseval theorem, we can write for all $\tau \in [0,T]$
$$
\begin{aligned}
\int_{\R^3_+} | \nabla u(\tau,x) |^2 \, \mathrm{d} x &= \Vert A^{1/2}u(\tau) \Vert_{\Ld^2(\R^3_+)}^2 = \int_0^{+\infty} \lambda \,  \Vert \mathrm{d}E_{\lambda}u(\tau)\Vert_{\Ld^2(\R^3_+)}^2\\
&\geq  \int_{\sqrt{\lambda}>g(\tau)}\lambda \,  \Vert \mathrm{d}E_{\lambda}u(\tau)\Vert_{\Ld^2(\R^3_+)}^2 \\
&\geq g^2(\tau) \| u(\tau)\|_{\Ld^2(\R^3_+)}^2 - g^2(\tau) \int_0^{g^2(\tau)}  \,  \Vert \mathrm{d}E_{\lambda}u(\tau)\Vert_{\Ld^2(\R^3_+)}^2.
\end{aligned}
$$
%Thanks to the fact that 
%\begin{align*}
%2\int_s^t \langle F(\tau), u(\tau) \rangle \, \mathrm{d}\tau \leq \int_s^t g^{-2}(\tau) \Vert F(\tau) \Vert_{\Ld^2(\R^3_+)}^2 \, \mathrm{d}\tau + \int_s^t g^{2}(\tau) \Vert u(\tau) \Vert_{\Ld^2(\R^3_+)}^2 \, \mathrm{d}\tau,
%\end{align*}
Then, we can combine the strong energy inequality (\ref{ineq-energyNS}) for the Navier-Stokes equations with the Cauchy-Schwarz inequality to get the following key inequality: for all $t \in [0,T]$ and almost every $s \in [0,t]$ (including $s=0$)
\begin{multline}\label{key}
\Vert u(t) \Vert_{\Ld^2(\R^3_+)}^2 + \int_s^t g^2(\tau) \Vert u(\tau) \Vert_{\Ld^2(\R^3_+)}^2 \, \mathrm{d} \tau   \leq \Vert u(s) \Vert_{\Ld^2(\R^3_+)}^2 + 
 \int_s^t g^{-2}(\tau) \Vert F(\tau) \Vert_{\Ld^2(\R^3_+)}^2 \, \mathrm{d}\tau \\
 +2\int_s^t g^2(\tau) \int_0^{g^2(\tau)}  \,  \Vert \mathrm{d}E_{\lambda}u(\tau)\Vert_{\Ld^2(\R^3_+)}^2 \, \mathrm{d} \tau.
\end{multline}
Before going further, we introduce the notation $E_{\leq g^2(\tau)}:=\mathbf{1}_{[0,g^2(t)]}(A)$ (in the sense of the bounded functional calculus).
%\begin{align*}
%E_{\leq g^2(\tau)}:=\mathbf{1}_{[0,g^2(t)]}(A)=\int_0^{g^2(\tau)} \mathrm{d}E_{\lambda}.
%\end{align*}
In order to estimate the last term in the right-hand side of \eqref{key}, we fix $ \tau \in [s,t]$ and follow \cite{BM-goodGron}: we choose $\psi \in \mathscr{D}_{\mathrm{div}}(\R^3_+)$ and then take $\varphi(\sigma):=e^{-(\tau-\sigma)A} E_{\leq g^2(\tau)}\psi$ as a test function in the weak formulation \eqref{weak:formulNS} of the Navier-Stokes equations between time $\underline{\tau}$ and $\tau$, for some $\underline{\tau} \leq \tau$ (see e.g. \cite{Masuda} for a discussion about this available procedure).
%Note that this is possible because $\varphi \in D(A) \cap \Ld^3(\R^3_+)$ and that we can extend the definition of \eqref{weak:formulNS} to this class of test functions, on any interval $[\underline{\tau}, \tau]$ because of the weak continuity of $u$ in $\Ld^2_{\mathrm{div}}(\R^3_+)$ (see \cite{Masuda} for a more complete discussion about this fact: here this is a straightforward consequence of the density of $\mathscr{D}_{\mathrm{div}}(\R^3_+)$ in $\H^1_{0, \mathrm{div}}(\R^3_+) \cap \Ld^3(\R^3_+)=\H^1_{0, \mathrm{div}}(\R^3_+)$). 
This yields
\begin{align}\label{testGoodtronqfreq}
\int_{\underline{\tau}}^{\tau} \Big\langle (u\cdot \nabla )u(\sigma), \varphi(\sigma) \Big\rangle \, \mathrm{d}\sigma=\int_{\underline{\tau}}^{\tau} \Big\langle F(\sigma), \varphi(\sigma) \Big\rangle \, \mathrm{d}\sigma-\Big\langle u(\tau), E_{\leq g^2(\tau)} \psi \Big\rangle + \Big\langle u(\underline{\tau}), e^{-(\tau-\underline{\tau})A} E_{\leq g^2(\tau)}\psi \Big\rangle.
\end{align}
%Indeed, note that two terms have cancelled because
%\begin{align*}
%\Big\langle u(\sigma),\partial_\sigma \varphi(\sigma) \Big\rangle=\Big\langle u(\sigma),A \varphi(\sigma)\Big\rangle= \Big\langle \nabla u(\sigma),\nabla \varphi(\sigma) \Big\rangle.
%\end{align*}
We now take $\underline{\tau}=0$ and by observing that $(u\cdot \nabla )u=\mathrm{div} \,(u \otimes u)$, we have
\begin{align*}
\int_{0}^{\tau} \Big\langle (u\cdot \nabla )u(\sigma), e^{-(\tau-\sigma)A} E_{\leq g^2(\tau)}\psi \Big\rangle \, \mathrm{d}\sigma = -\int_{0}^{\tau} \Big\langle u(\sigma), u(\sigma) \cdot \nabla  E_{\leq g^2(\tau)} e^{-(\tau-\sigma)A} \psi \Big\rangle \, \mathrm{d}\sigma,
\end{align*}
therefore \eqref{testGoodtronqfreq} yields
\begin{align*}
\left\vert \Big\langle u(\tau), E_{\leq g^2(\tau)} \psi \Big\rangle \right\vert &\leq \left\vert  \Big\langle u_0, e^{-\tau A} E_{\leq g^2(\tau)}\psi \Big\rangle \right\vert +  \int_{0}^{\tau} \left\vert\Big\langle F(\sigma), E_{\leq g^2(\tau)} e^{-(\tau-\sigma)A} \psi \Big\rangle \right\vert \, \mathrm{d}\sigma \\
& \quad +  \int_{0}^{\tau} \left\vert \Big\langle u(\sigma), u(\sigma) \cdot \nabla  E_{\leq g^2(\tau)} e^{-(\tau-\sigma)A} \psi \Big\rangle \right\vert \, \mathrm{d}\sigma \\
& :=(\mathrm{I}) + (\mathrm{II})+(\mathrm{III}).
\end{align*}
%For $(\mathrm{I})$, we merely write
%\begin{align*}
%(\mathrm{I}) \lesssim \Vert \mathbf{1}_{[0,g^2(\tau)]}(A) e^{-\tau A} u_0 \Vert_{\Ld^2(\R^3_+)} \Vert \psi \Vert_{\Ld^2(\R^3_+)} \leq \Vert e^{-\tau A} u_0\Vert_{\Ld^2(\R^3_+)} \Vert \psi \Vert_{\Ld^2(\R^3_+)},
%\end{align*}
%where $\lesssim$ hides a universal constant. For $(\mathrm{II})$, we use Lemma \ref{decayStokes:LpLq} and get
%\begin{align*}
%(\mathrm{II}) & \leq \int_{0}^{\tau} \Vert F(\sigma) \Vert_{\Ld^2(\R^3_+)} \Vert E_{\leq g^2(\tau)} e^{-(\tau-\sigma)A} \psi \Vert_{\Ld^2(\R^3_+)} \, \mathrm{d}\sigma  \lesssim   \Vert \psi \Vert_{\Ld^2(\R^3_+)} \int_{0}^{\tau} \Vert F(\sigma) \Vert_{\Ld^2(\R^3_+)} \, \mathrm{d}\sigma.
%\end{align*}
For $(\mathrm{I})$ and $(\mathrm{II})$, we use Lemma \ref{decayStokes:LpLq} and get
\begin{align*}
(\mathrm{I}) &\lesssim \Vert \mathbf{1}_{[0,g^2(\tau)]}(A) e^{-\tau A} u_0 \Vert_{\Ld^2(\R^3_+)} \Vert \psi \Vert_{\Ld^2(\R^3_+)} \leq \Vert e^{-\tau A} u_0\Vert_{\Ld^2(\R^3_+)} \Vert \psi \Vert_{\Ld^2(\R^3_+)}, \\
(\mathrm{II}) & \leq \int_{0}^{\tau} \Vert F(\sigma) \Vert_{\Ld^2(\R^3_+)} \Vert E_{\leq g^2(\tau)} e^{-(\tau-\sigma)A} \psi \Vert_{\Ld^2(\R^3_+)} \, \mathrm{d}\sigma  \lesssim   \Vert \psi \Vert_{\Ld^2(\R^3_+)} \int_{0}^{\tau} \Vert F(\sigma) \Vert_{\Ld^2(\R^3_+)} \, \mathrm{d}\sigma,
\end{align*}
where $\lesssim$ refers to a universal constant. For $(\mathrm{III})$, we rely on the following lemma, which can be found in the proof of \cite[Lemma 4.3]{BM}.
\begin{lem}
For any $w \in \H^1_{\mathrm{div}}(\R^3_+)$, $z \in \Ld^2_{\mathrm{div}}(\R^3_+)$ and $\gamma : \R^+ \rightarrow \R^+$, we have for all $\tau \geq 0$
\begin{align*}
\Big\langle w, w \cdot \nabla  E_{\leq \gamma(\tau)} z \Big\rangle  \leq C \gamma^{\frac{5}{4}}(\tau) \Vert w \Vert_{\Ld^2(\R^3_+)}^2 \Vert z \Vert_{\Ld^2(\R^3_+)},
\end{align*}
where $\langle \cdot , \cdot \rangle$ stands for the inner product on $\Ld^2(\R^3_+)$ and where $C>0$ is a universal constant.
\end{lem}
We apply this Lemma with $w=u(\sigma)$, $z=e^{-(\tau-\sigma)A} \psi$ and $\gamma=g^2$ to obtain
\begin{align*}
(\mathrm{III}) \lesssim \Vert \psi \Vert_{\Ld^2(\R^3_+)} g^{5/2}(\tau)\int_0^t \Vert u(\sigma) \Vert_{\Ld^2(\R^3_+)}^2 \, \mathrm{d}\sigma.
\end{align*}
Since all the previous estimates are valid for any $\psi \in \mathscr{D}_{\mathrm{div}}(\R^3_+)$, we infer that
\begin{align*}
\int_0^{g^2(\tau)}  \,  \Vert \mathrm{d}E_{\lambda}u(\tau)\Vert_{\Ld^2(\R^3_+)}^2  \leq  C\Vert e^{-\tau A} u_0\Vert_{\Ld^2(\R^3_+)}^2 + C \left( \int_{0}^{\tau} \Vert F(r) \Vert_{\Ld^2(\R^3_+)} \, \mathrm{d}r \right)^2 +  C g(\tau)^{5} \left( \int_{0}^{\tau} \Vert u(r) \Vert_{\Ld^2(\R^3_+)}^2 \, \mathrm{d}r \right)^2,
\end{align*}
where $C>0$ stands for a universal constant.
In view of \eqref{key}, we end up with the following differential inequality: for almost all $s \geq 0$ and for all $t \geq s$
\begin{align*}
\Vert u(t) \Vert_{\Ld^2(\R^3_+)}^2 + \int_s^t g^2(\tau) \Vert u(\tau) \Vert_{\Ld^2(\R^3_+)}^2 \, \mathrm{d}\tau  &\leq \Vert u(s) \Vert_{\Ld^2(\R^3_+)}^2 +2\int_s^t C g(\tau)^{7} \left( \int_{0}^{\tau} \Vert u(r) \Vert_{\Ld^2(\R^3_+)}^2 \, \mathrm{d}r \right)^2 \, \mathrm{d}\tau \\
& \quad +2\int_s^t  C g^2(\tau) \Bigg[  \Vert e^{-\tau A} u_0 \Vert^2_{\Ld^2(\R^3_+)} + g^{-4}(\tau) \Vert F(\tau)\Vert_{\Ld^2(\R^3_+)}^2  \\
& \qquad \qquad \qquad \qquad + \left( \int_{0}^{\tau} \Vert F(r) \Vert_{\Ld^2(\R^3_+)} \, \mathrm{d}r \right)^2 \Bigg]\, \mathrm{d}\tau.
\end{align*}
%Owing to the Grönwall Lemma-like \ref{delay:gron} stated in Appendix \ref{Gronwall-appendix}, with 
We then rely on the following Grönwall-like inequality (see \cite{BM-goodGron}).
\begin{lem}\label{delay:gron}
Let $t>0$. Let $y \in \Ld^{\infty}(0,t)$ and $\beta \in \Ld^1(0,t)$ be nonnegative functions. Suppose that the following differential inequality holds: for almost all $0 \leq s<t$
\begin{align*}
y(t) + \int_s^t \mathrm{g}(\tau) y(\tau) \, \mathrm{d} \tau \leq y(s) + \int_s^t \beta(\tau) \, \mathrm{d} \tau,
\end{align*}
where $\tau \mapsto \mathrm{g}(\tau)$ is a smooth positive function on $[0,t]$.
Then we have for almost all $0\leq s<t$
 $$
 \exp\left( \int_{0}^t \mathrm{g} (\tau) \, \mathrm{d} \tau\right)  y(t) \leq \exp\left( \int_{0}^s \mathrm{g} (\tau) \, \mathrm{d} \tau\right) y(s) + \int_{s}^t \exp\left(  \int_{0}^{\tau} \mathrm{g} (r) \, \mathrm{d} r\right) \beta(\tau) \, \mathrm{d} \tau.
 $$
 \end{lem}
%We apply this Lemma with
%\begin{align*}
% y(t) &= \Vert u(t) \Vert_{\Ld^2(\R^3_+)}^2, \\
% \beta(\tau) &=  C g^2(\tau) \left[  \Vert e^{-\tau A} u_0 \Vert^2_{\Ld^2(\R^3_+)} + g^{-4}(\tau) \Vert F(\tau)\Vert_{\Ld^2(\R^3_+)}^2 + \left( \int_{0}^{\tau} \Vert F(r) \Vert_{\Ld^2(\R^3_+)} \, \mathrm{d}r \right)^2  \right]\\
% &\quad + C g^{7}(\tau)  \left( \int_0^\tau \| u(\sigma)\|_{\Ld^2(\R^3_+)}^2 \, \mathrm{d} \sigma\right)^2,
% \end{align*}
%and we finally obtain
Applying this lemma, we finally obtain
\begin{equation}
  \begin{aligned}
 \Vert u(t) \Vert_{\Ld^2(\R^3_+)}^2  \exp\left(\int_0^t g^2(\tau) \, \mathrm{d} \tau  \right) & \leq 
\Vert u_0 \Vert_{\Ld^2(\R^3_+)}^2 \\
& \quad+ C  \int_0^t g^2(\tau) \left[ \Vert e^{-\tau A} u_0 \Vert_{\Ld^2(\R^3_+)}^2 + g^{-4}(\tau) \Vert F(\tau)\Vert_{\Ld^2(\R^3_+)}^2 \right] \exp\left(\int_0^\tau g^2(r) \, \mathrm{d} r  \right) \mathrm{d} \tau \\
&\quad + C  \int_0^t g^2(\tau)\left( \int_0^\tau \| F(\sigma)\|_{\Ld^2(\R^3_+)} \, \mathrm{d} \sigma \right)^2 \exp\left(\int_0^\tau g^2(r) \, \mathrm{d} r  \right) \mathrm{d} \tau \\
&\quad + C  \int_0^t g^7(\tau)\left( \int_0^\tau \| u(\sigma)\|_{\Ld^2(\R^3_+)}^2 \, \mathrm{d} \sigma \right)^2 \exp\left(\int_0^\tau g^2(r) \, \mathrm{d} r  \right) \mathrm{d} \tau.
\end{aligned}
  \end{equation}
We now take
 \begin{align*}
g^2(t):=\dfrac{\alpha}{1+t}, \ \ \ \ \exp\left(\int_0^t g^2(\tau) \, \mathrm{d} \tau \right) =(1+t)^{\alpha},
 \end{align*}
where $\alpha>0$ and this means that for all $t \in [0,T]$
\begin{equation}\label{eq-k}
  \begin{aligned}
 \Vert u(t) \Vert_{\Ld^2(\R^3_+)}^2  (1+t)^{\alpha}& \lesssim
\Vert u_0 \Vert_{\Ld^2(\R^3_+)}^2 
+  \int_0^t (1+\tau)^{\alpha-1} \left[ \Vert e^{-\tau A} u_0 \Vert_{\Ld^2(\R^3_+)}^2 + (1+\tau)^2 \Vert F(\tau)\Vert_{\Ld^2(\R^3_+)}^2 \right]  \mathrm{d} \tau \\
& \quad +   \int_0^t (1+\tau)^{\alpha-1}\left( \int_0^\tau \| F(\sigma)\|_{\Ld^2(\R^3_+)} \, \mathrm{d} \sigma \right)^2  \mathrm{d} \tau \\
&\quad +   \int_0^t (1+\tau)^{\alpha-7/2}\left( \int_0^\tau \| u(\sigma)\|_{\Ld^2(\R^3_+)}^2 \, \mathrm{d} \sigma \right)^2  \mathrm{d} \tau,
\end{aligned}
  \end{equation}
where $\lesssim$ is independent of $t$.
Using Lemma \ref{decayStokes:LpLq}, we get for all $\tau >0$
%$\theta \in (0,1)$
%\begin{align*}
%\Vert e^{-\tau A} u_0 \Vert_{\Ld^2(\R^3_+)} \leq  C\Vert  u_0 \Vert_{\Ld^2(\R^3_+)}, \ \ \ \
%%\Vert e^{-\tau A} u_0 \Vert_{\Ld^2(\R^3_+)} &\leq  \dfrac{C}{\tau^{\frac{3}{4}}}\Vert  u_0 \Vert_{\Ld^{\frac{3}{3-\theta}}(\R^3_+)},
%\Vert e^{-\tau A} u_0 \Vert_{\Ld^2(\R^3_+)} \leq  \dfrac{C}{\tau^{\frac{3}{2}(1-\frac{1}{2})}}\Vert  u_0 \Vert_{\Ld^{1}(\R^3_+)},
%\end{align*}
%therefore for all $\tau >0$ 
% il suffit juste d'utiliser la première inégalité si $t_0<\tau<1$ et l'autre si $+1<\tau$
\begin{align*}
\Vert e^{-\tau A} u_0 \Vert_{\Ld^2(\R^3_+)}^2 \leq \dfrac{C}{(1+\tau)^{3/2}} \left(\Vert  u_0 \Vert_{\Ld^{2}(\R^3_+)}^2+\Vert  u_0 \Vert_{\Ld^{1}(\R^3_+)}^2 \right),
\end{align*}
where $C$ is independent of $\tau$ and $u$. Observe now that, according to the assumption \eqref{HYPdecay:BRINK2}, the contribution of the source term $(1+\tau)^2 \Vert F(\tau)\Vert_{\Ld^2(\R^3_+)}^2 + \left( \int_{0}^{\tau} \Vert F(r) \Vert_{\Ld^2(\R^3_+)} \, \mathrm{d}r \right)^2$ is as $(1+\tau)^{-3/2}$ so that for all $t>0$
\begin{align}\label{gU_0+F}
\begin{split}
\Vert u_0 \Vert_{\Ld^2(\R^3_+)}^2 
+  \int_0^t (1+\tau)^{\alpha-1}& \left[  \Vert e^{-\tau A} u_0 \Vert^2_{\Ld^2(\R^3_+)} + g^{-4}(\tau) \Vert F(\tau)\Vert_{\Ld^2(\R^3_+)}^2 + \left( \int_{0}^{\tau} \Vert F(r) \Vert_{\Ld^2(\R^3_+)} \, \mathrm{d}r \right)^2\right]  \mathrm{d} \tau \\ 
& \lesssim \Vert u_0 \Vert_{\Ld^2(\R^3_+)}^2 +   \int_{0}^t  \dfrac{1}{(1+\tau)^{1-\alpha+3/2}} \, \mathrm{d}\tau, \\
 & \lesssim \Vert u_0 \Vert_{\Ld^2(\R^3_+)}^2 +(1+t)^{\alpha-3/2},
\end{split}
\end{align}
provided that $\alpha \neq 3/2$, and where $\lesssim$ depends on $\Vert u_0 \Vert_{\Ld^2(\R^3_+)}^2 + \Vert u_0 \Vert_{\Ld^1(\R^3_+)}^2+ \mathrm{C}$, for the constant $\mathrm{C}>0$ appearing in \eqref{HYPdecay:BRINK2}.
%We set
%$$\mathscr{E}_1(0):=\E(0)+\Vert  u_0 \Vert_{\Ld^{1}(\R^3_+)}^2$$

Let us assume that on $[0,T]$, we have
\begin{align}
\label{hyp:beta1}\Vert u(t) \Vert_{\Ld^2(\R^3_+)}^2  & \lesssim  \dfrac{1}{(1+t)^{\beta}},
\end{align}
for some $\beta >0$, different from $1$. 
%If $\beta<1$, then 
%\begin{align*}
%\int_{0}^t (1+\tau)^{\alpha-7/2}\left( \int_{0}^\tau \| u(r)\|_{\Ld^2(\R^3_+)}^2 \, \mathrm{d} r \right)^2  \, \mathrm{d} \tau 
%& \leq \int_{0}^t (1+\tau)^{\alpha-7/2+2(1-\beta)}\, \mathrm{d} \tau,
%\end{align*}
%and hence, as we have a bound by a constant if $\beta >1$ (provided that $\alpha<5/2$), 
From this estimate, we deduce the following inequality
\begin{align}\label{g7}
\begin{split}
\int_{0}^t (1+\tau)^{\alpha-7/2}\left( \int_{0}^\tau \| u(r)\|_{\Ld^2(\R^3_+)}^2 \, \mathrm{d} r \right)^2 \, \mathrm{d} \tau   
 \lesssim   \left\{
    \begin{array}{ll}
          1 \ \ & \text{if}  \ \ \beta<1, \ \ \alpha-2 \beta-1/2<0, \\
        (1+t)^{\alpha-2\beta-1/2} \ \ &  \text{if} \ \ \beta<1, \ \  \alpha-2 \beta-1/2 \geq 0, \\
        1 \ \ &  \text{if} \ \ \beta>1.
    \end{array}
\right.
\end{split}
\end{align}
Now, we start with $\beta=0$: the a priori estimate \eqref{hyp:beta1} indeed holds with this choice of exponent because the energy inequality \eqref{ineq-energyNS} can be rewritten, together with \eqref{HYPdecay:BRINK2}, as
\begin{align*}
\Vert u(t) \Vert_{\Ld^2(\R^3_+)}^2  +2\int_0^t  \Vert \nabla u(s) \Vert_{\Ld^2(\R^3_+)}^{2}  \mathrm{d}s &\lesssim  \Vert u_0 \Vert_{\Ld^2(\R^3_+)}^2 + \int_0^t \Vert F(s) \Vert_{\Ld^2(\R^3_+)} \, \mathrm{d}s  \lesssim \Vert u_0 \Vert_{\Ld^2(\R^3_+)}^2 + \int_0^t (1+s)^{-7/4} \, \mathrm{d}s \\
& \lesssim \Vert u_0 \Vert_{\Ld^2(\R^3_+)}^2 +1.
\end{align*}
We then take $\alpha=1$ and use \eqref{eq-k}, \eqref{gU_0+F} and \eqref{g7} to obtain
\begin{align*}
(1+t) \Vert u(t) \Vert_{\Ld^2(\R^3_+)}^2 \lesssim 1+(1+t)^{-1/2} +(1+t)^{1/2} ,
\end{align*}
therefore
\begin{align*}
 \Vert u(t) \Vert_{\Ld^2(\R^3_+)}^2 \lesssim (1+t)^{-1/2}.
\end{align*}
This means that the a priori estimate \eqref{hyp:beta1} now holds with $\beta=1/2$. We then take $\alpha=2$ and combining again \eqref{eq-k}, \eqref{gU_0+F} and \eqref{g7} yields
\begin{align*}
 \Vert u(t) \Vert_{\Ld^2(\R^3_+)}^2 &\lesssim (1+t)^{-2}+(1+t)^{-3/2} \lesssim (1+t)^{-3/2}.
\end{align*}
Note that in view of \eqref{gU_0+F} and \eqref{g7},  we cannot improve the exponent in the previous estimate. The proof of Theorem \ref{cond:decay} is therefore complete.
\end{proof}

\subsection{Parabolic regularization for the Navier-Stokes system on $\R^3_+$}\label{AnnexeParabNS}

\begin{thm}\label{RegParabNS}
For all $T>0$, there exists a universal constant $\mathrm{C}_{\star}>0$ such that the following holds. Consider $u_0 \in \H^1_{\mathrm{div}}(\R^3_+)$ and $F \in \Ld^2_{\mathrm{loc}}(\R^+;\Ld^2(\R^3_+))$ and $T>0$ such that
\begin{align}\label{Hyp:AnnexNSreg}
\Vert u_0 \Vert_{\H^1(\R^3_+)}^2 +  \int_0^T \left[ \Vert F(s) \Vert_{\Ld^2(\R^3_+)}^2 + \Vert F(s) \Vert_{\Ld^2(\R^3_+)} \right] \, \mathrm{d}s \leq \mathrm{C}_{\star}.
\end{align}
Then, there exists on $[0,T]$ a unique Leray solution to the Navier-Stokes system 
\begin{equation*}
\left\{
      \begin{aligned}
        \partial_t u +(u\cdot \nabla_x )u-\Delta_x u + \nabla_x p&=F,\\
        \mathrm{div}_x \, u&=0, \\
        u_{\mid x_3=0}&=0, \\
        u_{\mid t=0}&=u_0,
      \end{aligned}
    \right.
\end{equation*}
with initial data $u_0$ and source $F$. This solution $u$ belongs to $\Ld^{\infty}(0,T;\H^1_{\mathrm{div}}(\R^3_+)) \cap \Ld^{2}(0,T;\H^2(\R^3_+))$ and satisfies for almost every $t \in [0,T]$
\begin{align}\label{Ineg:AnnexNSreg}
\Vert \nabla u(t) \Vert_{\Ld^2(\R^3_+)}^2 + \dfrac{1}{2}\int_0^t \Vert  \mathrm{D}^2 u(s) \Vert_{\Ld^2(\R^3_+)}^2 \, \mathrm{d}s \leq  \widetilde{\mathrm{C}} \left( \Vert u_0 \Vert_{\H^1(\R^3_+)}^2 +  \Vert F \Vert_{\Ld^2(0,T;\Ld^2(\R^3_+))}^2 \right),
\end{align}
for some universal constant $\widetilde{\mathrm{C}}>0$.
\end{thm}
\begin{proof}
We argue in two different steps. First, if such a Leray solution to the Navier-Stokes equations exists, it satisfies in particular $u \in \Ld^{4}(0,T;\H^1_{\mathrm{div}}(\R^3_+))$, which is a well-known case of weak-strong uniqueness for this system (see e.g. \cite[Theorem 3.3]{NS-chemin}, which holds for general domains). Thus, this solution will be equal to any Leray solution of the Navier-Stokes system with source $F$ and initial value $u_0$.
%\textcolor{red}{In a second step, we have to prove that such a solution indeed exists on $[0,T]$: to do so, we rely on a standard approximation procedure by regularising the data $F$ and $u_0$ and by considering a Galerkine approximation $(u_N)_N$ of the corresponding Navier-Stokes system. We then obtain additionnal estimates for the sequence $(u_N)_N$ on $[0,T]$ which are available thanks to \ref{Hyp:AnnexNSreg}. Combining these parabolic estimates with the classical energy estimates for the Leray solution of the Navier-Stokes system on $[0,T]$, we use a classical compactness argument to get the existence of a Leray solution in $\Ld^{\infty}(0,T;\H^1_{\mathrm{div}}(\R^3_+)) \cap \Ld^{2}(0,T;\H^2(\R^3_+))$ and which satisfies the estimate (\ref{Ineg:AnnexNSreg}) on $[0,T]$.}

In a second step, we have to prove that such a solution indeed exists on $[0,T]$: to do so, we rely on a standard approximation procedure by regularizing the data $F$ and $u_0$ and by smoothing the convection term in the momentum equation (see e.g. \cite{MS} where the Yosida operator is used, as in Section \ref{Section:PBCauchy}). We essentially obtain a sequence $(u_N)_N$ of solutions to regularized Navier-Stokes systems. This sequence will classically be bounded in $\Ld^{\infty} \Ld^2 \cap \Ld^{2} \H^1$ and will satisfy additional parabolic estimates coming from \eqref{Hyp:AnnexNSreg}. A compactness argument allows to get the existence of a Leray solution in $\Ld^{\infty}(0,T;\H^1_{\mathrm{div}}(\R^3_+)) \cap \Ld^{2}(0,T;\H^2(\R^3_+))$ and which satisfies the estimate (\ref{Ineg:AnnexNSreg}) on $[0,T]$.

Hence in the following, we deal with a smooth solution $u$ with smooth data $u_0$ and we are looking for parabolic estimates for the velocity field $u$. The following strategy is based on an idea which is well-known in the context of the inhomogeneous incompressible Navier-Stokes equations (see \cite{Lady, DanchinINS}): we first apply the Leray projection $\mathbb{P}$ to the Navier-Stokes equations and then multiply by $\partial_t u$ to obtain, after integration by parts
\begin{align*}
\dfrac{1}{2}\dfrac{\mathrm{d}}{\mathrm{d}t}\Vert \nabla u \Vert_{\Ld^2(\R^3_+)}^2 + \Vert  \partial_t u \Vert_{\Ld^2(\R^3_+)}^2= -\langle (u \cdot \nabla )u, \partial_t u \rangle_{\Ld^2(\R^3_+)}+\langle F, \partial_t u \rangle_{\Ld^2(\R^3_+)},
\end{align*}
where we have dropped the time variable. We then use the Young inequality twice in order to absorb $ \Vert  \partial_t u \Vert_{\Ld^2(\R^3_+)}^2 $ in the l.h.s, that is on $[0,T]$
\begin{align}\label{mult:D_tu}
\dfrac{\mathrm{d}}{\mathrm{d}t}\Vert \nabla u \Vert_{\Ld^2(\R^3_+)}^2 + \Vert  \partial_t u \Vert_{\Ld^2(\R^3_+)}^2 \leq 2 \int_{\R^3_+} \vert (u \cdot \nabla )u \vert^2 \, \mathrm{d}x + 2 \int_{\R^3_+} \vert F\vert^2 \, \mathrm{d}x.
\end{align}
In order to estimate $\mathrm{D}^2 u$, we rewrite the Navier-Stokes system satisfied by $(u,p)$ as the following stationary Stokes system in the half-space $\R^3_+$ with source term  $F-(u \cdot \nabla )u-\partial_t u$
\begin{align*}
 \left\{
    \begin{array}{ll}
         -\Delta u + \nabla p= F-(u \cdot \nabla )u-\partial_t u,\\[2mm]
         \mathrm{div} \, u=0,\\[2mm]
        u_{\mid x_3=0}=0. 
    \end{array}
\right.
\end{align*}
According to \cite[Theorem IV.3.2]{Galdi}, there exists a universal constant $C>0$ such that on $[0,T]$
\begin{align*}
\Vert \mathrm{D}^2 u \Vert_{\Ld^2(\R^3_+)} \leq C \Vert F-(u \cdot \nabla )u-\partial_t u \Vert_{\Ld^2(\R^3_+)},
\end{align*}
and thus, denoting by $C>0$ another universal constant, we have on $[0,T]$
\begin{align*}
\Vert \mathrm{D}^2 u \Vert_{\Ld^2(\R^3_+)}^2 \leq C \left( \Vert F \Vert_{\Ld^2(\R^3_+)}^2 +  \Vert (u \cdot \nabla )u \Vert_{\Ld^2(\R^3_+)}^2 + \frac{1}{2}\Vert \partial_t u \Vert_{\Ld^2(\R^3_+)}^2 \right).
\end{align*}
Coming back to \eqref{mult:D_tu}, we get
\begin{align}\label{mult:D_tu2}
\dfrac{\mathrm{d}}{\mathrm{d}t}\Vert \nabla u \Vert_{\Ld^2(\R^3_+)}^2 +\frac{1}{2} \Vert  \partial_t u \Vert_{\Ld^2(\R^3_+)}^2 + \frac{1}{C}  \Vert \mathrm{D}^2 u \Vert_{\Ld^2(\R^3_+)}^2  \leq 3 \int_{\R^3_+} \vert (u \cdot \nabla )u \vert^2 \, \mathrm{d}x + 3 \int_{\R^3_+} \vert F\vert^2 \, \mathrm{d}x.
\end{align}
We now provide a bound for the convective term on [0,T] by using consecutively the Cauchy-Schwarz inequality, interpolation inequality for Lebesgue spaces and the Gagliardo-Nirenberg-Sobolev inequality (see Theorem \ref{gagliardo-nirenberg}) to write
\begin{align}\label{majConv1}
\begin{split}
\int_{\R^3_+} \vert (u \cdot \nabla )u \vert^2 \, \mathrm{d}x & \leq \Vert u \Vert_{\Ld^4(\R^3_+)}^2 \Vert \nabla u \Vert_{\Ld^4(\R^3_+)}^2 \\
& \leq \Vert u \Vert_{\Ld^2(\R^3_+)}^{1/2} \Vert u \Vert_{\Ld^6(\R^3_+)}^{3/2} \Vert \nabla u \Vert_{\Ld^2(\R^3_+)}^{1/2} \Vert \nabla u \Vert_{\Ld^6(\R^3_+)}^{3/2} \\
& \lesssim \Vert u \Vert_{\Ld^2(\R^3_+)}^{1/2} \Vert \nabla u \Vert_{\Ld^2(\R^3_+)}^{3/2} \Vert \nabla u \Vert_{\Ld^2(\R^3_+)}^{1/2} \Vert \mathrm{D}^2 u \Vert_{\Ld^2(\R^3_+)}^{3/2}.
\end{split}
\end{align}
Thanks to Young inequality, we deduce the existence of a universal constant $\overline{C}>0$ such that
%therefore the Young inequality ($4,4/3)$ yields the existence of a universal constant $\overline{C}>0$ such that
\begin{align}\label{majConv2}
\begin{split}
3 \Vert (u \cdot \nabla )u  \Vert_{\Ld^2(\R^3_+)}^{2} \leq \overline{C} \Vert u \Vert_{\Ld^2(\R^3_+)}^{2} \Vert \nabla u \Vert_{\Ld^2(\R^3_+)}^{8} + \frac{1}{2C} \Vert \mathrm{D}^2 u \Vert_{\Ld^2(\R^3_+)}^{2}.
\end{split}
\end{align}
We integrate \eqref{mult:D_tu2} between $0$ and $t \in [0,T]$ and get
\begin{multline*}
\Vert \nabla u(t) \Vert_{\Ld^2(\R^3_+)}^2 +\frac{1}{2C} \int_0^t \Vert \mathrm{D}^2 u(s) \Vert_{\Ld^2(\R^3_+)}^{2} \, \mathrm{d}s \\ 
\leq \Vert \nabla u_0 \Vert_{\Ld^2(\R^3_+)}^2 +3  \int_0^t \Vert F(s) \Vert_{\Ld^2(\R^3_+)}^{2} \, \mathrm{d}s + \overline{C}  \int_0^t \Vert u (s) \Vert_{\Ld^2(\R^3_+)}^{2} \Vert \nabla u(s) \Vert_{\Ld^2(\R^3_+)}^{8} \, \mathrm{d}s.
\end{multline*}
If we set
\begin{align}
\label{x(t)} x(t)&:=\Vert \nabla u(t) \Vert_{\Ld^2(\R^3_+)}^2 +\frac{1}{2C} \int_0^t \Vert \mathrm{D}^2 u(s) \Vert_{\Ld^2(\R^3_+)}^{2} \, \mathrm{d}s,\\
\label{h(t)} h(t)&:=\Vert \nabla u_0 \Vert_{\Ld^2(\R^3_+)}^2 +3  \int_0^t \Vert F(s) \Vert_{\Ld^2(\R^3_+)}^{2} \, \mathrm{d}s,\\[2mm]
g(s)&:=\overline{C} \Vert u (s) \Vert_{\Ld^2(\R^3_+)}^{2} \Vert \nabla u(s) \Vert_{\Ld^2(\R^3_+)}^{2},
\end{align}
then the previous inequality can be written as
\begin{align*}
x(t) \leq h(t) + \int_0^t g(s) x(s)^3 \, \mathrm{d}s,
\end{align*}
for all $t \in [0,T]$. We are now in position to apply Bihari's lemma, that we recall now (see \cite{Dannan} for a proof).
\begin{lem}\label{Bihari}
Let $x,g$ and $h$ three positive continuous functions on $\R^+$, such that $h$ is increasing . Let $w$ a continuous submultiplicative and nondecreasing function on $\R^+$ such that $w(u)>0$ for all $u>0$. Suppose that for all $t \geq 0$
\begin{align}
x(t) \leq h(t) + \int_0^t g(s) w(x(s)) \, \mathrm{d}s.
\end{align}
If the function
\begin{align*}
W(u):= \int_1^u \frac{\mathrm{d}s}{w(s)}, \ \ u \in (0,a) \ \ \text{for some} \ a>0,
\end{align*}
is a bijection on $(0,a)$, then the following inequality
\begin{align}\label{Bihari:ccl}
x(t) \leq h(t) W^{-1}\left(\int_0^t g(s) \dfrac{w(h(s))}{h(s)} \, \mathrm{d}s \right)
\end{align}
holds for all $t \geq 0$ such that $\int_0^t g(s) \frac{w(h(s))}{h(s)} \, \mathrm{d}s \in (0,a)$.
\end{lem}
We use this lemma with the function $w(s)=s^3$ associated to
\begin{align*}
W(u)=\int_1^u \dfrac{\mathrm{d}s}{s^3}=\frac{1}{2}\left( 1-\frac{1}{u^2} \right), \ \ u>0,
\end{align*}
and whose inverse, defined on $(0,1/2)$, is $W^{-1}(v)=(1-2v)^{-1/2}$. The inequality \eqref{Bihari:ccl} reads on $[0,T]$ as
\begin{align}
x(t) \leq h(t) \left( 1-2 \displaystyle \int_0^t g(s) h^2(s) \mathrm{d}s \right)^{-1/2},
\end{align}
provided that for all $t \in [0,T]$, we have $1-2  \int_0^t g(s) h^2(s) \mathrm{d}s >0$. Thus, this is enough to ensure for instance that
\begin{align}\label{cond:bihari:NSreg}
1-2  \int_0^T g(s) h^2(s) \mathrm{d}s >1/2,
\end{align}
is satisfied in order to get $x(t) \leq h(t)$ for all $t \in [0,T]$.  We observe that 
\begin{align*}
\int_0^T g(s) h^2(s) \mathrm{d}s \leq  \left[ \Vert \nabla u_0 \Vert_{\Ld^2(\R^3_+)}^2 +3  \int_0^T \Vert F(s) \Vert_{\Ld^2(\R^3_+)}^{2} \, \mathrm{d}s \right]^2 \int_0^T g(s) \mathrm{d}s,
\end{align*}
and
\begin{align*}
\int_0^T g(s) \mathrm{d}s \leq \overline{C} \Vert u \Vert_{\Ld^{\infty}(0,T;\Ld^2(\R^3_+))}^{2} \int_0^T  \Vert \nabla u(s) \Vert_{\Ld^2(\R^3_+)}^{2}  \mathrm{d}s.
\end{align*}
Furthermore, the following energy inequality for the Navier-Stokes system with source $F$ and initial data $u_0$ is satisfied by $u$
\begin{align*}
\Vert u \Vert_{\Ld^{\infty}(0,T;\Ld^2(\R^3_+))}^{2}  +2\int_0^T  \Vert \nabla u(s) \Vert_{\Ld^2(\R^3_+)}^{2}  \mathrm{d}s \leq \widetilde{C} \left[ \Vert u_0 \Vert_{\Ld^2(\R^3_+)}^2 + \left(\int_0^T \Vert F(s) \Vert_{\Ld^2(\R^3_+)} \, \mathrm{d}s \right)^2 \right],
\end{align*}
where $\widetilde{C}>0$ is a universal constant. The condition \eqref{cond:bihari:NSreg} can thus be satisfied provided that the quantity
\begin{align*}
 \left[ \Vert \nabla u_0 \Vert_{\Ld^2(\R^3_+)}^2 +  \int_0^T \Vert F(s) \Vert_{\Ld^2(\R^3_+)}^{2} \, \mathrm{d}s \right]^2  \left[ \Vert u_0 \Vert_{\Ld^2(\R^3_+)}^2 + \left(\int_0^T \Vert F(s) \Vert_{\Ld^2(\R^3_+)} \, \mathrm{d}s \right)^2 \right]^2
\end{align*}
is small enough. We observe that an assumption of the type of \eqref{Hyp:AnnexNSreg} can indeed provide such a smallness. Now, in view of the definitions \eqref{x(t)}-\eqref{h(t)}, the inequality $x(t) \leq h(t)$ for all $t \in [0,T]$ brings to the conclusion \eqref{Ineg:AnnexNSreg} and this completes the proof.
\end{proof}
\end{appendix}
\medskip

\section*{Acknowledgements}
\addcontentsline{toc}{section}{Acknowledgements}
I would like to acknowledge my PhD advisors Daniel Han-Kwan and Ayman Moussa for
their confidence, support and guidance, as well as their constructive criticism, scientific suggestions and careful rereadings during the preparation of this work. I would also like to thank Jean-Yves Chemin and Raphaël Danchin for discussions about Theorem \ref{RegParabNS}.
Finally, I am very grateful to Richard Höfer and to the anonymous referee for their suggestions and remarks, as well as for pointing out several inaccuracies in first versions of this paper.

%\nocite*{}
\bibliographystyle{abbrv}

\addcontentsline{toc}{section}{References}

\bibliography{biblio}

\end{document}